\documentclass[3p, final, hidelinks]{elsarticle}
\usepackage{framed,multirow}

\usepackage{lmodern}
\usepackage{microtype}
\usepackage[dvipsnames]{xcolor}
\usepackage{tikz}
\usepackage{hyperref}

\usepackage{printlen}
\microtypesetup{activate={true,nocompatibility},
                final,
                tracking=true,
                kerning=true,
                spacing=true,
                factor=1000,
                stretch=20,
                shrink=20}
\microtypecontext{spacing=nonfrench}

\usepackage{mathtools}
\usepackage{bm}
\usepackage{tabularx}
\usepackage{booktabs}

\graphicspath{{Figures_arxiv/}}

\usepackage{amsmath}
\newcommand{\pnpm}[2]{\ensuremath{\mathbb{P}_{#1}\mathbb{P}_{#2}}}
\usepackage{amssymb}
\usepackage{latexsym}
\usepackage{graphicx}  
\usepackage{subfig}
\usepackage{amsthm}

\newcommand{\up}[1]{\ensuremath{\mathrm{#1}}}
\DeclarePairedDelimiter{\norm}{\lVert}{\rVert}

\DeclareMathOperator{\erfc}{erfc}
\DeclareMathOperator{\trace}{tr}
\DeclareMathOperator{\de}{d\!}
\renewcommand{\vec}[1]{\bm{\mathrm{#1}}}
\newcommand{\transpose}[1]{\ensuremath{{#1}^\mathsf{T}}}
\newcommand{\keps}{k_{\epsilon}}
\newcommand{\rs}{r_\ast}
\newcommand{\rx}{R_x}
\newcommand{\ry}{R_y}
\newcommand{\aderdg}{ADER-DG\ }
\newcommand{\aderweno}{ADER-WENO\ }
\newcommand{\finvol}{Finite Volume\ }
\newcommand{\disgal}{Discontinuous Galerkin\ }
   \newcolumntype{C}{>{\centering\let\newline\\\arraybackslash}X}
   \newcolumntype{R}{>{\raggedleft\let\newline\\\arraybackslash}X}
   \newcolumntype{L}{>{\raggedright\let\newline\\\arraybackslash}X}
   \newcolumntype{Z}{>{\raggedright\let\newline\\\arraybackslash}X}
\newcommand{\od}[2]{\frac{\de{#1}}{\vphantom{l^l}\de{#2}}}
\newcommand{\sinq}{\sin^2 \psi}
\newcommand{\cosq}{\cos^2 \psi}
\newcommand{\babs}{\norm{\vec{b}}}
\newcommand{\divergence}{\nabla\cdot}

\newcommand{\sch}{Schmidmayer \etal\  \cite{Schmidmayer2017}\ }
\newcommand{\vphantomi}{}
\newcommand{\ev}[2]{\left. #1\right|_{#2}}
\newcommand{\pd}[3][0pt]{\frac{\raisebox{#1}{$\partial#2$}}{\vphantom{l^l}\partial#3}}
\newcommand{\abs}[1]{\left|#1\right|}
\newcommand{\uvec}[1]{\bm{\mathrm{\hat{#1}}}}
\newcommand{\interior}[1]{#1^{\circ}}

\newcommand{\intelemt}{\int \limits_{t_n}^{t_{n+1}}}
\newcommand{\intelemomega}{\int \limits_{\Omega_{ijk}}}

\newcommand{\intelemome}{\int \limits_{\Omega_{ijk}}}
\newcommand{\intelemiome}{\int \limits_{\interior{\Omega_{ijk}}}}
\newcommand{\intelemdome}{\int \limits_{\partial\Omega_{ijk}}}
\newcommand{\qhm}{\vec{q}_h^\textsc{l}}
\newcommand{\qhp}{\vec{q}_h^\textsc{r}}
\newcommand{\fhm}{\vec{f}_h^\textsc{l}}
\newcommand{\fhp}{\vec{f}_h^\textsc{r}}
\newcommand{\qhl}{\vec{q}_{\textsc{l}}}
\newcommand{\qhr}{\vec{q}_{\textsc{r}}}
\newcommand{\shlll}{S_{\textsc{l}}}
\newcommand{\shllr}{S_{\textsc{r}}}

\newcommand{\psix}{{\psi_{m_x}}}
\newcommand{\psiy}{{\psi_{m_y}}}
\newcommand{\psiz}{{\psi_{m_z}}}
\newcommand{\psit}{{\psi_{m_t}}}

\newcommand{\dimi}{k_{\rho_1}}
\newcommand{\dimii}{k_{\rho_2}}
\newcommand{\dimviii}{k_{b_1}}
\newcommand{\assq}{a_\sigma^2}
\newcommand{\asq}{a^2}
\newcommand{\kfiv}{k_4}
\newcommand{\kfv}{k_5}
\newcommand{\kfvi}{k_6}
\newcommand{\kfvii}{k_7}
\newcommand{\kfviii}{k_8}
\newcommand{\kfix}{k_9}
\newcommand{\kfx}{k_{10}}
\newcommand{\kfxi}{c_h^2\,\kfx}
\newcommand{\kfxii}{\assq\,\kfx}

\newcommand{\khi}{h_1}
\newcommand{\khii}{h_2}
\newcommand{\khiii}{h_3}
\newcommand{\khiv}{h_4}
\newcommand{\khv}{h_5}
\newcommand{\khvi}{h_6}
\newcommand{\figurescalefactor}{0.94}
\newcommand{\figurescalefactorb}{0.915}
\newcommand{\figurescalefactorc}{0.94}
\newcommand{\etal}{\textit{et al}}
\newcommand{\aposteriori}{\textit{a posteriori\ }}
\newcommand{\halb}{\frac{1}{2}}

\journal{Journal of Computational Physics}

\begin{document}

\begin{frontmatter}

\title{\Large
\textbf{\textls[-4]{High order ADER schemes and GLM curl cleaning for a first 
order hyperbolic formulation of compressible flow with surface tension}}}

\author[1]{Simone Chiocchetti$^{*}$}

\ead{simone.chiocchetti@unitn.it}
\cortext[cor1]{Corresponding author}

\author[1]{Ilya Peshkov}

\author[2]{Sergey Gavrilyuk}

\author[1]{Michael Dumbser}

\address[1]{{Laboratory of Applied Mathematics, Department of Civil, Environmental and Mechanical Engineering}, 
    {University of Trento}, {Via Mesiano 77, 38123 Trento, {Italy}}}

\address[2]{{IUSTI UMR 7343-CNRS, D\'epartement de M\'ecanique}, 
    {Aix-Marseille Universit\'{e}}, {5 rue Enrico Fermi, 13453 Marseille, {France}}}

\begin{abstract}
In this work, we introduce two novel reformulations of the weakly hyperbolic model 
for two-phase
flow with surface tension, recently forwarded by Schmidmayer \etal. 
In the model, the tracking of phase
boundaries is achieved by using a \emph{vector interface field}, rather than a scalar tracer, so
that the surface-force stress tensor can be expressed directly as an algebraic function of the state variables,
without requiring the computation of gradients of the scalar tracer. An interesting and 
important feature of
the model is that this interface field obeys a \emph{curl involution constraint}, that is, the
vector field is required to be curl-free at all times. 

The proposed modifications are intended to restore the strong hyperbolicity of the model, and are closely
related to divergence-preserving numerical approaches developed in the field of 
numerical magnetohydrodynamics (MHD).
The first strategy is based on the theory of Symmetric Hyperbolic and 
Thermodynamically Compatible (SHTC) systems forwarded by Godunov in the 60s and 70s
and yields a modified system of governing equations which includes some symmetrisation terms, in 
analogy to the approach adopted later by Powell \etal\  in the 90s
for the ideal MHD equations. The second technique is an extension of the hyperbolic Generalized 
Lagrangian Multiplier (GLM) divergence cleaning approach, 
forwarded by Munz \etal\  in applications to the Maxwell and MHD equations.

We solve the resulting nonconservative hyperbolic PDE systems with high order ADER \disgal (DG) methods
with \aposteriori \finvol subcell limiting 
and carry out
a set of numerical tests concerning flows dominated by surface tension as well as shock-driven
flows. We also provide a new exact solution to the equations, show convergence of the schemes for 
orders of accuracy up to ten in space and time, and investigate the role of hyperbolicity and of 
curl constraints in the long-term stability of the computations.
\end{abstract}


\end{frontmatter}

\section{Introduction}\label{sec:introduction}

The necessity to develop mathematical models and numerical methods for explicit treatment of
liquid-gas or liquid-liquid interfaces and accounting for the surface tension is well recognized in
computational fluid mechanics. In particular, it is required in many applications relying on
numerical simulations such as atomization \cite{Popinet2010}, boiling and cavitation
\cite{Saurel2014}, additive manufacturing~\cite{Khairallah2016,Francois2017}, multi-phase flows in
porous media \cite{Raeini2012} and microfluidics \cite{Worner2012}.

A great number of models and numerical methods have been proposed, and reviewing the vast literature
existing on the topic is outside the scope of this paper. Instead, we refer the reader to some 
recent review articles in
the field, \textit{e.g.} \cite{Saurel2018,Popinet2018}. Each model and method has its own pros and
cons, and usually the choice is conditioned by many factors. In particular, our preference is to
develop a first-order hyperbolic model for surface tension which is motivated by our intention to
further develop the unified first-order hyperbolic formulation for continuum mechanics forwarded in 
\cite{HPR2016,dumbser2016b,dumbser2017,SHTC-GENERIC-CMAT}. Hyperbolic PDE systems have several
attractive features for the numerical treatment. Namely, \textit{(i)} the initial value problem for
such systems is locally well posed; \textit{(ii)} a quite general theory of hyperbolic PDEs can be 
established, \textit{e.g.} \cite{Serre2007,Dafermos2005}, and hyperbolic systems are subject to a
unified numerical treatment (one may use one and the same numerical solver for various hyperbolic
models representing totally different physics). In particular, this suggests a possibility for the 
development of a general purpose code \cite{Dumbser2018a}; \textit{(iii)} first-order hyperbolic systems are
less sensitive to the quality of the computational grid~\cite{Nishikawa2014a}, have less
restrictive stability conditions on the time-step \cite{dumbser2016b}, and in principle, the same nominal 
order of accuracy is achievable for all quantities, including those representing gradients in 
the high-order PDE models.
In the last decade, models and computational approaches relying on hyperbolic equations have
been developed for viscous Newtonian \cite{HPR2016,dumbser2006b,Nishikawa1} and non-Newtonian
\cite{Jackson2018} flows, for non-linear elastic and elastoplastic solids
\cite{GodRom2003,Hyper-Hypo2018,HYP2016}, for heat conduction \cite{dumbser2016b,SHTC-GENERIC-CMAT},
for resistive electrodynamics \cite{dumbser2017}, for dispersive equations
\cite{Romenski2011,Mazaheri2016,Dhaouadi2018}, and for the  Einstein field equations of gravity 
\cite{ccz4,Dumbser2019}.

Concerning the treatment of the material interfaces, one may in general distinguish the following
computational approaches \cite{Saurel2018}. In Lagrangian schemes for compressible multi-material 
flows \cite{MaireMM2,ShashkovMultiMat3,ShashkovRemap3,ShashkovRemap5} the interface is
fully resolved and has thickness equal to zero. These algorithms are therefore also called sharp 
interface methods. In these approaches, 
the computational grid is aligned with the interface and deforms following the interface kinematics. 
These methods can be very 
accurate if the deformation of the interface is not too large, but they might be limited by mesh
distortion at large deformations. Topological changes (break up, coalescences) are very difficult to
implement in the Lagrangian framework, but they are at least in principle possible, see \textit{e.g.} 
\cite{Springel,HOArepo}. Particle-based methods such as Smoothed Particle  
Hydrodynamics (SPH), see \cite{GingoldMonaghan}, can be also attributed to this class of methods 
and, while they solve the
issues of mesh tangling by removing it completely, they lack the sharp description of material 
interfaces that characterise mesh-based Lagrangian methods and are more alike diffuse interface
schemes in this regards. Furthermore, many SPH schemes also lack even zeroth order consistency 
and require artificial viscosity for stabilization, see \cite{GingoldMonaghan,Monaghan1994,
Monaghan2005}, although this issue can be resolved using Riemann solver based SPH schemes \cite{Vila,SPH3D}. 

A different alternative strategy to model the interface dynamics consists in employing a method
based on an Eulerian description of the governing equations, in order to avoid difficulties related
to the mesh handling. In this class of methods, the interfaces are tracked implicitly via a new
state variable that is often called colour function. Examples include the Level Set method (LSM) 
\cite{levelset1,FedkiwEtAl1,sussman1994level,sethian2003level} 
and the Volume of Fluid (VoF) approach \cite{HirtNichols,scardovelli1999direct}. Eulerian 
methods are characterised by the fact that the interface (or rather,
its implicit representation via a scalar tracer) is actually distributed (diffused) onto several
grid cells. However, LSM and VoF are computational techniques that allow a subsequent 
reconstruction of a localised representation of the material boundaries (that is, with segments,
polynomials, or other curves) and do not deal with the thermodynamics of the mixing zone
(interface). Due to their explicit interface reconstruction, also level-set and VoF methods are 
considered \textit{sharp interface methods}, although the interface is not explicitly resolved by the 
computational mesh, as in Lagrangian schemes. 

Diffuse interface methods (DIM), instead, do \textit{not} use any explicit interface reconstruction 
technique, but they are capable, by specifying directly at the PDE level the governing principles 
for mixed states, to  represent complex non-equilibrium thermodynamics of the interface in many 
physical settings, \textit{e.g.} phase transition and cavitation, mass transport, heat conduction, 
fluid-structure interaction, see \textit{e.g.} \cite{Anderson1998,saurel1999,kapila2001,perigaud2005,
petitpas2009a,favrie2012,Saurel2018,Galenko2018,Gomez2019,Gomez2017a}.

A fourth class of methods is that of Arbitrary Lagrangian Eulerian (ALE) schemes, see 
\cite{LagrangeNC,Lagrange3D,LagrangeDG,HOArepo} and references therein, which combine the 
Lagrangian mesh motion with an implicit representation of material discontinuities. This type of
methods are particularly attractive as the mesh motion significantly reduces the smearing of contact
discontinuities and of the scalar tracer variables that are employed for the description of material
boundaries.

\paragraph{Purpose of the paper}
The schemes employed in this work fall into the class of Eulerian methods and the proposed mathematical
models follow a diffuse interface approach. In particular, the interface capturing between different
phases is achieved by introducing a new \emph{vector interface field} $\vec{b}$, which represents the
gradient of a colour function. This new vectorial quantity, \textit{i.e.} the gradient of the colour function, is 
then explicitly evolved in time as a new independent state variable of the system, rather than 
performing a reconstruction of the phase boundaries from the colour function. In other works dealing 
with surface tension effects, the geometrical information on the direction of normal vectors and 
curvature of the interface is commonly obtained via post-processing of the scalar colour function, see 
\cite{brackbill1992, noh1976slic, youngs1982time, sussmanHF, popinetHF, popinetzaleski1999, fusterpopinet2018allmach}. 

In a sense, this methodology takes the diffuse interface approach further in the direction of
decoupling the task of tracking material interfaces from a specific solution algorithm or code, in
that we not only account for separation between different materials by means of a colour function to
then recover geometrical information via a reconstruction technique, but we directly include the
concept of an interface, a diffuse one specifically, at the PDE level. This significantly
simplifies the use of general purpose high order methods to solve the resulting governing equations, 
as it allows writing the governing equations as a system of first order PDEs and
eliminates the need for interface reconstruction procedures or finite-difference-type approximations
in the computation of local curvature values and interface normal vectors. On the other hand, the
downside of any first-order reduction of a high-order PDE system is the presence of so-called
involution constraints in the former \cite{Gundlach04a,ccz4,Dumbser2019}, which in general are \emph{stationary}
differential equations that are satisfied by the governing PDE system \emph{for all times} if they
are satisfied by the initial data. In particular in this paper we shall deal with the curl-type
involution constraint on the vector field 
$\vec{b}$
\begin{equation*}
\nabla \times \vec{b} = \vec{0}.
\end{equation*}
Thus, from the numerical viewpoint, one may emphasise that the physical consistency of the numerical
solution might be completely lost if the involution constraint is violated. The development of
involution-constraint-preserving numerical methods is one of the key aspects in dealing with
first-order reductions of high-order PDE systems. 

The motivation for this paper is thus two-fold. First, we propose two separate ways of recovering 
a strongly hyperbolic formulation from the original weakly hyperbolic model for surface tension in compressible
two-phase flow advanced by Gavrilyuk \etal\  in \cite{Berry2008a} and further developed by \sch. Both
reformulations strongly rely on the curl-free constraint of the interface vector field $\vec{b}$.
Therefore, our second motivation for this paper was to systematically test the ability of the family of 
high-order ADER
\disgal (ADER-DG) \cite{dumbser2008b,dumbser2009,dumbser2014,zanotti2014} and ADER \finvol (ADER-FV) schemes 
\cite{titarev2002,toro2002,titarev2005,dumbser2007b,dumbser2008a} designed for general systems of first 
order balance laws to deal with hyperbolic PDE systems with curl involution constraints, 
and to compare stability, accuracy and physical consistency of the proposed hyperbolic approaches 
with the original weakly hyperbolic formulation.

The first of the two reformulations of the model of \sch is based on the theory of Symmetric
Hyperbolic and Thermodynamically Compatible (SHTC) systems \cite{God1961,Godunov1996,SHTC-GENERIC-CMAT},
and consists in modifying the momentum and energy equations by adding some \emph{symmetrising}
nonconservative terms which are linear combinations of the involution constraint and thus, are formally zero
at the continuum level. The approach is analogous to that used by Powell \etal\  in 
\cite{Powell1997,powell1999}
for the equations of ideal Magnetohydrodynamics and based on the ideas of Godunov
\cite{godunov1972}, thus in this work we will refer to these non-conservative symmetrising 
terms as Godunov--Powell-type nonconservative products.

The second modified model is again based on ideas developed in the context of numerical
Magnetohydrodynamics and specifically on the hyperbolic Generalized Lagrangian Multiplier (GLM) 
divergence-cleaning approach of Munz \etal\  \cite{MunzCleaning, Dedneretal}. Thus, analogously to the
GLM formulation for ideal MHD equations, which include an involution constraint on the divergence of
the magnetic field, the surface tension model in consideration is \textit{augmented} by a similar evolution 
equation for the additional curl-cleaning vector field. Such a GLM extension of the original model
of \sch also allows to fix the issue with its weak hyperbolicity and makes the augmented GLM system 
again strongly hyperbolic.

The paper is organised as follows. In Section\,\ref{sec:mathematicalmodel}, we recall the features
of the weakly hyperbolic \sch model and discuss the two modifications of the system which allow to recover
strong hyperbolicity. We then provide an exact
solution for cylindrically and spherically symmetric objects subject to surface tension
forces, that would mean droplets and water jets, with some considerations on the implications on the
diffuse interface treatment of such objects. We also explicitly compute the eigenvalues and a full
set of eigenvectors for both of the new models proposed in this paper. In Section\,\ref{sec:numericalmethod}, 
we provide the details of the high-order ADER-DG and ADER-FV methods employed in this paper.
Section\,\ref{sec:testproblems} presents a set of test problems allowing to validate the
new mathematical models and the numerical schemes proposed in this paper. Finally, the summary of 
the results as well as a discussion of further perspectives is given in Section\,\ref{sec:conclusions}. 



\section{Mathematical model}\label{sec:mathematicalmodel}
The original weakly hyperbolic two-phase, single velocity, single pressure model 
proposed in the paper of \sch can be written as follows:
\begin{subequations} 
\label{eq:gavrilyuk}
    \begin{align}
        &\partial_t\left(\alpha_1\,\rho_1\right) 
            + \divergence\left(\alpha_1\,\rho_1\,\vec{u}\right) = 0, \label{gavr.mass1}\\[1mm] %
        &\partial_t\left(\alpha_2\,\rho_2\right) 
            + \divergence\left(\alpha_2\,\rho_2\,\vec{u}\right) = 0, \label{gavr.mass2}\\[1mm] %
        %
        %
        &\partial_t\left(\rho\,\vec{u}\right) 
            + \divergence\left(\rho\,\vec{u}\otimes\vec{u} + p\,\vec{I} + \vec{\Omega}\right) = 
            \vec{0}, \label{gavr.mom}\\[1mm] 
        &\partial_t\left(\rho\,E\right) 
            + \divergence\left[\left(\rho\,E + p\right)\,\vec{u} + \vec{\Omega}\,\vec{u}\right] = 
            0, \label{gavr.energy}\\[1mm] 
        &\partial_t\left(\alpha_1\right) + \vec{u} \cdot \nabla \alpha_1  - K\,\divergence\vec{u} = 
        0,\label{gavr.alpha}\\[1mm] 
        &\partial_t\left(\vec{b}\right) + \left(\nabla\vec{b}\right)\,\vec{u} + 
        \transpose{(\nabla\vec{u})}\,\vec{b} = \vec{0},\qquad \nabla\times\vec{b} = 
        \vec{0},\label{gavr.interface}\\[1mm] 
        &\partial_t\left(c\right) + \vec{u}\cdot\nabla c = 0.
    \end{align}
\end{subequations}
The model consists of two mass conservation equations \eqref{gavr.mass1} and \eqref{gavr.mass2}, 
one for each of the two phases, 
a single (vector) equation \eqref{gavr.mom} for the conservation of mixture momentum 
$\rho\,\vec{u}$ and one 
scalar
equation \eqref{gavr.energy} for the conservation of the total energy of the mixture {${\rho\,E = \rho\,e +
\tfrac{1}{2}\,\rho\norm{\vec{u}}^2 + \sigma\,\norm{\vec{b}}}$}, which includes a surface energy
contribution $\sigma\,\norm{\vec{b}}$ to be added to the usual internal and kinetic energy terms. We
denote {${\vec{u} = \transpose{\left(u,\ v,\ w\right)}}$} and {${\vec{b} = \transpose{\left(b_1,\
b_2,\ b_3\right)}}$}, the former being the velocity field, the latter the \emph{interface field}
defined in the following paragraphs. 

The nonconservative equation \eqref{gavr.alpha} can be derived from the 
pressure equilibrium equation and the hypotheses about isentropic behaviour of each phase.
Because
the volume fractions $\alpha_1$ and $\alpha_2$ are subject to the constraint $\alpha_1 + \alpha_2 =
1$, one equation is sufficient for the description of both. The mixture density is defined as $\rho
= \rho_1\,\alpha_1 + \rho_2\,\alpha_2$, where $\rho_1$ and $\rho_2$ are the densities of the first
and the second phase, respectively. Neglecting surface tension effects, the first five equations of
\eqref{eq:gavrilyuk} are known as Kapila's \cite{kapila2001} five equation model, which is the 
stiff relaxation limit of the seven-equation Baer--Nunziato model \cite{baer1986}.
The system is closed if the specific energy $e$ is given as a function of the other variables.
In this work, we employ 
the so-called stiffened gas equation of state, in order to establish a biunivocal
relation between the pressure of each phase $p_1$ or $p_2$ 
and the corresponding internal energy density 
$\rho_1\,e_1$ or $\rho_2\,e_2$ as follows:  
\begin{equation}
    p_1 = \left(\gamma_1 - 1\right)\,\rho_1\,e_1 - \gamma_1\,\Pi_1, \qquad p_2 = 
        \left(\gamma_2 - 1\right)\,\rho_2\,e_2 - \gamma_2\,\Pi_2,
\end{equation}
with $\gamma_1$, $\gamma_2$, $\Pi_1$, $\Pi_2$ given parameters of the equation of state. Due to the
pressure equilibrium assumption $p_1 = p_2 = p$, the mixture equation of state then reads 
\begin{equation}
    p = \frac{
      \rho\,e\,\left(\gamma_1 - 1\right)\left(\gamma_2 - 1\right) - \alpha_1\,\gamma_1\,\Pi_1\,
          \left(\gamma_2 - 1\right) - \alpha_2\,\gamma_2\,\Pi_2\,\left(\gamma_1 - 1\right)
      }{
      \alpha_1\,\left(\gamma_2 - 1\right) + \alpha_2\,\left(\gamma_1 - 1\right)
      }.
\end{equation}
Furthermore, for this choice of closure relation, we have
\begin{equation} \label{eq:kapilaspeedofsound}
        K = \frac{\alpha_1\,\alpha_2\,\left(\rho_2\,a_2^2 - \rho_1\,a_1^2\right)}{\alpha_1\,
            \rho_2\,a_2^2 + \alpha_2\,\rho_1\,a_1^2},\qquad\textnormal{with}  \qquad 
        a_1 = \sqrt{\frac{\gamma_1\,\left(p + \Pi_1\right)}{\rho_1}}, \qquad
          a_2 = \sqrt{\frac{\gamma_2\,\left(p + \Pi_2\right)}{\rho_2}}.
\end{equation}

For the purpose of capturing the evolution of the interface geometry, a passively-advected scalar
quantity $c$, commonly termed colour function, is introduced; this quantity, similarly to the volume
fraction and mass fraction functions, ranges from zero to one and indicates, in a diffused sense,
the position of the interface.

Forces due to surface tension are taken into account by means of a tensor $\vec{\Omega}$ which can
be written in terms of the gradient of the colour function and of a constant surface tension
coefficient $\sigma$ as 
\begin{equation} \label{eq:surfacetension}
  \vec{\Omega} = \sigma\,\norm{\nabla c}\,\left(\frac{
      \nabla c \otimes \nabla c}{\norm{\nabla c}^2} - \vec{I}\right).
\end{equation}
The associated surface energy density is given by 
$\sigma\,\norm{\nabla c}$, meaning that when the colour function $c$ is 
a Heaviside-type function, the surface energy 
is a Dirac-delta-type function.
Such a conservative formulation \cite{gueyffier1999} of surface tension is well established and
essentially equivalent to the very popular Continuum Surface Force (CSF) approach of Brackbill \etal\ 
\cite{brackbill1992}, but since the tensor $\vec{\Omega}$ depends non-linearly on the derivatives of
the state variables, it is difficult to certify the well-posedness of the resulting initial-value
problem. Moreover, in order to compute surface tension forces, one would have to reconstruct the
colour function data and deduce the interface position, interface normal vectors and the local
interface curvature from this reconstructed information. In order to obtain a first order hyperbolic
formulation, a new \emph{interface field} 
$\vec{b} = \nabla c = \transpose{\left(b_1,\ b_2,\ b_3\right)}$
was introduced in \sch, together with a corresponding evolution equation
in the form 
\begin{equation}
  \partial_t\left(\vec{b}\right) + \left(\nabla\vec{b}\right)\,\vec{u} + 
      \transpose{(\nabla\vec{u})}\,\vec{b} = \vec{0}, 
\end{equation}
in which all components of the interface field should be treated as independent state variables. 
Since the new field $\vec{b}$ has been defined as the gradient of a scalar function, it 
must satisfy the constraint
\begin{equation} \label{curl.b}
\nabla \times \vec{b} = \vec{0}.
\end{equation}

This procedure, besides allowing to write the governing equations as a system of first order PDEs,
completely avoids the computation of local curvature values and interface normal vectors. The 
surface tension tensor can now be directly computed from the state variables as a nonlinear algebraic 
function  
\begin{equation} \label{eq:surfacetension2}
  \vec{\Omega} = \sigma\,\norm{\vec{b}}\,\left(\frac{\vec{b} \otimes 
      \vec{b}}{\norm{\vec{b}}^2} - \vec{I}\right).
\end{equation}
In turn, a new difficulty is introduced by the necessity of properly representing and transporting
the interface field $\vec{b}$ itself, which can be extremely challenging for a numerical scheme due
to the presence of Dirac-delta-like features\footnote{Recall that $\mathbf{b}=\nabla c$, hence if
$c$ approximates a step function, its gradient is an approximation of the delta distribution},
requiring very high spatial resolution and low numerical dissipation. Nonetheless, the resolution
requirements can be slightly relaxed by initializing the interface field from a smoothed colour
function, without observing strong modifications of the pressure jumps across interfaces, as can be
noted in Figure \ref{fig:droplet-exact-solution-double}, in which we show the exact solutions for
the pressure profiles inside a droplet having a geometrical representation with varying degree of
interface smoothing.
To clarify, in this work, contrary to the numerical approach adopted by Schmidmayer \etal\ 
in \cite{Schmidmayer2017}, the colour function $c$ itself is \emph{never} used for computing the
capillary stress tensor, nor has it any effect on any of the other fields, that is, the evolution
equation of the colour function is only coupled passively with the rest of the system and could be, in
principle, removed altogether from the computation. 

A major difficulty in the numerical discretisation of the evolution equation of $\mathbf{b}$ is 
obviously the curl involution constraint $\nabla \times \mathbf{b}=0$, which will be thoroughly treated  
and discussed in this paper.   

The last, but not the least, difficulty derives from the lack of hyperbolicity of the original \sch
model. In particular, it is shown in \cite{Schmidmayer2017} that model~\eqref{eq:gavrilyuk} is only
\textit{weakly hyperbolic}, that is all its eigenvalues are real, but one 
cannot find a full set of linearly independent eigenvectors. In
other words, system~\eqref{eq:gavrilyuk} is not diagonalizable. It is known that the initial-value
problem (Cauchy problem) is not necessarily well-posed in a conventional sense (\textit{i.e.} in
$C^\infty$) for weakly hyperbolic systems, \textit{e.g.} see \cite{Colombini2018, Colombini1982}.
The simplest example is the following. Consider a linearised system describing a pressureless gas 
\begin{equation}
\label{eq:whsystemexample}
\begin{aligned}
    &\partial_t \, \rho'  + u_0\,\partial_x \rho' + \rho_0\,\partial_x u' = 0,\\
    &\partial_t \, u' + u_0\,\partial_x u'  = 0,
\end{aligned}
\end{equation}
with the fluctuations $\rho'$, $u'$ and $\rho_0$, $u_0$ given constants. 
The system is weakly hyperbolic with eigenvalues $\lambda_1 = \lambda_2 = u_0$ and a single-parameter
eigenspace spanned by $\vec{R}_1 = \transpose{[1,\ 0]}$.
With suitable initial conditions, the solution of \eqref{eq:whsystemexample} is given by
\begin{equation}
    \begin{aligned}
        &u'(t,\ x) = f(x - u_0\,t),\\
        &\rho'(t,\ x) = t\,g(x - u_0\,t), \qquad \text{with} \qquad g(y) = 
        -\rho_0\,f^\prime(y),\quad y = x - u_0\,t.
    \end{aligned}
\end{equation}

Even for bounded initial data, the solution is not: it is growing linearly in time.
The lack of well-posedness of the Cauchy problem may cause instabilities and may require the design
of very specific numerical methods which help to stabilise the solution. More precisely: in order to
discretise the original weakly hyperbolic model \cite{Schmidmayer2017}, a special
structure-preserving numerical scheme would be needed, which is able to preserve the curl-free
condition \textit{exactly} at the discrete level for all times, similar to the exactly
divergence-free schemes developed in the last decades for the Maxwell and MHD equations, see \textit{e.g.}
\cite{Yee66,DeVore,BalsaraSpicer1999,BalsaraAMR,GardinerStone,BalsaraMultiDRS,ADERdivB,
balsarahlle2d,balsarahlle3d,MUSIC1,MUSIC2,BalsaraCED,HazraBalsara} and references therein. Much less
is known, instead, on exactly curl-preserving schemes. A rather general framework for the
construction of structure-preserving schemes (including curl-preserving methods) was developed by
Hyman and Shashkov in \cite{HymanShashkov1997} and Jeltsch and Torrilhon in
\cite{JeltschTorrilhon2006,Torrilhon2004}. Further work on mimetic and structure-preserving finite
difference schemes can be found \textit{e.g.} in \cite{Margolin2000,Lipnikov2014,Carney2013}. For families of
compatible finite element methods, the reader is referred to \cite{Nedelec1,Nedelec2,Cantarella,
Hiptmair,Monk,Arnold,Alonso2015}. However, to the very best of our knowledge, exactly
curl-preserving schemes for the PDE systems considered in this article have never been developed. 

Our intention here, however, is \textit{not}  
to develop specific (problem dependent) structure-preserving numerical techniques, but we rather prefer 
to use general purpose methods of the ADER-FV and ADER-DG family, which can be applied to very general  
hyperbolic systems with non-conservative products and (stiff) algebraic source terms. Therefore, 
one of the main goals of this paper is to modify the \sch model in order to obtain a 
\textit{strongly hyperbolic} version, with a full set of linearly independent eigenvectors. 
In the following Section we will discuss that at least two different ways of achieving such a goal exist. 

\subsection{Recovering hyperbolicity of the model with Godunov--Powell-type symmetrising
terms} Since the original weakly hyperbolic form of the model is not suitable for the solution with
explicit Godunov-type schemes, and motivated by the theory of Symmetric Hyperbolic and Thermodynamically
Compatible (SHTC) equations~\cite{God1961,Godunov1996,Rom2001,SHTC-GENERIC-CMAT}, we introduce some
formal modifications to system \eqref{eq:gavrilyuk} in such a way that the new system can be written
in the symmetric Godunov form and the eigenvector that was reported missing in the paper of
Schmidmayer \etal\  \cite{Schmidmayer2017} is recovered. Note that in \cite{Schmidmayer2017}, this
issue was circumvented by discretising the colour function equation and computing the interface
field as its gradient, instead of directly solving the vector equation for the interface field
$\vec{b}$. It is necessary to emphasise that the applied modifications are valid on smooth solutions, 
while the validity of the obtained non-conservative  hyperbolic model on discontinuities
requires further investigation.

The modifications are applied by introducing in the momentum and energy equations two
nonconservative terms that, at the continuum level at least, are identically null, thanks to the
curl constraint~\eqref{curl.b}, which is nothing else but Schwarz's rule of symmetry of second order 
derivatives 
\begin{equation} \label{eq:compatibilitytensor}
    \nabla \vec{b} = \transpose{\left(\nabla \vec{b}\right)}.
\end{equation}
In the paper of Schmidmayer \etal\  \cite{Schmidmayer2017} this property was used in the
hyperbolicity study, following Ndanou \etal\  \cite{ndanou2014}, to write the evolution
equation for the gradient of the colour function $\vec{b}$ in a Galilean-invariant form
\begin{equation}\label{b.Galilean}
    \partial_t \vec{b} + (\nabla \vec{b})\,\vec{u} + \transpose{(\nabla\vec{u})}\,\vec{b} = \vec{0},
\end{equation}
rather than a non-Galilean-invariant conservation form $
    \partial_t \vec{b} + \nabla\cdot\left[(\vec{b}\cdot\vec{u})\,\vec{I}\right] = \vec{0}. $

In order to conduct our mathematical and numerical study of system \eqref{eq:gavrilyuk}, we make use
of the same compatibility condition and rewrite the fully non-conservative equation \eqref{b.Galilean} 
in a semi-conservative form 
\begin{equation} 
\label{eq:nonconservativeb}
  \partial_t\vec{b} + \divergence \left[(\vec{u}\cdot\vec{b})\,\vec{I}\right] + 
      \left[\nabla\vec{b} - \transpose{(\nabla\vec{b})}\right]\,\vec{u} = \vec{0}.
\end{equation}
Furthermore, we add a similar nonconservative contribution as the last term on the left-hand side 
of \eqref{eq:nonconservativeb} to the momentum equation, which then becomes
\begin{equation} 
\label{eq:nonconservativeu}
    \partial_t\left(\rho\,\vec{u}\right) 
            + \divergence\left(\rho\,\vec{u}\,\vec{u} + p\,\vec{I} + \vec{\Omega}\right) + 
            \left[\transpose{(\nabla\vec{b})} - \nabla\vec{b}\right]\,
            \sigma\,\dfrac{\vec{b}}{\norm{\vec{b}}} = \vec{0}  
\end{equation}
and to the energy equation, formally accounting for the work due to the newly introduced forces
\begin{equation} 
\label{eq:nonconservativeE}
  \partial_t\left(\rho\,E\right) 
            + \divergence\left[\left(\rho\,E + p\right)\,\vec{u} + \vec{\Omega}\cdot\vec{u}\right] +
            \left[\transpose{(\nabla\vec{b})} - 
            \nabla\vec{b}\right]\,\sigma\,\dfrac{\vec{b}}{\norm{\vec{b}}}\cdot\vec{u} = 0.
\end{equation}
The modified model with Godunov--Powell-type symmetrising nonconservative products is
then written compactly as
\begin{equation} 
\label{eq:model}
    \partial_t
        \begin{pmatrix*}[c]
                   \vphantomi \alpha_1\,\rho_1\\[1.0mm]
                   \vphantomi \alpha_2\,\rho_2\\[1.0mm]
                   \vphantomi \rho\,\vec{u}\\[1.0mm]
                   \vphantomi \rho\,E\\[1.0mm]
                   \vphantomi \alpha_1\\[1.0mm]
                   \vphantomi \vec{b}\\[1.0mm]
                   \vphantomi c\\[1.0mm]
        \end{pmatrix*} + 
    \divergence
        \begin{pmatrix*}[c]
                   \vphantomi \alpha_1\,\rho_1\,\vec{u}\\[1.0mm]
                   \vphantomi \alpha_2\,\rho_2\,\vec{u}\\[1.0mm]
                   \vphantomi \rho\,\vec{u}\otimes\vec{u} + p\,\vec{I} + \vec{\Omega}\\[1.0mm]
                   \vphantomi \left(\rho\,E + p\right)\,\vec{u} + \vec{\Omega}\,\vec{u}\\[1.0mm]
                   \vphantomi 0\\[1.0mm]
                   \vphantomi (\vec{u}\cdot\vec{b})\,\vec{I}\\[1.0mm]
                   \vphantomi 0\\[1.0mm]
        \end{pmatrix*} + 
    \begin{pmatrix*}[c]
               \vphantomi 0\\[1.0mm]
               \vphantomi 0\\[1.0mm]
               \vphantomi \left[\transpose{(\nabla\vec{b})} - 
                   \nabla\vec{b}\right]\,\sigma\,{\vec{b}}/{\norm{\vec{b}}}\\[1.0mm]
               \vphantomi \left[\transpose{(\nabla\vec{b})} - 
                   \nabla\vec{b}\right]\,\sigma\,{\vec{b}}/{\norm{\vec{b}}}\cdot\vec{u}\\[1.0mm]
               \vphantomi \vec{u} \cdot \nabla \alpha_1 - K\,\divergence\vec{u}\\[1.0mm]
               \vphantomi \left[\nabla\vec{b} - \transpose{(\nabla\vec{b})}\right]\,\vec{u}\\[1.0mm]
               \vphantomi \vec{u} \cdot \nabla c\\[1.0mm]
    \end{pmatrix*} = \vec{0}.
\end{equation}
As a result of the fact that all the new nonconservative terms in \eqref{eq:nonconservativeb},
\eqref{eq:nonconservativeu}, and \eqref{eq:nonconservativeE} evaluate to zero if the compatibility
condition \eqref{eq:compatibilitytensor} is fulfilled, the formulation
\eqref{eq:model} is, at least for smooth solutions on the continuum level, entirely equivalent to the model
\eqref{eq:gavrilyuk}. Yet, the important difference is that now a full set of linearly independent eigenvectors
(given in the following subsection) can be obtained for this new form of the system, and thus one
can prove the strong hyperbolicity of the model. 
We will then discuss in Section~\ref{sec:testproblems} the different behaviour that the 
two formulations exhibit at the discrete level.  

Finally, it is necessary to emphasise that the transformations described above do not
ruin the thermodynamic compatibility of the governing PDEs, that is the over-determined
system~\eqref{eq:model}, together with an appropriate entropy equation, 
still forms a compatible system, \textit{e.g.} see \cite{SHTC-GENERIC-CMAT}. 

\subsection{Eigenstructure of the strongly hyperbolic Godunov--Powell-type model} 
\label{sec:hyperbolicity}
By defining a vector of conserved variables $\vec{Q}$ and a vector of primitive variables $\vec{V}$ as 
\begin{equation}
  \vec{Q} =
  \transpose{\left(
  \alpha_1\,\rho_1,\ \alpha_2\,\rho_2,\ \rho\,\transpose{\vec{u}},\ \rho\,E,\ \alpha_1,\ 
  \transpose{\vec{b}},\ c\right)},\qquad
  \vec{V} =
  \transpose{\left(
  \rho_1,\ \rho_2,\ \transpose{\vec{u}},\ p,\ \alpha_1,\ \transpose{\vec{b}},\ c\right)
  }
\end{equation}
the governing PDE system \eqref{eq:model} can be written in compact matrix-vector notation as 
\begin{equation}
    \frac{\partial \vec{Q}}{\partial t } + \nabla \cdot \vec{F}(\vec{Q}) + 
        \vec{B}(\vec{Q}) \, \nabla \vec{Q}=\vec{0},  
\end{equation}
where $\vec{F}(\vec{Q})$ is a nonlinear flux tensor and $\vec{B}(\vec{Q}) \, \nabla
\vec{Q}$ accounts for the non-conservative products. The quasi-linear form of the PDE in terms of
the conservative variables $\vec{Q}$ reads 
\begin{equation}
  \frac{\partial \vec{Q}}{\partial t } + \vec{A}(\vec{Q}) \, \nabla \vec{Q}=\vec{0}, \qquad 
  \text{with} \qquad \vec{A}(\vec{Q}) = \frac{\partial \vec{F}}{\partial \vec{Q}} \, 
  \nabla \vec{Q} + \vec{B}(\vec{Q}).
\end{equation}
In terms of the vector of primitive variables $\vec{V}$ it can be rewritten as
\begin{equation} \label{eq:primitivequasilinear}
   \frac{\partial \vec{V}}{\partial t } + \vec{C}(\vec{V}) \, \nabla \vec{V}=\vec{0}, 
   \qquad \text{with}  \qquad  \vec{C}(\vec{V}) = \pd{\vec{V}}{\vec{Q}} \left(\pd{\vec{F}}{\vec{V}} + 
   \vec{B}\,\pd{\vec{Q}}{\vec{V}}\right). 
\end{equation}
Due to the rotational invariance of \eqref{eq:model}, in order to compute its eigenstructure, and thus 
assess its
hyperbolicity, it will be sufficient to project the equations along a generic $x$ direction
specified by a unit vector $\hat{\mathbf{e}}_x$, so that the matrix of coefficients appearing in
\eqref{eq:primitivequasilinear} has a projection ${\mathbf{C}_1 = \mathbf{C} \, \hat{\mathbf{e}}_x}$
which reads 
\begin{equation}
\begingroup
\setlength{\arraycolsep}{3.9pt} 
\vec{C}_1 = 
    \begin{pmatrix*}
        \vphantomi u & 0 & \dfrac{(\alpha_1 + K)\,\rho_1}{\alpha_1} &   0 &   0 &               0 & 0 & 0 & 0 & 0 & 0 \\[4.0mm]
        \vphantomi 0 & u & \dfrac{(\alpha_2 - K)\,\rho_2}{\alpha_2} &   0 &   0 &               0 & 0 & 0 & 0 & 0 & 0 \\[4.0mm]
        \vphantomi 0 & 0 &                                        u &   0 &   0 & \dfrac{1}{\rho} & 0 & \dfrac{\sigma\,b_1\,(b_2^2 + b_3^2)}{\rho\,\babs^3} & \dfrac{\sigma\,b_1^2\,b_2}{\rho\,\babs^3} & \dfrac{\sigma\,b_3\,b_1^2}{\rho\,\babs^3} & 0 \\[4.0mm]
        \vphantomi 0 & 0 &                                        0 &   u &   0 &               0 & 0 & \dfrac{\sigma\,b_2\,(b_2^2 + b_3^2)}{\rho\,\babs^3} & \dfrac{\sigma\,b_1\,b_2^2}{\rho\,\babs^3} & \dfrac{\sigma\,b_1\,b_2\,b_3}{\rho\,\babs^3} & 0 \\[4.0mm]
        \vphantomi 0 & 0 &                                        0 &   0 &   u &               0 & 0 & \dfrac{\sigma\,b_3\,(b_2^2 + b_3^2)}{\rho\,\babs^3} & \dfrac{\sigma\,b_1\,b_2\,b_3}{\rho\,\babs^3} & \dfrac{\sigma\,b_1\,b_3^2}{\rho\,\babs^3} & 0 \\[4.0mm]
        \vphantomi 0 & 0 &                                \rho\,a^2 &   0 &   0 &               u & 0 & 0 & 0 & 0 & 0 \\[1.0mm]
        \vphantomi 0 & 0 &                                       -K &   0 &   0 &               0 & u & 0 & 0 & 0 & 0 \\[1.0mm]
        \vphantomi 0 & 0 &                                      b_1 & b_2 & b_3 &               0 & 0 & u & 0 & 0 & 0 \\[1.0mm]
        \vphantomi 0 & 0 &                                        0 &   0 &   0 &               0 & 0 & 0 & u & 0 & 0 \\[1.0mm]
        \vphantomi 0 & 0 &                                        0 &   0 &   0 &               0 & 0 & 0 & 0 & u & 0 \\[1.0mm]
        \vphantomi 0 & 0 &                                        0 &   0 &   0 &               0 & 0 & 0 & 0 & 0 & u \\[1.0mm]
    \end{pmatrix*}.
\endgroup
\end{equation}
The matrix collects the flux Jacobian with respect to the primitive variables as well as the matrix
of coefficients for the nonconservative products, also written in terms of gradients of the
primitive variables. For $\vec{C}_1$ we computed the following eigenvalues
\begin{equation}
    \vec{\lambda} = \begin{pmatrix*}[c]
    u\\[1.0mm]
    u\\[1.0mm]
    u\\[1.0mm]
    u\\[1.0mm]
    u\\[1.0mm]
    u\\[1.0mm]
    u\\[1.0mm]
    u - \sqrt{k_1 + k_3}\\[1.0mm]
    u + \sqrt{k_1 + k_3}\\[1.0mm]
    u - \sqrt{k_1 - k_3}\\[1.0mm]
    u + \sqrt{k_1 - k_3}\\[1.0mm]
    \end{pmatrix*}, \qquad \textnormal{with}  \qquad 
    \left\{
    \begin{aligned}
        & k_1 = \dfrac{a^2 + a_\sigma^2}{2},\\[1.0mm]
        & k_2 = \dfrac{a^2 - a_\sigma^2}{2},\\[1.0mm]
        & k_3 = \sqrt{k_1^2 - (1-\beta_1^2)\,a^2\,a_\sigma^2},\\[3.0mm]
        & a_\sigma^2 = \frac{\sigma}{\rho}\,\babs\,(1 - \beta_1^2),\\[1.0mm]
        & \beta_1 = {\dfrac{b_1}{\babs}},\quad \beta_2 = {\dfrac{b_2}{\babs}},\quad 
            \beta_3 = {\dfrac{b_3}{\babs}},\\[1.0mm]
    \end{aligned}
    \right.
\end{equation}
and with $a$ being the Wood \cite{wood1930} speed of sound for the mixture defined by
\begin{equation}
    a = \sqrt{\frac{\rho_1\,a_1^2\,\rho_2\,a_2^2}{\rho\,\left(\alpha_1\,\rho_2\,
        a_2^2 + \alpha_2\,\rho_1\,a_1^2\right)}}.
\end{equation}
The model includes seven contact waves moving with the fluid velocity $u$, 
and four mixed capillarity/pressure waves.
The eleven linearly independent right eigenvectors of the $\vec{C}_1$ matrix are 
\begin{equation}
\vec{R} = 
 \begin{pmatrix*}[c]
0 & 0              & 0              & 0 & 0        & 0 & 1 & -\beta_1\,\dimi        & -\beta_1\,\dimi        & -\beta_1\,\dimi        & -\beta_1\,\dimi\\[1.0mm]
0 & 0              & 0              & 0 & 0        & 1 & 0 & \beta_1\,\dimii        & \beta_1\,\dimii        & \beta_1\,\dimii        & \beta_1\,\dimii\\[1.0mm]
0 & 0              & 0              & 0 & 0        & 0 & 0 & \beta_1\,\kfiv         & -\beta_1\,\kfiv        & \beta_1\,\kfv          & -\beta_1\,\kfv\\[1.0mm]
0 & 0              & 0              & 0 & -\beta_3 & 0 & 0 & -\beta_2\,\kfvi        & \beta_2\,\kfvi         & -\beta_2\,\kfvii       & \beta_2\,\kfvii\\[1.0mm]
0 & 0              & 0              & 0 & \beta_2  & 0 & 0 & -\beta_3\,\kfvi        & \beta_3\,\kfvi         & -\beta_3\,\kfvii       & \beta_3\,\kfvii\\[1.0mm]
0 & 0              & 0              & 0 & 0        & 0 & 0 & -\beta_1\,\rho\,\asq   & -\beta_1\,\rho\,\asq   & -\beta_1\rho\,\asq     & -\beta_1\,\rho\,\asq\\[1.0mm]
0 & 0              & 0              & 1 & 0        & 0 & 0 & \beta_1\,K             & \beta_1\,K             & \beta_1\,K             & \beta_1\,K\\[1.0mm]
0 & \beta_1\beta_3 & \beta_1\beta_2 & 0 & 0        & 0 & 0 & \kfiv\,\kfvi\,\dimviii & \kfiv\,\kfvi\,\dimviii & \kfv\,\kfvii\,\dimviii & \kfv\,\kfvii\,\dimviii\\[1.0mm]
0 & 0              & 1 - \beta_1^2  & 0 & 0        & 0 & 0 & 0                      & 0                      & 0                      & 0\\[1.0mm]
0 & 1 - \beta_1^2  & 0              & 0 & 0        & 0 & 0 & 0                      & 0                      & 0                      & 0\\[1.0mm]
1 & 0              & 0              & 0 & 0        & 0 & 0 & 0                      & 0                      & 0                      & 0\\[1.0mm]
\end{pmatrix*},
\end{equation}
with
\begin{equation}
    \begin{array}{llll}
        \kfiv = \sqrt{k_1 + k_3},&
        \kfv = \sqrt{k_1 - k_3},&
        \kfvi = (k_2 - k_3)/\kfiv,&
        \kfvii = (k_2 + k_3)/\kfv,\\[2mm]
        \dimi = (K + \alpha_1)\,\rho_1/\alpha_1,&
        \dimii = (K - \alpha_2)\,\rho_2/\alpha_2,&
        \dimviii = \babs/\assq. & 
    \end{array}
\end{equation} 
We can then conclude that, on smooth solutions, the hyperbolicity  of the surface tension model 
forwarded in
\cite{Berry2008a, Schmidmayer2017} can be restored by including the proposed Godunov--Powell
symmetrising nonconservative products.

\subsection{Hyperbolic curl cleaning with the generalized Lagrangian multiplier (GLM) approach} 

The modified PDE system discussed in the previous sections, which allows to restore strong
hyperbolicity compared to the original model of Schmidmayer \etal\  \cite{Schmidmayer2017},
very closely follows the ideas of Godunov \cite{God1972MHD} and Powell \etal\ 
\cite{Powell1997,powell1999} concerning the symmetrisation and the numerical treatment of the 
divergence-free
condition of the magnetic field in the MHD system, respectively. 

An alternative and very successful numerical treatment of the divergence-free condition of the
magnetic field for the Maxwell and MHD equations is the so-called generalized Lagrangian multiplier
(GLM) approach of Munz \etal\  \cite{MunzCleaning,Dedneretal}, which consists in a
\textit{hyperbolic} divergence \textit{cleaning} achieved by adding a new auxiliary scalar field to
the PDE system, whose role is to transport divergence errors out of the computational domain via 
acoustic-type waves, so that they cannot accumulate locally. In the following, we adapt the GLM approach to the system
\eqref{eq:gavrilyuk} with the curl involution $\nabla \times \mathbf{b}= 0$. The augmented GLM version 
of the system reads  
\begin{subequations} \label{eq:gavrilyuk.glm}
    \begin{align}
    &\partial_t\left(\alpha_1\,\rho_1\right) 
    + \divergence\left(\alpha_1\,\rho_1\,\vec{u}\right) = 0,\label{eq:glmrho1} \\[1mm] 
    &\partial_t\left(\alpha_2\,\rho_2\right) 
    + \divergence\left(\alpha_2\,\rho_2\,\vec{u}\right) = 0, \\[1mm] 
    &\partial_t\left(\rho\,\vec{u}\right) 
    + \divergence\left(\rho\,\vec{u}\otimes\vec{u} + p\,\vec{I} + \vec{\Omega}\right) = 
    \vec{0}, \\[1mm] 
    &\partial_t\left(\rho\,E\right) 
    + \divergence\left[\left(\rho\,E + p\right)\,\vec{u} + \vec{\Omega}\,\vec{u}\right] = 
    0,\label{eq:glmenergyconservative}\\[1mm] 
    &\partial_t\left(\alpha_1\right) + \vec{u} \cdot \nabla \alpha_1  - K\,\divergence\vec{u} = 
    0,\\[1mm] 
    &\partial_t\left(\vec{b}\right) + \nabla \cdot \left[ \left(\vec{u} \cdot \vec{b}\right)\,\vec{I} \right] + 
    \left[ \nabla \vec{b} - \transpose{\left(\nabla \vec{b}\right)} \right]\, \vec{u} + 
    c_h\, \nabla \times \vec{\psi} = \vec{0}, \\[1mm] 
    &\partial_t\left(c\right) + \vec{u}\cdot\nabla c = 0, \\[1mm] 
    &\partial_t \left(\vec{\psi}\right) + \vec{u} \cdot \nabla{\vec{\psi}}  -
     c_h \, \nabla \times \vec{b} = - \kappa \, \vec{\psi},\label{eq:glmpsi} 
    \end{align}
\end{subequations}
where $c_h$ is the artificial wave speed associated with the hyperbolic curl cleaning process and
$\kappa$ is a small damping parameter, which in the present work is always set as $\kappa = 0$.
For a similar approach applied to a first order hyperbolic reduction of the Einstein field equations, 
see \cite{Dumbser2019}. Note the curl-curl structure in the equations for $\mathbf{b}$ and the 
cleaning field $\boldsymbol{\psi}$, which have a Maxwell-type form, i.e. in the augmented GLM curl cleaning 
system, the constraint violations in the vector field $\mathbf{b}$ are transported away via 
electromagnetic-type waves. 
It is easy to see that in the limit $c_h \to \infty$ one obtains $\nabla \times \vec{b} \to 0$. Due to the
presence of the transport term $\vec{u}\cdot\nabla\vec{\psi}$ in the evolution equation
\eqref{eq:glmpsi}, which is needed in order to have a Galilean invariant system, the cleaning 
vector field $\vec{\psi}$, unlike in \cite{Dumbser2019}, does not obey an additional linear 
divergence-free involution, and thus we chose not to enforce any additional constraints on the 
cleaning field itself. 

Note that, similar to \cite{Derigs2018}, in order to account for the effects of curl-cleaning on 
the total energy balance, one 
should in principle replace the energy conservation equation \eqref{eq:glmenergyconservative} with
\begin{equation} \label{eq:glmenergync}
    \partial_t\left(\rho\,E\right) 
    + \divergence\left[\left(\rho\,E + p\right)\,\vec{u} + \vec{\Omega}\,\vec{u}\right] + 
        c_h\,\frac{\vec{b}\cdot\nabla\times\vec{\psi}}{\norm{\vec{b}}} = 0. 
\end{equation}
Nonetheless, the computations shown in this work are carried out retaining the formally conservative
equation \eqref{eq:glmenergyconservative}. In preliminary tests, we found negligible differences
between the results from the energy-consistent equation \eqref{eq:glmenergync} and from the formally
conservative system which neglects the correction given in Eq.~\eqref{eq:glmenergync}, and the basic
properties of the two systems are the same (namely both systems are hyperbolic, have the same
eigenvalues, and a full set of eigenvectors can be found in both cases). Likewise, formulations
including the Godunov--Powell nonconservative products, in combination with the GLM curl cleaning
equations have been tested and yielded results that are comparable with those obtained with GLM curl
cleaning alone. Furthermore, we also tested the equivalence at the discrete level of the interface
field equation in its original fully-nonconservative form \eqref{b.Galilean} with its partially
conservative discretisation according to Eq.~\eqref{eq:nonconservativeb}.

Hyperbolicity of the augmented GLM curl cleaning system \eqref{eq:gavrilyuk.glm} can be shown by
repeating the procedure carried out in Section~\ref{sec:hyperbolicity} to compute explicitly a set
of fourteen eigenvalues
\begin{equation} \label{eq:glmeigenvalues}
    \vec{\lambda} = \begin{pmatrix*}[c]
    u\\[1.0mm]
    u\\[1.0mm]
    u\\[1.0mm]
    u\\[1.0mm]
    u\\[1.0mm]
    u\\[1.0mm]
    u-c_h\\[1.0mm]
    u-c_h\\[1.0mm]
    u+c_h\\[1.0mm]
    u+c_h\\[1.0mm]
    u - \sqrt{k_1 + k_3}\\[1.0mm]
    u + \sqrt{k_1 + k_3}\\[1.0mm]
    u - \sqrt{k_1 - k_3}\\[1.0mm]
    u + \sqrt{k_1 - k_3}\\[1.0mm]
    \end{pmatrix*}, \quad \textnormal{with}  \quad 
    \left\{
    \begin{aligned}
        & \beta_1 = {\dfrac{b_1}{\babs}},\quad \beta_2 = 
            {\dfrac{b_2}{\babs}},\quad \beta_3 = {\dfrac{b_3}{\babs}},\qquad
        a_\sigma^2 = \frac{\sigma}{\rho}\,\babs\,(1 - \beta_1^2),\\[1.0mm]
        & k_1 = \dfrac{a^2 + a_\sigma^2}{2},\qquad
        k_2 = \dfrac{a^2 - a_\sigma^2}{2},\qquad
        k_3 = \sqrt{k_1^2 - (1-\beta_1^2)\,a^2\,a_\sigma^2},\\[4.0mm]
        & \kfiv = \sqrt{k_1 + k_3},\qquad
        \kfv = \sqrt{k_1 - k_3},\\[6mm]
        &\kfvi = (k_2 - k_3)/\kfiv,\qquad
        \kfvii = (k_2 + k_3)/\kfv,\\[6mm]
        & \kfviii = \left(k_2 + \beta_1^2\,\assq + k_3\right)/k_4,\qquad
        \kfix = \left(k_2 + \beta_1^2\,\assq - k_3\right)/k_5,\\[2.0mm]
        & \kfx = \frac{%
        \beta_1^2\,\asq + c_h^2\,\left\{1 + {(\gamma_1 - 1)\,(\gamma_2 - 
            1)}\,{\left[\alpha_2\,(\gamma_1 - 1) + \alpha_1\,(\gamma_2 - 1)\right]}^{-1}\right\}}{%
        c_h^2\,(\assq - c_h^2) + \asq\,\left[c_h^2 - \assq\,(1 - \beta_1^2)\right]},
    \end{aligned}
    \right.
\end{equation}
here reported together with some 
auxiliary variables, which, supplemented with
\begin{equation}
\label{eq:convenience2}
    \begin{array}{lll}
        \dimi = (K + \alpha_1)\,\rho_1/\alpha_1,&
        \dimii = (K - \alpha_2)\,\rho_2/\alpha_2,&
        \dimviii = \babs/\assq,\\[2mm]
           \khi = \beta_1^2 - \kfxi + (1 - \beta_1^2)\,\kfxii,\phantom{Mx} &
           \khii = 1 + \beta_2^2\,(\kfxii - 1),\phantom{Mx} &
          \khiii = 1 + \beta_3^2\,(\kfxii - 1) , \\[2mm]
           \khiv = 1 - \kfxii,\phantom{Mx} &
            \khv = 1 + \kfxi - \kfxii,\phantom{Mx} &
          \khvi = \beta_1\,\kfxi.\\
    \end{array}
\end{equation}
are used to write compactly the set of fourteen linearly independent right eigenvectors. 
The wave structure includes six transport fields (contact waves), four \emph{cleaning waves} with eigenvalues 
$\lambda_{c_h} = u \pm c_h$, and four waves of mixed capillary/acoustic nature with eigenvalues 
$\lambda_{p\,\sigma} = u \pm \sqrt{k_1 \pm k_3}$, which are the same obtained from the previous variants
of the mathematical model.
Recalling the definitions given in Equations~\eqref{eq:glmeigenvalues}~and~\eqref{eq:convenience2}, 
the first ten eigenvectors of the augmented GLM curl-cleaning model, associated with the transport 
and cleaning eigenvalues, are
\begin{equation}
\vec{R}_{1-10} = 
\begin{pmatrix*}[c]
0 & 0 & 0 & 0    & 0 & 1 & \khi\,\dimi                   & \khi\,\dimi                    & \khi\,\dimi                    & \khi\,\dimi                   \\[1.0mm]
0 & 0 & 0 & 0    & 1 & 0 & \khi\,\dimii                  & \khi\,\dimii                   & \khi\,\dimii                   & \khi\,\dimii                  \\[1.0mm]
0 & 0 & 0 & 0    & 0 & 0 & \khi\,c_h                     & \khi\,c_h                      & \khi\,c_h                      & \khi\,c_h                     \\[1.0mm]
0 & 0 & 0 & -b_3 & 0 & 0 & -\khii\,c_h\,\beta_1/\beta_2  & -\khiv\,c_h\,\beta_1\,\beta_2  & -\khii\,c_h\,\beta_1/\beta_2   & -\khiv\,c_h\,\beta_1\,\beta_2 \\[1.0mm]
0 & 0 & 0 & b_2  & 0 & 0 & \khiv\,c_h\,\beta_1\,\beta_3  & \khiii\,c_h\,\beta_1/\beta_3   & \khiv\,c_h\,\beta_1\,\beta_3   & \khiii\,c_h\,\beta_1/\beta_3  \\[1.0mm]
0 & 0 & 0 & 0    & 0 & 0 & \khv\,\rho\,c_h^2             & -\khv\,\rho\,c_h^2             & -\khv\,\rho\,c_h^2             & \khv\,\rho\,c_h^2             \\[1.0mm]
0 & 0 & 1 & 0    & 0 & 0 & \khi\,K                       & -\khi\,K                       & -\khi\,K                       & \khi\,K                       \\[1.0mm]
0 & 0 & 0 & 0    & 0 & 0 & \khvi\,\babs                  & -\khvi\,\babs                  & -\khvi\,\babs                  & \khvi\,\babs                  \\[1.0mm]
0 & 0 & 0 & 0    & 0 & 0 & \rho\,c_h^2/(\beta_2\,\sigma) & 0                              & -\rho\,c_h^2/(\beta_2\,\sigma) & 0                             \\[1.0mm]
0 & 0 & 0 & 0    & 0 & 0 & 0                             & -\rho\,c_h^2/(\beta_3\,\sigma) & 0                              & \rho\,c_h^2/(\beta_3\,\sigma) \\[1.0mm]
0 & 1 & 0 & 0    & 0 & 0 & 0                             & 0                              & 0                              & 0                             \\[1.0mm]
1 & 0 & 0 & 0    & 0 & 0 & 0                             & 0                              & 0                              & 0                             \\[1.0mm]
0 & 0 & 0 & 0    & 0 & 0 & 0                             & \rho\,c_h^2/(\beta_3\,\sigma)  & 0                              & \rho\,c_h^2/(\beta_3\,\sigma) \\[1.0mm]
0 & 0 & 0 & 0    & 0 & 0 & \rho\,c_h^2/(\beta_2\,\sigma) & 0                              & \rho\,c_h^2/(\beta_2\,\sigma)  & 0                             \\[1.0mm]
\end{pmatrix*},
\end{equation}
while the remaining four eigenvectors, corresponding to the capillary/acoustic waves are
\begin{equation}
\vec{R}_{11-14} = 
\begin{pmatrix*}[c]
\kfviii\,\dimi           & \kfviii\,\dimi          & \kfix\,\dimi             & \kfix\,\dimi\\[1.0mm]
-\kfviii\,\dimii         & -\kfviii\,\dimii        & -\kfix\,\dimii           & -\kfix\,\dimii\\[1.0mm]
-\kfviii\,\kfiv          & \kfviii\,\kfiv          & -\kfix\,\kfv             & \kfix\,\kfv\\[1.0mm]
-\assq\,\beta_1\,\beta_2 & \assq\,\beta_1\,\beta_2 & -\assq\,\beta_1\,\beta_2 & \assq\,\beta_1\,\beta_2\\[1.0mm]
-\assq\,\beta_1\,\beta_3 & \assq\,\beta_1\,\beta_3 & -\assq\,\beta_1\,\beta_3 & \assq\,\beta_1\,\beta_3\\[1.0mm]
\kfviii\,\rho\,\asq      & \kfviii\,\rho\,\asq     & \kfix\,\rho\,\asq        & \kfix\,\rho\,\asq\\[1.0mm]
-\kfviii\,K              & -\kfviii\,K             & -\kfix\,K                & -\kfix\,K\\[1.0mm]
\kfiv\,\beta_1\,\babs    & \kfiv\,\beta_1\,\babs   & \kfv\,\beta_1\,\babs     & \kfv\,\beta_1\,\babs\\[1.0mm]
0                        & 0                       & 0                        & 0\\[1.0mm]
0                        & 0                       & 0                        & 0\\[1.0mm]
0                        & 0                       & 0                        & 0\\[1.0mm]
0                        & 0                       & 0                        & 0\\[1.0mm]
0                        & 0                       & 0                        & 0\\[1.0mm]
0                        & 0                       & 0                        & 0\\[1.0mm]
\end{pmatrix*}.
\end{equation}
We can therefore conclude that the augmented GLM system \eqref{eq:gavrilyuk.glm} is strongly hyperbolic. 
However, its \textit{major advantage} over the Godunov-Powell-type system is that the GLM system is 
fully conservative, while the Godunov-Powell system is not, at least when standard general-purpose 
schemes are used that do not satisfy the curl involution constraint exactly at the discrete level. 

\subsection{Exact equilibrium solution for a symmetric droplet with diffuse interface} 
\label{sec:exactsol}
A steady state solution can be easily obtained for a two-dimensional water column or a
three-dimensional droplet (hereafter we will take the liberty to call \textit{droplets}  
the two-dimensional objects as well) by first assigning a radial profile $c(r) = c(\norm{\vec{x}})$ for the
interface between the two phases, then computing the corresponding interface field $\vec{b} = \nabla
c$ and balancing the surface tension forces, which are known once a specific geometrical
configuration is chosen, with the pressure field.

For convenience, we define the dimensionless radial coordinate {$\rs = r/R = \norm{\vec{x}}/R$},
with $R$ being the radius of the water column or droplet. Here, with the notation $\vec{x}$ we
indicate the Cartesian position vector, independently from the number of space dimensions $d$. We
then set the radial profile of the colour function to be a smoothed {Heaviside} step function of the
form
\begin{equation}
    c(r_\ast) = \frac{1}{2}\,\erfc\left(\frac{\rs - 1}{\keps}\right),
\end{equation}
with the dimensionless interface thickness parameter $\keps$ controlling the intensity of the
smoothing. The Cartesian gradient of the colour function is immediately computed as
\begin{equation} \label{eq:exactb}
    \vec{b} (\vec{x}) = -\frac{\vec{x}}{\sqrt\pi\,k_\epsilon\,R\,\norm{\vec{x}}}
        \exp\left[-{\left(\frac{r_\ast - 1}{k_\epsilon}\right)}^2\right].
\end{equation}
It can be verified by easy calculations that 
the dimensionless interface scaling parameter $\keps$ corresponds to four times the standard 
deviation of the Gaussian curve
describing the profile of the interface energy $\sigma\,\norm{\vec{b}}$ along the radial direction, 
rendered dimensionless with respect to the nominal radius of the droplet. To give a clear physical 
meaning to the quantity, one can say that in the region of 
space bounded by $1-\keps/2 \le r_\ast \le 1 + \keps/2$, about $95.5\%$ of the surface energy is stored.
In a uniform flow, all the governing equations are satisfied for any choice of the density and
volume fraction fields, and one can compute the radial pressure profile from the momentum equation
by requiring that the pressure gradient be balanced with the divergence of the surface tension
tensor $\vec{\Omega}$. Clearly from a physical/geometrical standpoint, the colour function and volume fraction
variables are closely related and cannot be set independently.
One can then derive from the momentum equation 
\begin{equation}
    \nabla p + \sigma\,\vec{b}\,\nabla\cdot\frac{\vec{b}}{\norm{\vec{b}}} = 0
\end{equation}
a simple ordinary differential equation 
\begin{equation}
    \od{p}{r}(r) = \partial_r p(r) = -\frac{\vec{x}}{r}\cdot\left(\nabla\cdot\vec{\Omega}\right),
\end{equation}
that, by evaluating the divergence of the capillarity tensor $\vec{\Omega}$ from \eqref{eq:surfacetension2} 
and substituting the ansatz for the interface field \eqref{eq:exactb} yields
\begin{equation} \label{eq:odp}
    \od{p}{r_\ast}(r_\ast) = -(d - 1)\,\frac{\sigma}{R}\,\frac{1}{\sqrt\pi\,\keps\,r_\ast}\,
    \exp\left[-{\left(\frac{r_\ast - 1}{k_\epsilon}\right)}^2\right],
\end{equation}
Note that the ODE \eqref{eq:odp} and thus the pressure profile depend parametrically on the group $(d-1)\,\sigma/R$ 
and are otherwise solely functions of the geometry expressed through Eq.~\eqref{eq:exactb}.
One can then directly integrate \eqref{eq:odp} with atmospheric pressure $p_\up{atm}$
as a far field boundary condition in order to 
obtain the equilibrium pressure field
\begin{equation} \label{eq:compactpressureprofile}
    p\left(\rs\right) = p_\up{atm} + (d - 1)\,\frac{\sigma}{R}\,\int_{\rs}^{\infty} 
    \frac{1}{\sqrt\pi\,\keps\,\rs^\prime}\,\exp\left[-{\left(\frac{\rs^\prime - 
        1}{\keps}\right)}^2\right]\,\de{\rs^\prime},
\end{equation}
where $\rs^\prime$ is an auxiliary integration variable. 
The integral can be computed to machine precision with the aid of a Gauss--Legendre quadrature rule
with the precaution of defining a sufficiently refined integration mesh in the interface region.

\begin{figure}[!b]
    \includegraphics[draft=false, scale=\figurescalefactor]{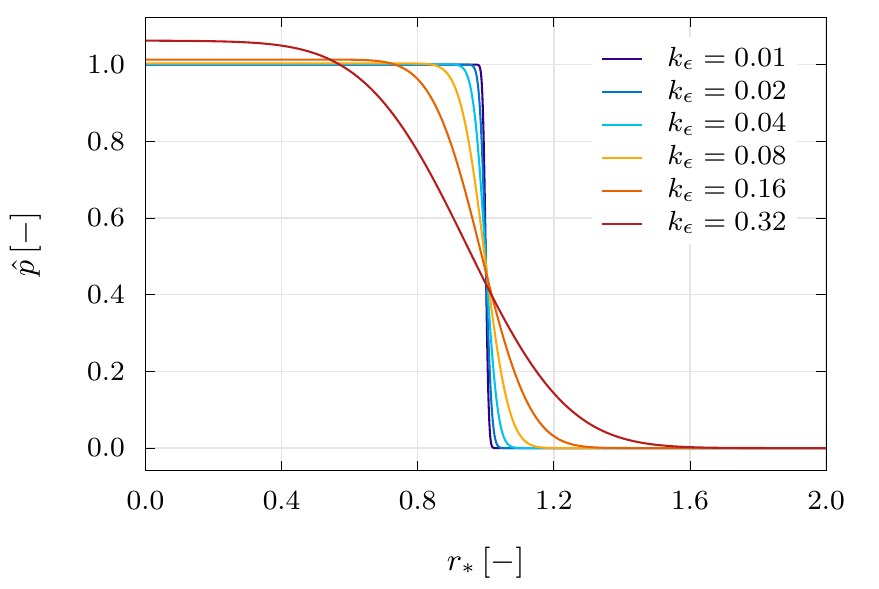}%
    \includegraphics[draft=false, scale=\figurescalefactor]{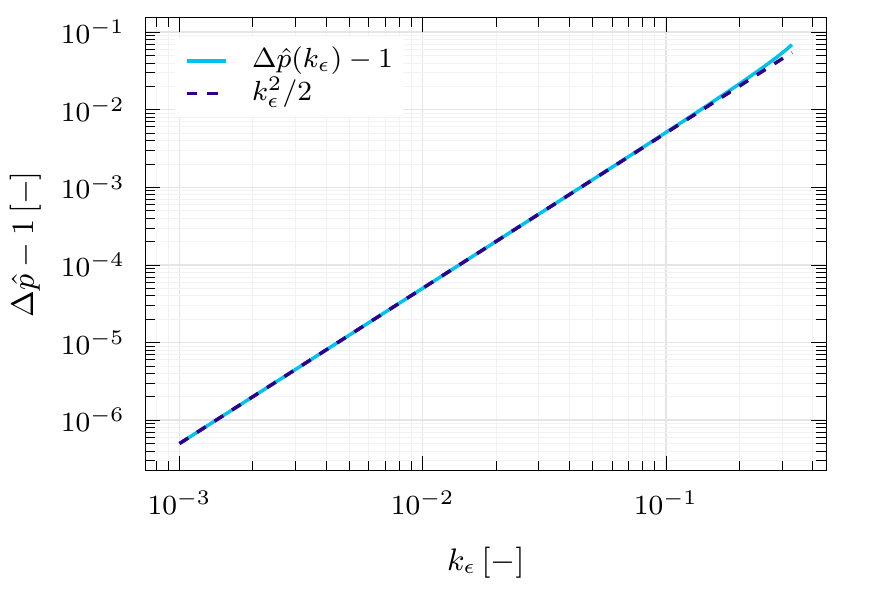}
    \caption{Exact pressure profiles for two-dimensional droplets with diffuse interface. On the
    left, the radial profile of the dimensionless pressure $\hat{p} = (p - p_\up{atm})\,R/[(d -
    1)\,{\sigma}]$ is shown for a range of values of the interface thickness parameter $k_\epsilon$. On the
    right, we plot the error of the dimensionless pressure jump $\Delta\hat{p}$
    of a \emph{smooth droplet} with respect to the dimensionless Young--Laplace pressure jump (unity) as a function of the 
    interface thickness parameter $k_\epsilon$. The dashed line represents a power law approximation of
    the curve.} 
    \label{fig:droplet-exact-solution-double}
\end{figure}

\subsubsection{Consistency with the Young--Laplace law.}
In the limit of vanishing interface thickness ($k_\epsilon \rightarrow 0$), one can verify 
that Eq.~\eqref{eq:compactpressureprofile}
yields a pressure jump between atmospheric condition and the centre of the droplet of the form
\begin{equation}
    \lim_{k_\epsilon \to 0} \Delta p(k_\epsilon) = \lim_{k_\epsilon \to 0} \left[
    \ev{p\left(r_\ast = 0\right)}{k_\epsilon} - p_\up{atm}\right] = (d - 1)\,\frac{\sigma}{R},
\end{equation}
so that the well-known Young--Laplace formula is recovered in the limit of a sharp interface. Also, we can
point out that, even for heavily smoothed droplets, the Young--Laplace formula provides similar
values to the ones obtained from Eq.~\eqref{eq:compactpressureprofile}, as can be seen in the left
panel of Figure~\ref{fig:droplet-exact-solution-double}, and that both estimates for the pressure
jump converge to the same value quite quickly: from the right panel of
Figure~\ref{fig:droplet-exact-solution-double} it is apparent that the following approximation holds
\begin{equation}
  \Delta p(k_\epsilon) \simeq (d - 1)\,\frac{\sigma}{R}\,\left(1+ \frac{1}{2}\,k_\epsilon^2\right),
  \label{eqn.DIYL} 
\end{equation}
which means that the pressure jump is overall affected only by relatively small deviations
from the Young--Laplace law, even for droplets with rather large diffuse interface region and 
converges to the sharp interface reference solution quadratically as the interface thickness 
$k_\epsilon$ vanishes 

\section{Numerical method}
\label{sec:numericalmethod}
In this section, we summarise the key elements of the family of numerical methods employed in this
work, which are the high order ADER \disgal (\pnpm{N}{N}) and ADER-WENO \finvol (\pnpm{0}{M})
schemes with \aposteriori subcell limiting. These methods can be applied to 
general nonconservative
hyperbolic systems of balance laws. Discontinuous Galerkin schemes for nonconservative hyperbolic systems 
have been introduced in \cite{dumbser2009, dumbser2010, rhebergen2008}, following the path-conservative 
approach of Castro and Par\'es, originally developed for the \finvol framework \cite{castro2006, pares2006}
and based on the theory of Dal Maso, Le Floch and Murat \cite{DLMtheory}. 
A recent review of the history of the development of 
ADER schemes can be found in \cite{Busto2019}.
For details on some of the first 
fully discrete one-step Lax-Wendroff-type time discretisations proposed for DG methods, see 
\cite{qiu2005, dumbser2006b, stedg1, gassner2011}.

Modern explicit ADER schemes follow a fully discrete predictor-corrector procedure, which can be
regarded as a high order extension of the simple and successful MUSCL-Hancock approach 
\cite{vanleer1974,vanleer1979,toro2009}, 
rather than using a semi-discrete formulation in conjunction with multi-stage timestepping, as 
in Runge--Kutta DG methods \cite{CBS-convection-dominated}. First, a predictor step evolves
the polynomial data \emph{in the small} to obtain an approximate space-time polynomial solution in 
each cell, without taking coupling with neighbouring cells into account. Then, volume integrals
arising from a weak formulation of the differential problem can be easily evaluated with the aid of
appropriate quadrature formulas, and quadrature at space-time faces are used to compute averaged
numerical fluxes corresponding to the Riemann problems arising from extrapolation of the predictor
solution from adjacent cells.

At each timestep, 
the cell-local space-time predictor solution $\vec{q}_h(\vec{x},\ t)$ is computed from a piecewise polynomial
reconstruction $\vec{w}_h(\vec{x})$ of cell average data (for FV methods), or is directly available from the
evolved piecewise polynomial data $\vec{w}_h(\vec{x}) = \vec{u}_h(\vec{x})$ for DG methods. Since the two 
families of schemes (FV and DG) use 
the same discrete data representation (nodal degrees of freedom of a Gauss--Legendre--Lagrange
polynomial), the space-time \disgal predictor can be formulated in full generality
for both, or even for the more general family of $\pnpm{N}{M}$ schemes \cite{dumbser2008b, 
dumbser2010pnpmcns}.

From the space-time predictor solution one can immediately compute all the volume integrals
appearing in the fully discrete, one-step update formulas \eqref{eq:dgupdate}; in particular, this
operation can be carried out quite conveniently thanks to the choice of employing a nodal basis
where the nodal values are located at the Gauss--Legendre quadrature nodes 
and the basis functions
{${\phi_m(\vec{\xi}) = \psix (\xi)\, \psiy (\eta)\, \psiz (\zeta)}$} are two-dimensional or
three-dimensional tensor products of the Lagrange polynomials interpolating the Gauss--Legendre
quadrature nodes.

In this work, spurious oscillations that typically occur when employing higher than first order linear 
schemes, see \cite{godunov1959}, are minimised as follows: in the case of \finvol methods, we employ 
a nonlinear WENO reconstruction procedure, while for \disgal schemes we adopt the \textit{a posteriori} 
subcell \finvol limiting strategy \cite{dumbser2014}, that is, at each timestep a candidate solution 
is computed without any limiter, and then afterwards, if this candidate solution violates one or more 
physical and numerical  admissibility criteria (floating point exceptions, violation of positivity, 
violation of a discrete maximum
principle), then it is discarded and a new discrete solution is recomputed, starting again from valid data 
at the previous time step. This data is obtained by projecting the DG polynomial on a fine subcell 
grid, or directly from the subcell average representation of the data if it was already available 
at the previous timestep. 
Afterwards, the discrete solution is reconstructed back from subcell averages to a DG polynomial.  

\subsection{Data representation and notation}
The computational domain is partitioned in conforming Cartesian elements
\begin{equation}
\Omega_{ijk} = \left[ x_i - \frac{\Delta x_i}{2},\ x_i + \frac{\Delta x_i}{2} \right] \times
                \left[ y_j - \frac{\Delta y_j}{2},\ y_j + \frac{\Delta y_j}{2} \right] \times
                \left[ z_k - \frac{\Delta z_k}{2},\ z_k + \frac{\Delta z_k}{2} \right],
\end{equation}
and for each element a reference frame of coordinates is defined by
\begin{equation}
         \xi   = \frac{x - x_i}{\Delta x_i} + \halb,\qquad
         \eta  = \frac{y - y_j}{\Delta y_j} + \halb,\qquad
         \zeta = \frac{x - z_k}{\Delta z_k} + \halb.\qquad
\end{equation}
Discrete data are given as the degrees of freedom of a $d$-dimensional polynomial (in this
exposition we will use $d=3$ without loss of generality) represented by means of a set of nodal
basis functions in the form of three-dimensional tensor products of the Lagrange polynomials $\psix
(\xi)$, $\psiy (\eta)$, and $\psiz (\zeta)$ satisfying, at Gauss--Legendre quadrature node locations 
$\xi_{p_x}^\textsc{gl}$, $\eta_{p_y}^\textsc{gl}$, and $\zeta_{p_z}^\textsc{gl}$, 
the interpolation conditions 
$\psix(\xi_{p_x}^\textsc{gl}) = \delta_{m_x\,p_x}$, 
$\psiy(\eta_{p_y}^\textsc{gl}) = \delta_{m_y\,p_y}$, and 
$\psiz(\zeta_{p_z}^\textsc{gl}) = \delta_{m_z\,p_z}$, 
with $\delta_{ij}$ being the usual Kronecker symbol.

Throughout the paper we will use a compact multi-index notation so that the three-dimensional 
position of the generic Gauss--Legendre quadrature node of index $p$ is written 
${\vec{\xi}_p^\textsc{gl} = (\xi_{p_x}^\textsc{gl},\ \eta_{p_y}^\textsc{gl},\
\zeta_{p_z}^\textsc{gl})}$ and the three dimensional basis function of index $m$ can be expressed as
{${\phi_m(\vec{\xi}) = \psix (\xi)\, \psiy (\eta)\, \psiz (\zeta)}$}. Note that the interpolation
property can be written in an entirely analogous fashion with respect to the one-dimensional case,
that is, with the notation ${\phi_m(\vec{\xi}_p^\textsc{gl}) = \delta_{mp}}$.

Within each cell $\Omega_{ijk}$, at a given time $t = t^n$, the discrete solution is then written
(dropping for convenience of notation the cell indices $i,j,k$) as a 
polynomial of order $M$ in each direction 

\begin{equation}  \label{eq:uhansatz}
    \vec{u}_h(\vec{x},t^n) = \vec{u}_h(\vec{\xi}(\vec{x}),\ t^n) = 
    \sum_{m_x = 1}^{M+1}\sum_{m_y = 1}^{M+1}\sum_{m_z = 1}^{M+1}
    \psix(\xi)\,\psiy(\eta)\,\psiz(\zeta)\, \hat{\vec{u}}_{m_x\,m_y\,m_z}^n 
    = \sum_{m = 1}^{{(M+1)}^{d}} \phi_m(\boldsymbol{\xi})\, \hat{\vec{u}}_m^n.
\end{equation}
Furthermore, we will use the Einstein summation convention over repeated indices so that the discrete
solution can be expressed as $\vec{u}_h(\vec{x},\ t^n) =
\phi_m(\vec{\xi})\, \hat{\vec{u}}_m^n$. The \finvol data representation (cell
average values ${\vec{Q}}_{ijk}^n)$ can be regarded as a special case of \eqref{eq:uhansatz} in
which $M=0$ and the single basis function is the constant function $\phi_m^\textsc{fv}(\vec{\xi}) =
1$ within each element.

\subsection{Polynomial WENO reconstruction}
\label{sec:reconstruction}
In order to obtain a high order data reconstruction from \finvol cell averages, we employ a 
\emph{full polynomial}
WENO reconstruction, introduced in \cite{dumbser2007a} for unstructured meshes and employed, for
example, in \cite{dumbser2013a} on Cartesian grids. 
The most prominent difference between this approach and the original formulation by Jiang and 
Shu \cite{shu_efficient_weno}
is that instead of computing \emph{pointwise}
values of the conserved variables at the aid of optimal linear weights, here we seek to obtain the 
degrees of freedom of the \textit{entire reconstruction polynomial}, 
to be used in the computation of fluxes, non-conservative products and source terms via high order quadrature 
formulae. At this point we would also like to stress that entire WENO reconstruction polynomials with a 
reconstruction stencil of \textit{optimal compactness} can be achieved via the elegant CWENO approach 
forwarded by Puppo, Russo and Semplice \etal\  in \cite{LPR:99,LPR:2001,SCR:CWENOquadtree,cravero2018cweno,ADER_CWENO}.  
The first step is to define, for a
generic element $\Omega_{ijk} = [x_{i - {1}/{2}},\ x_{i + {1}/{2}}]\times[y_{j - {1}/{2}},\
y_{j + {1}/{2}}]\times[z_{k - {1}/{2}},\ z_{k + {1}/{2}}]$, the three sets of reconstruction
stencils, one for each space dimension. Each stencil will be identified by the triplet of subscripts
($i,\ j,\ k$) of the cell in which the reconstruction is sought, together with the superscript
describing the spatial direction ($x$, $y$ or $z$) of the reconstruction and an integer superscript
$s$ that identifies the specific stencil in the set. The three generic elements of such
reconstruction stencil sets will be then written as 
\begin{equation} \label{eq:stencils}
    S_{ijk}^{x,\, s} = \bigcup_{m = i-L}^{i+R}\Omega_{mjk},\qquad
    S_{ijk}^{y,\, s} = \bigcup_{m = j-L}^{j+R}\Omega_{imk},\qquad
    S_{ijk}^{z,\, s} = \bigcup_{m = k-L}^{k+R}\Omega_{ijm},\qquad
\end{equation}
where $L = L(M,\ s)$ and $R = R(M,\ s)$ are positive integers representing the number of elements in
the stencil respectively to the left and to the right of the principal cell. In each space
direction, for even values of the the reconstruction degree $M$, one always has $N_\up{s} = 3$
stencils, one central and two off-centre, with left and right extensions given by
\begin{equation}
    L(M,\ s) = 
    \begin{cases}
        M/2 & \text{if } s = 1, \\
        M & \text{if } s = 2, \\
        0 & \text{if } s = 3, \\
    \end{cases}\qquad
    R(M,\ s) = 
    \begin{cases}
        M/2 & \text{if } s = 1, \\
        0 & \text{if } s = 2, \\
        M & \text{if } s = 3, \\
    \end{cases}\qquad
\end{equation}
while for odd values of $M$ we define $N_\up{s} = 4$ types of stencils, two central, two off-centre, 
having extensions
\begin{equation}
    L(M,\ s) = 
    \begin{cases}
        (M+1)/2 & \text{if } s = 1, \\
        (M-1)/2 & \text{if } s = 2, \\
        M & \text{if } s = 3, \\
        0 & \text{if } s = 4, \\
    \end{cases}\qquad
    R(M,\ s) = 
    \begin{cases}
        (M-1)/2 & \text{if } s = 1, \\
        (M+1)/2 & \text{if } s = 2, \\
        0 & \text{if } s = 3, \\
        M & \text{if } s = 4. \\
    \end{cases}\qquad
\end{equation}

This choice ensures that each stencil be composed by a number of elements equal to the nominal order
of the scheme, which is $M+1$.

The dimension-by-dimension reconstruction is carried out by repeated application (over each space
dimension) of a one-dimensional-sweep procedure which constructs, in each cell, the $M+1$ degrees of
freedom of a polynomial of degree $M$, first solving a set of linear reconstruction equations
(imposing conservation of cell averages on each element of a given stencil), and then combining the
solutions of the reconstruction equations in a data-dependent, nonlinear fashion in order to ensure
the non-oscillatory character of the reconstructed polynomials.

In the generic element $\Omega_{ijk}$, the reconstruction polynomials obtained in each of the
three subsequent passes are expressed in terms of their degrees of freedom
$\hat{\vec{w}}^\textsc{1d}_{i\,j\,k,\,p}$, $\hat{\vec{w}}^\textsc{2d}_{i\,j\,k,\,p\,q}$, and
$\hat{\vec{w}}^\textsc{3d}_{i\,j\,k,\,p\,q\,r}$, relative to a tensor-product-type
Gauss--Legendre--Lagrange basis function. The degrees of freedom and thus the polynomials are
obtained as nonlinear convex combinations, written as
\begin{align}
    &\vec{w}_h^\textsc{1d}(x) = 
        \psi_p(\xi)\,\hat{\vec{w}}^\textsc{1d}_{ijk,\,p}, &&\hspace{-0cm}\text{with}\quad
    \hat{\vec{w}}^\textsc{1d}_{ijk,\,p} = 
        \omega_{s}^\textsc{1d}\,\hat{\vec{w}}^{\textsc{1d},\, s}_{ijk,\,p},\\
    &\vec{w}_h^\textsc{2d}(x,y) = 
        \psi_p(\xi)\,\psi_q(\eta)\,
        \hat{\vec{w}}^\textsc{2d}_{ijk,\,pq}, &&\hspace{-0cm}\text{with}\quad
    \hat{\vec{w}}^\textsc{2d}_{ijk,\,pq} = 
        \omega_s^\textsc{2d}\,\hat{\vec{w}}^{\textsc{2d},\, s}_{ijk,\,pq},\\
    &\vec{w}_h^\textsc{3d}(x,y,z) = 
        \psi_p(\xi)\,\psi_q(\eta)\,\psi_r(\zeta)\,
        \hat{\vec{w}}^\textsc{3d}_{ijk,\,pqr}, &&\hspace{-0cm}\text{with}\quad
    \hat{\vec{w}}^\textsc{3d}_{ijk,\,pqr} = 
        \omega_s^\textsc{3d}\,\hat{\vec{w}}^{\textsc{3d},\, s}_{ijk,\,pqr},
\end{align}
having defined, for each stencil $S_{ijk}^{x,\, s}$, $S_{ijk}^{y,\, s}$, and,
$S_{i\,j\,k}^{z,\, s}$ the linear (in the sense that they are not affected by limiting)
reconstruction polynomials computed from the solution of the reconstruction equations
\begin{align}
   &\vec{w}_h^{\textsc{1d},\, s}(x) = 
       \psi_p(\xi)\,\hat{\vec{w}}_{ijk,\,p}^{\textsc{1d},\, s},\\
   &\vec{w}_h^{\textsc{2d},\, s}(x,y) = 
       \psi_p(\xi)\,\psi_q(\eta)\,\hat{\vec{w}}_{i\,j\,k,\,p\,q}^{\textsc{2d},\, s},\\
   &\vec{w}_h^{\textsc{3d},\, s}(x,y,z) = 
       \psi_p(\xi)\,\psi_q(\eta)\,\psi_r(\zeta)\,\hat{\vec{w}}_{ijk,\,pqr}^{\textsc{3d},\, s},
\end{align}
where the indices $p$, $q$, and $r$ run from $0$ to $M$ (covering the number of degrees of freedom
to be reconstructed in the in each space dimension) and where we adopted Einstein's summation convention
over repeated indices.

Each pass is analogous to the first one, in that a one-dimensional stencil is used and only a
one-dimensional oscillation indicator has to be computed, but it must be remarked that, as
a result of performing the nonlinear stencil selection procedure in a given dimension before
operating the linear reconstructions in the next one, each one of the final three-dimensional
degrees of freedom is subject to composite limiting in the three space dimensions, which includes
information not only from the direct face neighbours, but from the node neighbours as well, and the
reconstruction is thus genuinely multi-dimensional. For alternative and more efficient multi-dimensional
finite volume WENO reconstructions on Cartesian meshes, see \cite{Buchmuller2014,Buchmuller2015}. 

More in detail, first, one constructs for each cell a one-dimensional polynomial
$\vec{w}_h^\textsc{1d}(x)$, while maintaining the data in the remaining directions ($y$ and $z$) in
piecewise constant cell-averaged form. 

The linear reconstruction equations, enforcing integral conservation on all elements in the stencil 
$S_{ijk}^{x,\, s}$, constitute a linear system whose solutions are the $M+1$ unknown degrees of
freedom $\hat{\vec{w}}_{ijk,\,p}^{1\textsc{d},s}$, and in the first space dimension read
\begin{equation}
   \dfrac{1}{\Delta x_m} \int_{x_{m -\frac{1}{2}}}^{x_{m +\frac{1}{2}}} 
   \psi_p(\xi(x))\,\hat{\vec{w}}^{\textsc{1d},\, s}_{ijk,\,p}\,\de{x} = 
   {\vec{Q}}_{mjk}^n,\quad\forall\,\Omega_{mjk} \in S_{ijk}^{x,\, s}.
\end{equation}
Then the nonlinear coefficients for the combination of the polynomials obtained from each stencil
are computed at each line sweep as
\begin{equation}
   \omega_s^{n\textsc{d}} = \tilde{\omega}_s^{n\textsc{d}}\,{\left(
   \sum_{s=1}^{N_s}\tilde{\omega}_s^{n\textsc{d}}\right)}^{-1},
   \quad\text{with}\quad \tilde{\omega}_s^{n\textsc{d}} = \lambda_s\,
   {\left(\sigma_s^{n\textsc{d}} + \epsilon\right)}^{-r},
\end{equation}
where the oscillation indicator $\sigma_s^{n\textsc{d}}$ is defined as 
\begin{equation}
   \sigma_s^{n\textsc{d}} = 
       \sum_{l=0}^M\sum_{m=0}^M
       \Sigma_{l\,m}\,\frac{\hat{\vec{w}}_{ijk,\,l}^{n\textsc{d},\,s}}{\vec{w}_0}\,
       \frac{\hat{\vec{w}}_{ijk,\,m}^{n\textsc{d},\,s}}{\vec{w}_0}, \quad \text{with} \quad
       \Sigma_{l\,m} = \sum_{\alpha = 1}^{M}\int_0^1
       \dfrac{\partial^\alpha\psi_l(\xi)}{\partial\xi^\alpha}\,
       \dfrac{\partial^\alpha\psi_m(\xi)}{\partial\xi^\alpha}\de{\xi}.
\end{equation}
The numerical parameters used for the computation of the nonlinear weights are $\lambda_s = 1$ for
off-centre stencils and $\lambda_s = 10^5$ for central stencils, and we set $\epsilon = 10^{-14}$
and $r = 8$. An important remark is that since the oscillation indicators are highly nonlinear,
particular care should be taken in dividing the input values by an appropriate scaling factor
$\vec{w}_0$. As a practical example, it is often the case that when using the stiffened gas equation
of state, very large values for the mixture energy variable appear even in standard pressure
conditions, which could lead to catastrophic loss of precision in the computation of the weights.
Such a scaling factor can be computed for example as
\begin{equation}
   \vec{w}_0 = \epsilon_0 + 
       \sum_{s}\sum_{m=0}^M \abs{\hat{\vec{w}}_{ijk,\,m}^{n\textsc{d},\,s}},
\end{equation}
that is, by evaluating, variable by variable, the sum of the absolute values of all the degrees of 
freedom of the input data over all stencils and adding a new constant $\epsilon_0 = 10^{-14}$ to 
avoid division by zero.

In the second pass, one obtains data in
the two-dimensional polynomial form $\vec{w}_h^{\textsc{2d},\, s}(x,\ y)$, by first solving 
\begin{equation}
   \dfrac{1}{\Delta y_m} \int_{y_{m -\frac{1}{2}}}^{y_{m +\frac{1}{2}}} 
   \psi_q(\eta(y))\,\hat{\vec{w}}^{\textsc{2d},\, s}_{ijk,\,pq}\,\de{y} = 
   \hat{\vec{w}}^{\textsc{1d}}_{imk,\,p},\quad\forall\,\Omega_{imk} \in S_{ijk}^{y,\, s},
\end{equation}
for each degree of freedom $\hat{\vec{w}}^{\textsc{1d}}_{imk,\,p}$ and then carrying out the
nonlinear selection as in the first pass. Analogously, in the third space dimension, conservation
over each element of the stencil gives
\begin{equation}
   \dfrac{1}{\Delta z_m} \int_{z_{m -\frac{1}{2}}}^{z_{m +\frac{1}{2}}} 
   \psi_r(\zeta(z))\,\hat{\vec{w}}^{\textsc{3d},\, s}_{ijk,\,pqr}\,\de{z} = 
   \hat{\vec{w}}^{\textsc{2d}}_{ijm,\,pq},\quad\forall\,\Omega_{ijm} \in S_{ijk}^{z,\, s},
\end{equation}
to be solved once for each degree of freedom $\hat{\vec{w}}^{\textsc{2d}}_{ijm,\,pq}$, 
and finally one obtains the sought three-dimensional weighted essentially-non-oscillatory polynomial after 
the nonlinear combination of the individual stencil polynomials has been applied.

\subsection{Reconstruction in primitive variables} 
\label{sec:primrec}
In this work, we employ a primitive variable reconstruction in order to better treat some of the
peculiar issues that are typically encountered in the numerical solution of multiphase flow models,
namely the presence of a complex, volume fraction-dependent equation of state and/or other issues
due to different material interfaces evolving separately, already reported in \cite{abgrall2001} and
\cite{larouturou1991}: this separation between interfaces might give rise to nonphysical
discontinuities in the velocity and density fields as well as positivity violations in the mass
fraction. 
The use of a primitive variable reconstruction for the TVD second order MUSCL-Hancock scheme, or the
\aderweno \pnpm{0}{2} scheme used in the troubled elements as subcell limiter schemes was found to
significantly mitigate these problems.

The primitive variable reconstruction procedure used for the \aderweno \pnpm{0}{2} limiter was
introduced in \cite{zanotti2015}, along with a predictor step formulated in terms of the primitive
form of the governing equations. The reconstruction is performed as follows: first, a conservative
polynomial WENO reconstruction is computed and the polynomials obtained from this step are
\emph{evaluated} at the cell centres so to obtain high order accurate point values for the conserved
variables. Then one can convert the point values of conserved variables to primitive variables and
perform a second WENO reconstruction to achieve a high order polynomial reconstruction of the
primitive data. This second reconstruction step repeats the same steps described in 
Section~\ref{sec:reconstruction}, with the difference that the reconstruction equations are not based
on directly enforcing conservation on a stencil, but rather they are obtained by requiring that the primitive
variable reconstruction polynomials interpolate the cell centre value where the conversion from 
conserved to primitive variables has taken place. For alternative high order WENO reconstructions in 
primitive variables, see \cite{BalsaraKim,PuppoRussoPrim}.  

\subsection{One-step, fully discrete, explicit update formulas}
We consider a general nonconservative hyperbolic system written as
\begin{equation} \label{eq:generalpde}
    \partial_t\vec{Q} + \nabla\cdot\vec{F}(\vec{Q}) + \vec{B}(\vec{Q}) \cdot \nabla\vec{Q} = \vec{S}(\vec{Q}).
\end{equation}
in a space time control volume $\Omega_{ijk}\times[t_n,\ t_{n+1}]$; we then define the differential volume
element $\de{\vec{x}} = \de{x}\,\de{y}\,\de{z}$ for compactly writing integrals over the control
volume {$\Omega_{ijk}$} and the surface element $\de{\vec{s}}$ for compactly writing integrals over its
boundary $\partial\Omega_{ijk}$. By multiplying each term of the PDE \eqref{eq:generalpde} with a test
function $\phi_p$, formally integrating over the space-time element and applying Gauss's theorem for
integrating the divergence of fluxes in space, we have a weak formulation
\begin{equation}
\begin{aligned}
    \intelemt \intelemome \phi_p\,\partial_t\vec{Q}\,\de{\vec{x}}\,\de{t} + 
    \intelemt \intelemdome \phi_p\,\vec{F}(\vec{Q})\cdot \uvec{n}\,\de{\vec{s}}\,\de{t}\ 
    -\intelemt \intelemome \nabla\phi_p \cdot \vec{F}(\vec{Q})\,\de{\vec{x}}\,\de{t}\ + &\\
    + \intelemt \intelemome \phi_p\,\left[\vec{B}(\vec{Q})\cdot \nabla\vec{Q} - 
    \vec{S}(\vec{Q})\right]\,\de{\vec{x}}\,\de{t} = \vec{0}, &
\end{aligned}
\end{equation}
with $\uvec{n}$ defined as the outward unit normal vector on the element boundary. Then, by
substituting the sought polynomial solution $\vec{u}_h(\vec{x},t^n) =
\phi_m(\vec{\xi})\, \hat{\vec{u}}_m^n$, 
as well as the polynomials $\vec{q}_h(\vec{x},t) = \theta_m(\vec{\xi},\tau)\, \hat{\vec{q}}_m$ 
obtained from the local space-time predictor detailed in the next section, we have the fully-discrete
one-step update formula
\begin{equation} \label{eq:dgupdate}
\begin{aligned}
    \left( \, \intelemome\phi_p\,\phi_q \de{\vec{x}} \right)\,\left( \hat{\vec{u}}_q^{n+1} 
    -\hat{\vec{u}}_q^{n} \right) + 
    \intelemt\intelemdome\phi_p\,\vec{F}_{\textsc{rp}}\left(\qhm,\ \qhp\right)\cdot \uvec{n}\,\de{\vec{s}}\,\de{t} 
    \intelemt\intelemdome\phi_p\,\vec{D}_\Psi\left(\qhm,\ \qhp\right)  \cdot \uvec{n}\,\de{\vec{s}}\,\de{t}\ +&\\
    -\intelemt\intelemiome\nabla\phi_p \cdot \vec{F}(\vec{q}_h)\,\de{\vec{x}}\,\de{t}
    +\intelemt\intelemiome\phi_p\,\left[\vec{B}(\vec{q}_h) \cdot \nabla\vec{q}_h - 
    \vec{S}(\vec{q}_h)\right]\,\de{\vec{x}}\,\de{t} = \vec{0},&
\end{aligned}
\end{equation}
where we denoted with $\vec{F}_{\textsc{rp}}$ the generic numerical flux function, that would be, for this
work, the Rusanov flux \eqref{eq:frus} or the HLL flux \eqref{eq:fhll}, but also other approximate Riemann 
solvers could be used, such as the generalized Osher and HLLEM methods forwarded in \cite{OsherNC,HLLEMNC}. 
Analogously, we define the path-conservative fluctuation term as 
\begin{equation} \label{eq:fluctuations}
   \vec{D}_\Psi\left(\qhm,\ \qhp\right) \cdot \uvec{n} = 
   \omega\,\int_0^1\vec{B}\left[\Psi\left(\qhm,\ \qhp,\ s\right)\right] \cdot \uvec{n}\,\de{s}\,\left(\qhp - \qhm\right), 
   \quad \text{with}\quad \Psi\left(\qhm,\ \qhp,\ s\right) = \qhm + s\,(\qhp - \qhm),
\end{equation}
in which $\Psi\left(\qhm,\ \qhp,\ s\right)$ is a simple segment path function connecting the left
and right states, and the path integral can be computed with a three-point Gauss--Legendre
quadrature regardless of the order of the scheme.
The coefficients $\omega$ must be chosen so to enforce the consistency condition
\cite{castro2006, pares2006}
\begin{equation} \label{eq:consistencyconditiondg}
\begin{aligned}
    \vec{D}_\Psi\left(\qhm,\ \qhp\right) \cdot \uvec{n} & - \vec{D}_\Psi\left(\qhp,\ \qhm\right) \cdot \uvec{n} = 
    \int_0^1\vec{B}\left[\Psi\left(\qhm,\ \qhp,\ s\right)\right] \cdot \uvec{n}\,\pd{\Psi}{s}\,\de{s}\,
\end{aligned}
\end{equation}
and simple expressions are provided in Sections~\ref{sec:hll}~and~\ref{sec:rus} for the HLL and
Rusanov fluxes.
The inversion of the mass matrix integrating the products $\phi_p\,\phi_q$ is trivial, as the choice
of basis yields an orthogonal basis and thus a diagonal mass matrix. The volume integrals appearing in \eqref{eq:dgupdate} may 
be directly evaluated by Gauss--Legendre quadrature using the nodes on which the degrees of freedom
of the space-time predictor solution $\vec{q}_h$ are defined, while for face integrals one has to
extrapolate $\qhm$ and $\qhp$ from two adjacent cells onto the Gauss--Legendre quadrature nodes at a
face, then evaluate the two-state numerical fluxes at each one of the quadrature nodes, and finally
operate the weighted sum of all the numerical fluxes.

Since numerical flux functions can be in principle computationally quite expensive, an attractive
alternative choice for the integration of fluxes at space-time cell boundaries, with respect to the
tensor-product quadrature rule, is the following: during the space-time predictor step, automatically
a polynomial approximation of the physical fluxes $\vec{f}_h = \vec{f}{(\vec{q}_h)}$ is computed
within each cell. When performing the extrapolation of $\qhm$ and $\qhp$ to the space-time boundaries, 
one may also directly extrapolate the approximation of the physical fluxes to the boundaries, obtaining 
thus at each space-time cell boundary $\fhm$ and $\fhp$. 

Here we denoted with $\vec{f}$ the projection of the physical flux tensor $\vec{F}$ on one
of the three canonical basis vectors indicating the orientation of the face-normal onto which the
flux is to be extrapolated, that is $\vec{f} = \vec{F}\cdot\uvec{e}_x$, $\vec{g} =
\vec{F}\cdot\uvec{e}_y$, or $\vec{h} = \vec{F}\cdot\uvec{e}_z$ in the first, second, or third direction, 
respectively.

Then one can treat $\qhm$, $\qhp$, $\fhm$, and $\fhp$ as four independent variables and recognise
that the numerical fluxes employed in the present work can be seen, if wavespeed estimates are
considered fixed, as split into a centred part (solely function of $\fhm$ and $\fhp$) and a
diffusive part (function of $\qhm$ and $\qhp$). Moreover, such a four-variable numerical flux with
fixed wavespeed estimates is linear in its arguments and in order to exploit this property, the
coefficients may be evaluated only once at a space-time-face-averaged state and employed for all
space-time face integration points. Thanks to the simple choice of a linear segment path, one can
apply the same approach to the computation of the path integral of nonconservative products, and
compute the average nonconservative product coefficient matrix 
\begin{equation} \label{eq:averagebmatrix}
    \vec{B}_{\vec{\Psi}} = \int_0^1\vec{B}\left[\Psi\left(\qhm,\ \qhp,\ 
        s\right)\right] \cdot \uvec{n}\,\de{s}
\end{equation}
only once, integrating between the averaged
states at the two faces, then multiplying \eqref{eq:averagebmatrix} by the two weights 
$\omega_{LR}$ and $\omega_{RL}$
and by the space-time-face average jump between $\qhm$ and $\qhp$, yielding 
$\vec{D}_\Psi\left(\qhm,\ \qhp\right) \cdot \uvec{n}$ 
and $\vec{D}_\Psi\left(\qhp,\ \qhm\right) \cdot \uvec{n}$, respectively. 

This means that only one nonlinear computation of the wavespeed estimates and other nonlinearities
in the Riemann solver has to be performed (with the face-averaged state of $\vec{q}_h$), while the
central part of the flux can be integrated directly, as well as the jump term in conserved variables. An
added benefit of this approach is that the scheme need not to retain information regarding the
space-time degrees of freedom of the predictor solution, making it possible and easy to implement
low-storage schemes that are of uniform arbitrary high order in space and time.

Finally, in order to guarantee stability of the explicit timestepping, 
in this work we restrict the timestep size by 
\begin{equation} \label{eq:stability}
    \Delta t = \textsc{CFL} \,\frac{k_{N}\min{(\Delta x,\ \Delta y,\ \Delta z)}}{d\,\lambda_\up{max}},
\end{equation}
with $N$ being the degree of the piecewise polynomial data representation, $d$ the number of space dimensions,
and $\lambda_\up{max}$ the maximum absolute value of all eigenvalues found in the domain (more
specifically, searching over all the quadrature nodes, i.e. where the degrees of freedom of
the nodal basis are collocated). With $\textsc{CFL} \leq 1$ we denote a Courant-type number that is
typically chosen as $\textsc{CFL} = 0.9$ for all the simulations presented in this work. The function 
$k_{N}$ was defined by numerical Von Neumann stability analysis in \cite{dumbser2008b} for polynomials  
of degree up to four, while for higher values of $N$, we refer to an experimental determination
based on numerical tests with linear advection. 

The first five values of $k_{N}$ are $1.0$, $0.33$, $0.17$, $0.10$, and $0.069$ starting from \finvol ($N=0$)
up to $N=4$ (fifth order \pnpm{4}{4} ADER-DG scheme), while for $5 \leq N \leq 9$ we choose, $k_{N}$ from
the vector $k_{5-9} = \transpose{[0.045,\ 0.038,\ 0.03,\ 0.02,\ 0.01]}$. We conclude by 
pointing out that condition \eqref{eq:stability} follows the same behaviour of the common $\Delta
t_\up{max} \sim 1/(2\,N + 1)$ hyperbola for RKDG methods, but is slightly more restrictive. 

\subsection{Space-time \disgal predictor}
We now describe the procedure to obtain the space-time predictor polynomials, which are defined as 
\begin{equation} \label{eq:qhansatz}
    \vec{q}_h(\xi,\eta,\zeta,\tau) =
    \psix(\xi)\,\psiy(\eta)\,\psiz(\zeta)\,\psit(\tau)\, \hat{\vec{q}}_{m_x\,m_y\,m_z\,m_t}
    = \theta_m(\xi,\eta,\zeta,\tau)\, \hat{\vec{q}}_m,  
\end{equation}

again formally allowing referencing to the components of $\theta$ with mono-indexing or
multi-indexing. The first step for the local time evolution starting from the polynomial data
$\vec{w}_h(\vec{x},\ t^n)$ is to write the governing PDE \eqref{eq:model} in a weak integral form in space and time as 
\begin{equation}
\begin{aligned}
    \intelemt \intelemomega \theta_p\,\partial_t\vec{q}_h\,\de{\vec{x}}\,\de{t} + 
    \intelemt \intelemomega \theta_p\,\nabla\cdot\vec{F}(\vec{q}_h)\,\de{\vec{x}}\,\de{t}
    + \intelemt \intelemomega \theta_p\,\left[\vec{B}(\vec{q}_h) \cdot \nabla\vec{q}_h - \vec{S}(\vec{q}_h)\right]\,\de{\vec{x}}\,\de{t} = 0,
\end{aligned}
\end{equation}
and then integrating by parts in time the first term (and by upwinding in time the value
of $\vec{q}_h(\vec{x},\ t^n)$ from the reconstruction polynomial $\vec{w}_h(\vec{x},t^n)$), we can write
\begin{equation} \label{eq:predictorsystemquick}
\begin{aligned}
    \phantom{+} \intelemomega \theta_p(\vec{x},\ t^{n+1})\,\vec{q}_h(\vec{x},\ t^{n+1})\,\de{\vec{x}}
    - \intelemomega \theta_p(\vec{x},\ t^n)\,\vec{w}_h(\vec{x},\ t^n)\,\de{\vec{x}}
    - \intelemt \intelemomega \partial_t\theta_p\,\vec{q}_h\,\de{\vec{x}}\,\de{t}\ +&\\
    +\intelemt \intelemomega \theta_p\,\nabla\cdot\vec{F}(\vec{q}_h)\,\de{\vec{x}}\,\de{t} +
    \intelemt \intelemomega \theta_p\,\vec{B}(\vec{q}_h)\cdot\nabla\vec{q}_h\,\de{\vec{x}}\,\de{t} = 
       \intelemt \intelemomega \theta_p\,\vec{S}(\vec{q}_h)\,\de{\vec{x}}\,\de{t}.&
\end{aligned}
\end{equation}
By then substituting the ansatz \eqref{eq:qhansatz} in \eqref{eq:predictorsystemquick} one obtains a
system of ${(M+1)}^{d+1}$ nonlinear algebraic equations which one can solve by means of a discrete
Picard iteration with appropriate initial guess, as discussed in \cite{dumbser2008a,Dumbser2018a,Busto2019}. 

\subsubsection{Predictor step in primitive variables}
In conjunction with the primitive variable WENO reconstruction described in 
Section~\ref{sec:reconstruction}, as well as for pure ADER \disgal schemes, for which 
primitive variable polynomials can be obtained by simply evaluating the primitive state vector
in correspondence of each quadrature node location (nodal degree of freedom), 
we also carry out the local time evolution procedure with a primitive variable formulation, 
as per the methodology introduced in \cite{zanotti2015}.
This variant of the local space-time predictor step is based on a primitive variable 
version of the governing equations, which directly evolves the primitive state vector $\vec{V}$,
uses only gradients of the primitive variables $\nabla\vec{V}$ and is recovered by applying
the chain rule to the governing equations in the form \eqref{eq:generalpde} to obtain
\begin{equation}
    \pd{\vec{V}}{t} + {\left(\pd{\vec{Q}}{\vec{V}}\right)}^{-1}\left(\pd{\vec{F}}{\vec{V}} + 
    \vec{B}\cdot\pd{\vec{Q}}{\vec{V}}\right)\cdot\nabla\vec{V} = 
    {\left(\pd{\vec{Q}}{\vec{V}}\right)}^{-1}\vec{S},
\end{equation}
We now assign the notation $\vec{w}_h^\ast$ to represent
the discrete reconstruction data \emph{in primitive variables}, obtained either by the primitive
variable WENO reconstruction, or by a straightforward conversion of the nodal degrees of 
freedom for ADER-DG schemes, and define $\vec{v}_h$ to be 
the discrete space-time predictor solution \emph{in primitive variables},
we can write a weak form of the governing equations as
\begin{equation}
    \intelemt \intelemomega \theta_l\,\partial_t\vec{v}_h\,\de{\vec{x}}\,\de{t}
    + \intelemt \intelemomega \theta_l\,
    {\left(\pd{\vec{Q}}{\vec{V}}\right)}^{-1}
    \left(
    \pd{\vec{F}}{\vec{V}} + 
    \vec{B}\cdot\pd{\vec{Q}}{\vec{V}}
    \right)
    \cdot\nabla\vec{v}_h\,\de{\vec{x}}\,\de{t} =
    \intelemt \intelemomega \theta_l\,{\left(\pd{\vec{Q}}{\vec{V}}\right)}^{-1}
    \vec{S}\,\de{\vec{x}}\,\de{t},
\end{equation}
and again integrating by parts in time one obtains a nonlinear algebraic system of $(M+1)^{d+1}$ equations
\begin{equation}
\begin{aligned}
\intelemomega \theta_l(\vec{x},\ t^{n+1})\,\vec{v}_h(\vec{x},\ t^{n+1})\,\de{\vec{x}} -  
  \intelemomega \theta_l(\vec{x},\ t^{n})\,\vec{w}_h^\ast(\vec{x},t^n)\,\de{\vec{x}}
    - \intelemt \intelemomega \partial_t\theta_l\,\vec{v}_h\,\de{\vec{x}}\ +&\\
    + \intelemt \intelemomega \theta_l\,
{\left(\pd{\vec{Q}}{\vec{V}}\right)}^{-1}
    \left(
    \pd{\vec{F}}{\vec{V}} + 
    \vec{B}\cdot\pd{\vec{Q}}{\vec{V}}
    \right)
    \cdot\nabla\vec{v}_h\,\de{\vec{x}} 
     = \intelemt \intelemomega \theta_l\,{\left(\pd{\vec{Q}}{\vec{V}}\right)}^{-1}\vec{S}\,\de{\vec{x}},&
\end{aligned}
\end{equation}
again to be solved via a discrete Picard iteration \cite{dumbser2008a} and then extrapolated to 
the cell boundaries to compute the numerical fluxes and fluctuations, as well as the volume integrals
of the explicit update formulas \eqref{eq:dgupdate}.

\subsubsection{The path-conservative Harten--Lax--Van Leer flux} \label{sec:hll}

We denote with $\vec{f}$, $\vec{g}$ and $\vec{h}$ the relevant projections of the physical flux tensor 
$\vec{F}$ onto the Cartesian coordinate directions, i.e. $\vec{f} = \vec{F}\cdot\uvec{e}_x$, 
$\vec{g} = \vec{F}\cdot\uvec{e}_y$ and $\vec{h} = \vec{F}\cdot\uvec{e}_z$, 
according to the direction normal to the face/edge along which the solution
of the Riemann problem is sought. With reference to two generic input states $\qhl$ and $\qhr$, the
HLL flux reads as follows, 
\begin{equation} \label{eq:fhll}
    \vec{F}_{\textsc{rp}}^{\textsc{hll}}\left(\qhl,\ \qhr\right) = 
    \dfrac{\shllr\,\vec{f}\left(\qhl\right) - \shlll\,\vec{f}\left(\qhr\right)}{\shllr - \shlll} + 
    \dfrac{S_\textsc{r}\,\shlll}{\shllr - \shlll}\,\left(\qhr - \qhl\right),
\end{equation}
and we give the estimates of the minimum and maximum wave speeds as
\begin{equation}
    \shlll = \min{\left[0,\ \lambda_\up{min}\left(\qhl\right),\ \lambda_\up{min}
        \left(\overline{\vec{q}}\right)\right]},\quad 
    \shllr = \max{\left[0,\ \lambda_\up{max}\left(\qhr\right),\ \lambda_\up{max}
        \left(\overline{\vec{q}}\right)\right]},\quad \text{with}
     \quad \overline{\vec{q}} = \dfrac{1}{2}\,\left(\qhl + \qhr\right),
\end{equation}
where $\lambda_\up{min}(\vec{q})$ and $\lambda_\up{max}(\vec{q})$ are functions computing, 
respectively, the minimum and the maximum eigenvalue of the system of equations for a given vector of
conserved variables $\vec{q}$.
Given an outward  unit normal vector $\uvec{n}$ such that the scalar product with the positive
generic direction vector $\uvec{e}_k$ can be either positive or negative unity, upwinding of the
nonconservative terms is accounted for by setting in Eq.~\eqref{eq:fluctuations}
\begin{equation}
    \omega = \frac{1}{2}\,\left(1+\frac{\shlll + \shllr}{\shllr - \shlll}\,\uvec{e}_k\cdot\uvec{n}\right).
\end{equation}
This means that, for a given face with jump states $\qhl$ and $\qhr$, in a Cartesian setting, we
will compute two weights $\vec{\omega}$ to associate with the two fluctuation terms, one associated with
a positive unit normal, one associated with a negative unit normal.  
\subsubsection{The Rusanov flux} \label{sec:rus}
The Rusanov flux is obtained from the HLL flux under the assumption that 
$\shlll = -S_\up{max}$ and $\shllr = S_\up{max}$ and can be written as  
\begin{equation} \label{eq:frus}
    \vec{F}_{\textsc{rp}}^{\textsc{rus}}\left(\qhl,\ \qhr\right) = \dfrac{1}{2}\left[\vec{f}\left(\qhl\right) + 
    \vec{f}\left(\qhr\right)\right] - \dfrac{1}{2}\,S_\up{max}\,\left(\qhr - \qhl\right).
\end{equation}
This flux only requires the computation of a single wave speed estimate which is 
 \begin{equation}
     S_\up{max} = \max{\left(\abs{\lambda_\up{min}\left(\qhl\right)},\ \abs{\lambda_\up{min}\left(\qhr\right)},\ 
     \abs{\lambda_\up{max}\left(\qhl\right)},\ \abs{\lambda_\up{max}\left(\qhr\right)}\right]};
\end{equation}
as for the conservative part, the nonconservative fluctuations associated with the Rusanov flux do
not account for upwinding and therefore, enforcing the generalized Rankine--Hugoniot consistency
condition \eqref{eq:consistencyconditiondg} \cite{castro2006, pares2006} we set $\omega = 1/2$.

\subsection{A posteriori subcell limiting (MOOD)} The \aposteriori subcell limiting approach \cite{dumbser2014} consists in 
\emph{first} computing a \emph{candidate solution} $\vec{u}_h^\ast(\vec{x},\ t^{n+1})$ from the
ADER-DG scheme, without applying any precaution for limiting spurious oscillations that are typical of high
order linear methods, and subsequently verifying the admissibility of such a solution by means of a
relaxed discrete maximum principle and other features that might characterise the solution as
locally not valid, such as violations of the positivity of density and pressure or floating-point
arithmetic exceptions. This novel \textit{a posteriori} limiting strategy for DG schemes follows the 
ideas of the MOOD approach, which was forwarded by Clain and Loub\`ere \etal\  in 
\cite{clain2011,diot2013,diot2014,loubere2014} within the \finvol framework. 
The \textit{relaxed} discrete maximum principle (DMP) is satisfied if, for all conserved (or primitive) 
variables, the solution is such that
\begin{equation}
    \min_{\vec{y}\in\mathcal{N}_{i\,j\,k}}\left[\vec{u}_h\left(\vec{y},\ t^n\right)\right] - \vec{\delta} \leq 
    \vec{u}_h^\ast(\vec{x},\ t_{n+1}) 
    \leq 
    \max_{\vec{y}\in\mathcal{N}_{i\,j\,k}}\left[\vec{u}_h\left(\vec{y},\ t^n\right)\right] + 
    \vec{\delta}, \quad \forall \vec{x} \in \Omega_{i\,j\,k},
\end{equation}
with 
\begin{equation}
\label{eq:dmptol}
    \vec{\delta} = \max\left(\epsilon_0,\ 
    \epsilon_1\,\left\{
    \max_{\vec{y}\in\mathcal{N}_{i\,j\,k}}\left[\vec{u}_h\left(\vec{y},\ t^n\right)\right] - 
    \min_{\vec{y}\in\mathcal{N}_{i\,j\,k}}\left[\vec{u}_h\left(\vec{y},\ t^n\right)\right]\right\},\ 
    \epsilon_2\,\min_{\vec{y}\in\mathcal{N}_{i\,j\,k}}\left[\abs{\vec{u}_h\left(\vec{y},\ t^n\right)}\right]\right).
\end{equation}
The three small constant parameters in Eq.~\eqref{eq:dmptol} are set as $\epsilon_0 = 10^{-4}$,
$\epsilon_1 = 10^{-3}$, and $\epsilon_2 = \epsilon_0/1000$, the last being intended to prevent
excessively restrictive requirements on the oscillations of variables which have typical magnitude
much larger than unity: by choosing $\epsilon_2 = \epsilon_0/1000$, we are prescribing that if, for
a given variable, all the values in $\mathcal{N}_{i\,j\,k}$ have absolute magnitude larger than
$1000$, then the dimensionless floor value of $\vec{\delta}$ for that variable will be comparable
to that of unit-scaled variables. This will typically be the case for liquid density or internal
energy, which otherwise might trigger the \aposteriori limiter unnecessarily. All of the cells where
the admissibility criteria are not satisfied are marked and the data from the previous timestep is
projected on a finer local \finvol subgrid; if a given cell was already marked during the
previous timestep, such data is recovered from the subcell-average representation directly, while
one must compute the local subcell averages of the polynomial data if the limiter state at the
previous timestep is not available. Then the solution is recomputed with a more robust \finvol
scheme and new polynomial data for the original element is reconstructed by solving an
overdetermined linear system of conservative reconstruction equations.

\section{Test problems}\label{sec:testproblems} 

In this section, we present the results obtained by applying the ADER family of methods to all 
three variants of the \sch model: the original weakly hyperbolic formulation \eqref{eq:gavrilyuk}, 
the hyperbolic non-conservative symmetrizable Godunov--Powell form \eqref{eq:model}, and the 
hyperbolic GLM curl-cleaning formulation \eqref{eq:gavrilyuk.glm}. As we have already mentioned, other 
variants of the model were also tested (\textit{e.g.} Godunov--Powell + GLM, or GLM with extra terms in the 
energy equation \eqref{eq:glmenergync}) but these variants show very similar results, at least for 
the considered test cases, to the first three formulations and therefore, are not presented here. 

\subsection{Numerical convergence results}
As a first benchmark, we conduct a numerical convergence study on a smooth test problem, for which we have 
derived the exact solution in Section~\ref{sec:exactsol}. 
The test is very similar to the well known isentropic vortex
advection problem \cite{shu1997} for the Euler equations of gasdynamics: a steady state solution is
initialised at time $t = 0$ in a uniform flow field $\vec{u} = \transpose{[u_0,\ v_0,\ 0]}$ and
evolved with periodic boundary conditions on a rectangular domain of edge lengths $L_x$ and $L_y$.
Due to the Galilean invariance of the governing equations, the exact solution is obtained by
transporting the initial condition with the uniform flow speed.
\begin{figure}[!b]
    \includegraphics[draft=false, scale=\figurescalefactor]{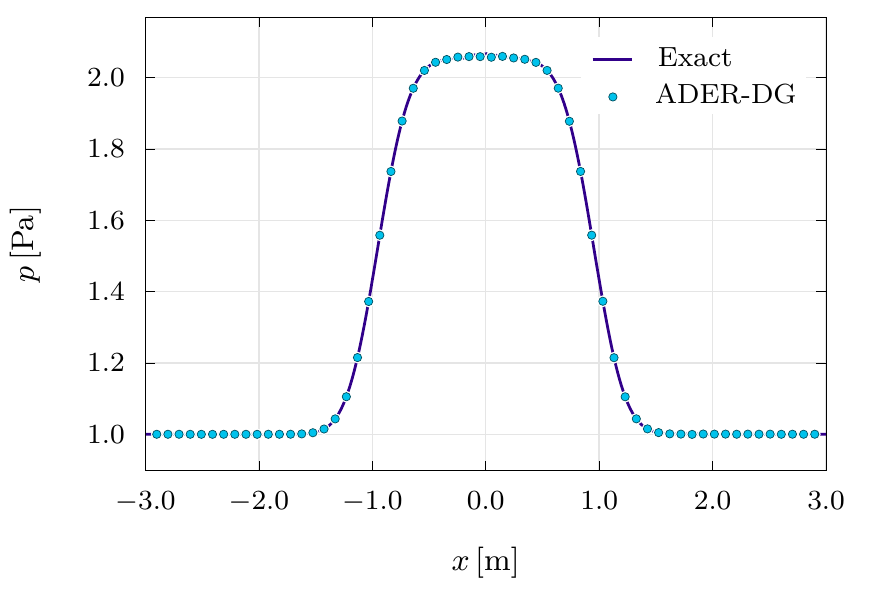}%
    \includegraphics[draft=false, scale=\figurescalefactor]{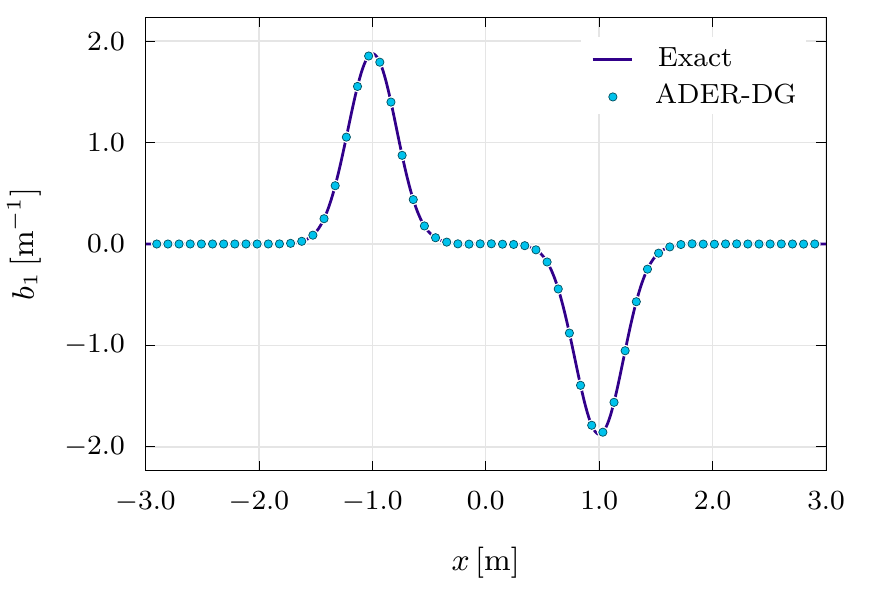}\\[2mm]
    \includegraphics[draft=false, scale=\figurescalefactor]{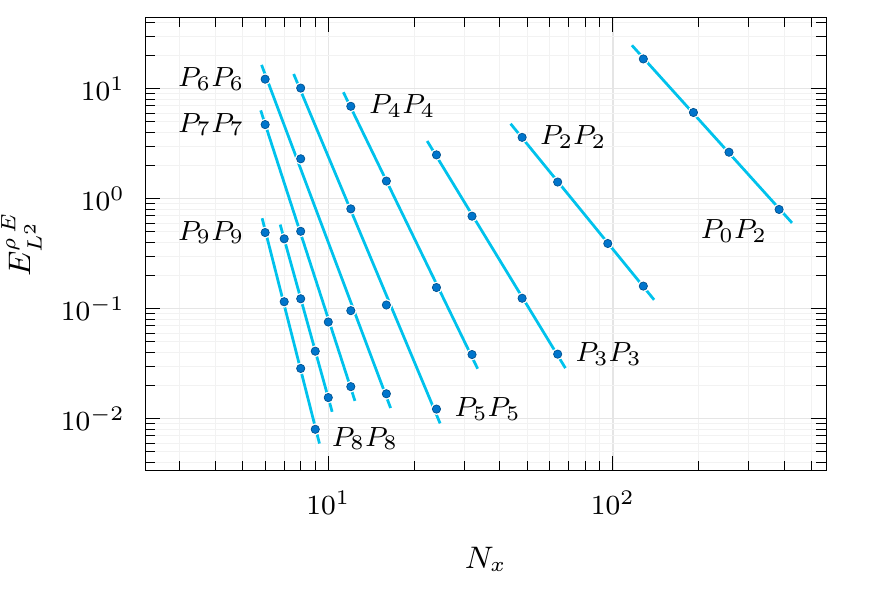}%
    \includegraphics[draft=false, scale=\figurescalefactor]{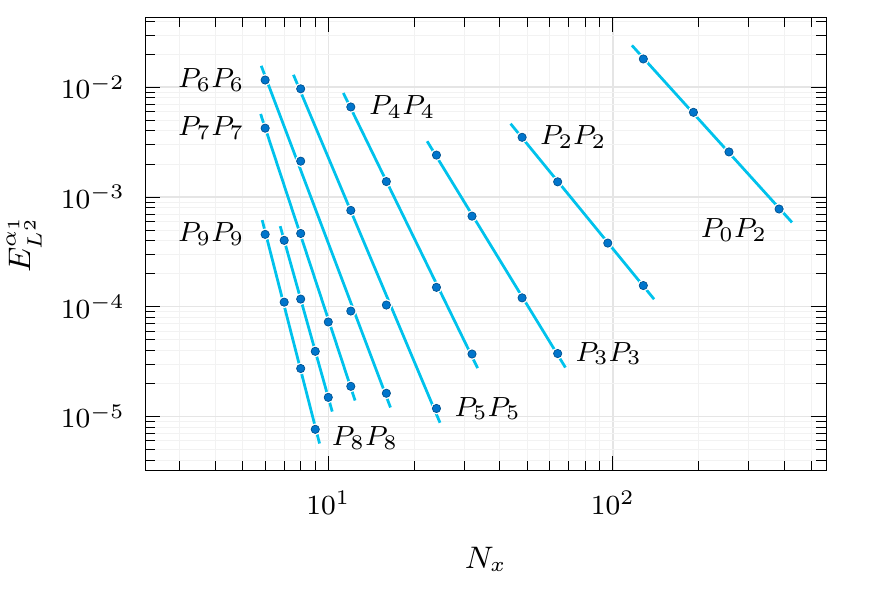}
    \caption{Numerical convergence results. In the top row, one-dimensional cuts (60 equidistant samples
    along the $x$ axis) of the pressure $p$ and the $x$ component of the interface field $b_1$ are presented at
    time $t=12\,\up{s}$ (at the end of the second advection cycle), computed with a $\pnpm{9}{9}$ \aderdg 
    scheme on a uniform grid of $6^2$ elements. In the bottom row we show logarithmic least-squares regression lines of
    the $L^2$ error norms of total energy density $\rho\,E$ and of the liquid volume fraction $\alpha_1$
    for the \aderdg \pnpm{N}{N} schemes of orders 3 to 10 and for the third order \pnpm{0}{2} \aderweno
    \finvol scheme.}
    \label{fig:droplet-convergence}
\end{figure}
\begin{table}[!t]
    \caption{Numerical convergence results regarding the ADER-DG $\pnpm{N}{N}$ schemes of nominal orders of
    accuracy 3 to 10 and the \aderweno $\pnpm{0}{2}$ \finvol scheme for all conserved variables. The
    values reported in the Table are computed from a logarithmic least-square fit of the $L^2$ error
    norms as shown in Figures~\ref{fig:droplet-convergence}.}
    \label{tab:convergencetable}
    \begin{tabularx}{\textwidth}{rlRRRRRRRR}
    \toprule
         &      &   $\alpha_1\,\rho_1$ & $\alpha_2\,\rho_2$ & $\rho\,u$ & $\rho\,v$ & $\rho\,E$ & $\alpha_1$ & $b_1$ & $b_2$ \\
    \midrule
$\pnpm{0}{2}$ & $\mathcal{O}_{L^1}$      &  2.9 &  3.1 &  2.9 &  2.9 &  2.9 &  2.9 &  2.8 &  2.8\\
              & $\mathcal{O}_{L^2}$      &  2.9 &  3.0 &  2.9 &  2.9 &  2.9 &  2.9 &  2.8 &  2.8\\
              & $\mathcal{O}_{L^\infty}$ &  2.8 &  3.0 &  2.8 &  2.8 &  2.8 &  2.8 &  2.7 &  2.7\\[2mm]
$\pnpm{2}{2}$ & $\mathcal{O}_{L^1}$      &  3.2 &  3.2 &  3.2 &  3.2 &  3.2 &  3.2 &  3.2 &  3.2\\
              & $\mathcal{O}_{L^2}$      &  3.2 &  3.2 &  3.2 &  3.2 &  3.2 &  3.2 &  3.2 &  3.2\\
              & $\mathcal{O}_{L^\infty}$ &  3.1 &  3.1 &  3.1 &  3.1 &  3.1 &  3.1 &  3.1 &  3.1\\[2mm]
$\pnpm{3}{3}$ & $\mathcal{O}_{L^1}$      &  4.3 &  4.3 &  4.3 &  4.3 &  4.3 &  4.3 &  4.3 &  4.3\\
              & $\mathcal{O}_{L^2}$      &  4.3 &  4.3 &  4.3 &  4.3 &  4.3 &  4.2 &  4.3 &  4.3\\
              & $\mathcal{O}_{L^\infty}$ &  4.2 &  4.1 &  4.2 &  4.2 &  4.2 &  4.2 &  4.3 &  4.3\\[2mm]
$\pnpm{4}{4}$ & $\mathcal{O}_{L^1}$      &  5.5 &  6.1 &  5.5 &  5.5 &  5.5 &  5.5 &  5.5 &  5.5\\
              & $\mathcal{O}_{L^2}$      &  5.3 &  5.6 &  5.3 &  5.3 &  5.3 &  5.3 &  5.3 &  5.3\\
              & $\mathcal{O}_{L^\infty}$ &  5.1 &  5.3 &  5.1 &  5.1 &  5.0 &  5.1 &  5.0 &  5.0\\[2mm]
$\pnpm{5}{5}$ & $\mathcal{O}_{L^1}$      &  6.5 &  7.3 &  6.5 &  6.5 &  6.5 &  6.4 &  6.4 &  6.4\\
              & $\mathcal{O}_{L^2}$      &  6.2 &  6.8 &  6.2 &  6.2 &  6.2 &  6.2 &  6.1 &  6.1\\
              & $\mathcal{O}_{L^\infty}$ &  5.8 &  6.6 &  5.8 &  5.8 &  5.8 &  5.7 &  5.8 &  5.8\\[2mm]
$\pnpm{6}{6}$ & $\mathcal{O}_{L^1}$      &  7.2 &  8.6 &  7.2 &  7.2 &  7.2 &  7.2 &  7.1 &  7.1\\
              & $\mathcal{O}_{L^2}$      &  6.9 &  8.0 &  6.9 &  6.9 &  6.9 &  6.9 &  6.8 &  6.8\\
              & $\mathcal{O}_{L^\infty}$ &  6.6 &  8.1 &  6.6 &  6.6 &  6.6 &  6.5 &  6.4 &  6.3\\[2mm]
$\pnpm{7}{7}$ & $\mathcal{O}_{L^1}$      &  8.3 & 10.0 &  8.3 &  8.3 &  8.3 &  8.2 &  8.1 &  8.1\\
              & $\mathcal{O}_{L^2}$      &  8.0 &  9.4 &  8.0 &  8.0 &  8.0 &  7.9 &  7.7 &  7.7\\
              & $\mathcal{O}_{L^\infty}$ &  7.8 &  9.4 &  7.8 &  7.8 &  7.8 &  7.4 &  7.1 &  7.1\\[2mm]
$\pnpm{8}{8}$ & $\mathcal{O}_{L^1}$      &  9.7 & 11.1 &  9.7 &  9.7 &  9.7 &  9.6 &  9.6 &  9.6\\
              & $\mathcal{O}_{L^2}$      &  9.3 & 10.2 &  9.3 &  9.3 &  9.3 &  9.2 &  9.2 &  9.2\\
              & $\mathcal{O}_{L^\infty}$ &  8.7 &  9.1 &  8.7 &  8.7 &  8.7 &  8.8 &  8.4 &  8.4\\[2mm]
$\pnpm{9}{9}$ & $\mathcal{O}_{L^1}$      & 10.7 & 11.7 & 10.7 & 10.7 & 10.7 & 10.6 & 10.6 & 10.6\\
              & $\mathcal{O}_{L^2}$      & 10.2 & 11.0 & 10.2 & 10.2 & 10.2 & 10.1 & 10.1 & 10.1\\
              & $\mathcal{O}_{L^\infty}$ &  9.9 & 11.2 &  9.9 &  9.9 &  9.9 & 10.0 &  9.7 &  9.7\\
    \bottomrule
    \end{tabularx}
\end{table}
\begin{table}[!t]
    \caption{Numerical convergence results regarding the ADER-DG $\pnpm{N}{N}$ schemes of nominal orders of
    accuracy 3 to 10 and the ADER-WENO $\pnpm{0}{2}$ \finvol scheme for the liquid volume fraction
    $\alpha_1$. With $N_x$ we indicate the number of cells in one row of the Cartesian computational
    grid.}
    \label{tab:convergencetablealpha}
    \begin{tabularx}{\textwidth}{rrRRRRRR}
    \toprule
            & $N_x$ & $E_{L^1}$ & $E_{L^2}$ & $E_{L^\infty}$ & $\mathcal{O}_{L^1}$ & $\mathcal{O}_{L^2}$ & $\mathcal{O}_{L^\infty}$ \\
    \midrule
   $\pnpm{0}{2}$ & 128 & $4.32\!\times\!10^{-2}$ &  $1.81\!\times\!10^{-2}$ &  $1.42\!\times\!10^{-2}$ &     $-$ &     $-$ &     $-$ \\
                 & 192 & $1.41\!\times\!10^{-2}$ &  $5.90\!\times\!10^{-3}$ &  $4.76\!\times\!10^{-3}$ &  $2.76$ &  $2.76$ &  $2.70$ \\
                 & 256 & $6.10\!\times\!10^{-3}$ &  $2.57\!\times\!10^{-3}$ &  $2.08\!\times\!10^{-3}$ &  $2.92$ &  $2.89$ &  $2.87$ \\
                 & 384 & $1.83\!\times\!10^{-3}$ &  $7.75\!\times\!10^{-4}$ &  $6.29\!\times\!10^{-4}$ &  $2.97$ &  $2.96$ &  $2.95$ \\[2mm]
   $\pnpm{2}{2}$ &  48 & $7.91\!\times\!10^{-3}$ &  $3.49\!\times\!10^{-3}$ &  $3.24\!\times\!10^{-3}$ &     $-$ &     $-$ &     $-$ \\
                 &  64 & $3.07\!\times\!10^{-3}$ &  $1.37\!\times\!10^{-3}$ &  $1.32\!\times\!10^{-3}$ &  $3.29$ &  $3.25$ &  $3.13$ \\
                 &  96 & $8.47\!\times\!10^{-4}$ &  $3.79\!\times\!10^{-4}$ &  $3.77\!\times\!10^{-4}$ &  $3.17$ &  $3.17$ &  $3.09$ \\
                 & 128 & $3.47\!\times\!10^{-4}$ &  $1.55\!\times\!10^{-4}$ &  $1.56\!\times\!10^{-4}$ &  $3.10$ &  $3.10$ &  $3.07$ \\[2mm]
   $\pnpm{3}{3}$ &  24 & $5.50\!\times\!10^{-3}$ &  $2.40\!\times\!10^{-3}$ &  $2.35\!\times\!10^{-3}$ &     $-$ &     $-$ &     $-$ \\
                 &  32 & $1.46\!\times\!10^{-3}$ &  $6.67\!\times\!10^{-4}$ &  $7.04\!\times\!10^{-4}$ &  $4.61$ &  $4.46$ &  $4.19$ \\
                 &  48 & $2.62\!\times\!10^{-4}$ &  $1.20\!\times\!10^{-4}$ &  $1.27\!\times\!10^{-4}$ &  $4.24$ &  $4.23$ &  $4.22$ \\
                 &  64 & $8.18\!\times\!10^{-5}$ &  $3.73\!\times\!10^{-5}$ &  $3.87\!\times\!10^{-5}$ &  $4.05$ &  $4.06$ &  $4.14$ \\[2mm]
   $\pnpm{4}{4}$ &  12 & $1.76\!\times\!10^{-2}$ &  $6.60\!\times\!10^{-3}$ &  $5.52\!\times\!10^{-3}$ &     $-$ &     $-$ &     $-$ \\
                 &  16 & $3.24\!\times\!10^{-3}$ &  $1.38\!\times\!10^{-3}$ &  $1.48\!\times\!10^{-3}$ &  $5.89$ &  $5.44$ &  $4.59$ \\
                 &  24 & $3.29\!\times\!10^{-4}$ &  $1.50\!\times\!10^{-4}$ &  $1.71\!\times\!10^{-4}$ &  $5.64$ &  $5.48$ &  $5.32$ \\
                 &  32 & $8.15\!\times\!10^{-5}$ &  $3.69\!\times\!10^{-5}$ &  $4.03\!\times\!10^{-5}$ &  $4.86$ &  $4.86$ &  $5.02$ \\[2mm]
   $\pnpm{5}{5}$ &   8 & $2.87\!\times\!10^{-2}$ &  $9.66\!\times\!10^{-3}$ &  $6.66\!\times\!10^{-3}$ &     $-$ &     $-$ &     $-$ \\
                 &  12 & $1.88\!\times\!10^{-3}$ &  $7.53\!\times\!10^{-4}$ &  $7.87\!\times\!10^{-4}$ &  $6.72$ &  $6.29$ &  $5.27$ \\
                 &  16 & $2.35\!\times\!10^{-4}$ &  $1.03\!\times\!10^{-4}$ &  $1.32\!\times\!10^{-4}$ &  $7.23$ &  $6.91$ &  $6.20$ \\
                 &  24 & $2.59\!\times\!10^{-5}$ &  $1.18\!\times\!10^{-5}$ &  $1.36\!\times\!10^{-5}$ &  $5.43$ &  $5.36$ &  $5.62$ \\[2mm]
   $\pnpm{6}{6}$ &   6 & $3.62\!\times\!10^{-2}$ &  $1.16\!\times\!10^{-2}$ &  $8.98\!\times\!10^{-3}$ &     $-$ &     $-$ &     $-$ \\
                 &   8 & $5.97\!\times\!10^{-3}$ &  $2.12\!\times\!10^{-3}$ &  $1.74\!\times\!10^{-3}$ &  $6.27$ &  $5.92$ &  $5.71$ \\
                 &  12 & $2.26\!\times\!10^{-4}$ &  $9.10\!\times\!10^{-5}$ &  $1.07\!\times\!10^{-4}$ &  $8.07$ &  $7.76$ &  $6.88$ \\
                 &  16 & $3.65\!\times\!10^{-5}$ &  $1.62\!\times\!10^{-5}$ &  $1.62\!\times\!10^{-5}$ &  $6.34$ &  $6.00$ &  $6.55$ \\[2mm]
   $\pnpm{7}{7}$ &   6 & $1.24\!\times\!10^{-2}$ &  $4.23\!\times\!10^{-3}$ &  $3.48\!\times\!10^{-3}$ &     $-$ &     $-$ &     $-$ \\
                 &   8 & $1.17\!\times\!10^{-3}$ &  $4.64\!\times\!10^{-4}$ &  $5.41\!\times\!10^{-4}$ &  $8.19$ &  $7.68$ &  $6.47$ \\
                 &  10 & $1.69\!\times\!10^{-4}$ &  $7.24\!\times\!10^{-5}$ &  $1.01\!\times\!10^{-4}$ &  $8.69$ &  $8.33$ &  $7.54$ \\
                 &  12 & $4.22\!\times\!10^{-5}$ &  $1.87\!\times\!10^{-5}$ &  $2.07\!\times\!10^{-5}$ &  $7.60$ &  $7.41$ &  $8.68$ \\[2mm]
   $\pnpm{8}{8}$ &   7 & $9.88\!\times\!10^{-4}$ &  $4.01\!\times\!10^{-4}$ &  $4.32\!\times\!10^{-4}$ &     $-$ &     $-$ &     $-$ \\
                 &   8 & $2.69\!\times\!10^{-4}$ &  $1.17\!\times\!10^{-4}$ &  $1.42\!\times\!10^{-4}$ &  $9.73$ &  $9.23$ &  $8.35$ \\
                 &   9 & $8.71\!\times\!10^{-5}$ &  $3.91\!\times\!10^{-5}$ &  $5.80\!\times\!10^{-5}$ &  $9.59$ &  $9.31$ &  $7.60$ \\
                 &  10 & $3.20\!\times\!10^{-5}$ &  $1.49\!\times\!10^{-5}$ &  $1.79\!\times\!10^{-5}$ &  $9.50$ &  $9.17$ & $11.16$ \\[2mm]
   $\pnpm{9}{9}$ &   6 & $1.19\!\times\!10^{-3}$ &  $4.55\!\times\!10^{-4}$ &  $5.21\!\times\!10^{-4}$ &     $-$ &     $-$ &     $-$ \\
                 &   7 & $2.53\!\times\!10^{-4}$ &  $1.10\!\times\!10^{-4}$ &  $1.49\!\times\!10^{-4}$ & $10.02$ &  $9.24$ &  $8.13$ \\
                 &   8 & $5.97\!\times\!10^{-5}$ &  $2.72\!\times\!10^{-5}$ &  $3.86\!\times\!10^{-5}$ & $10.83$ & $10.43$ & $10.11$ \\
                 &   9 & $1.62\!\times\!10^{-5}$ &  $7.59\!\times\!10^{-6}$ &  $9.05\!\times\!10^{-6}$ & $11.07$ & $10.84$ & $12.31$ \\
    \bottomrule
    \end{tabularx}
\end{table}

\subsubsection{Problem setup}

The initial condition for the liquid volume fraction $\alpha_1$ is given according to the chosen
colour function profile, but bounding it between the two values {$\alpha_\up{min} = 0.01$} and 
{$\alpha_\up{max} = 0.99$}, so that we have
\begin{equation}
   \label{eq:convergencesetup-alpha}
   \alpha_1(r) = \alpha_\up{min} + \frac{\alpha_\up{max} - \alpha_{\up{min}}}{2}\,
   \erfc\left(\frac{\rs^\prime - 1}{\keps}\right).\\
\end{equation}
Since also the density fields should be transported with uniform velocity regardless of their value,
we decided not to impose a constant value for $\rho_1$ and $\rho_2$, but rather specify a more
interesting periodic two-dimensional wave configuration in the form
\begin{align}
   \label{eq:convergencesetup-rho}
    &\rho_1(x,\ y) = \overline{\rho_1} + \delta\,\overline{\rho_1}\,\sin{\left[
        \omega\,(2\,x + y)\right]}\,\cos{\left[\omega\,(x - 2\,y)\right]},\\
    &\rho_2(x,\ y) = \overline{\rho_2} + \delta\,\overline{\rho_2}\,\sin{\left[
        \omega\,(x - 2\,y)\right]}\,\cos{\left[\omega\,(2\,x + y)\right]}.
\end{align}

The numerical values for the test are {$u_0 = v_0 = 3\,\up{m\,s^{-1}}$}, {$R = 1\,\up{m}$},
{$k_\epsilon = 0.3$}, {$\sigma = 1\,\up{N\,m^{-1}}$}, $p_\up{atm} = 1\,\up{Pa}$, $\overline{\rho_1} =
1000\,\up{kg\,m^{-3}}$, $\overline{\rho_2} = 1\,\up{kg\,m^{-3}}$, {$\delta = 0.1$}, {$\omega =
\pi/3\,\up{m^{-1}}$}, {$\Pi_1 = 20\,\up{Pa}$}, {$\Pi_2 = 0\,\up{Pa}$}, {$\gamma_1 = 4$}, {$\gamma_2 = 1.4$}.

The computational domain is the square $\Omega = [-3\,\up{m},\ 3\,\up{m}]\times[-3\,\up{m},\ 3\,\up{m}]$ so
that at $t = 12\,\up{s}$ we expect the water column to have completed two full advection cycles. We
evolve the system from time $t = 0$ to time $t = 12\,\up{s}$ for all 
$\pnpm{N}{N}$, $N = 2,\ 3,\ \hdots,\ 9$ schemes with local space-time DG predictor step performed in
the primitive variable variant, using the HLL flux. The employed mathematical model is the
nonconservative hyperbolic Godunov--Powell formulation. 
The results confirm that the error norms of the conserved variables decrease at a rate that is in agreement
with the nominal order of accuracy of the scheme, and are summarised in
Tables~\ref{tab:convergencetable} and \ref{tab:convergencetablealpha}. In
Table~\ref{tab:convergencetablealpha}, we report the error norms and convergence rates for the 
liquid
volume fraction $\alpha_1$, for numerical schemes of order up to 10. In
Table~\ref{tab:convergencetable}, numerical details concerning the regression lines 
of
the $L_2$ norms of the error for all variables are given. Since the interface field $\vec{b}$ is 
evolved as a vector of independent state variables, as opposed standard schemes which differentiate 
the colour function and thus lose one order of accuracy for the discrete gradient, in our scheme 
the nominal high order convergence rate is achieved for the gradient field $\mathbf{b}$ as  
well. The regression lines for the mixture energy density $\rho\,E$ and for the liquid volume
fraction $\alpha_1$ are plotted in Figure~\ref{fig:droplet-convergence}, where also 
one-dimensional cuts through the numerical solution are shown along the $x$ axis, 
comparing the exact solution derived in Section~\ref{sec:exactsol} with 60 uniformly 
spaced samples from a computation employing the $\pnpm{9}{9}$ \aderdg scheme on a very 
coarse uniform Cartesian grid composed of only $6^2$ total elements.

\subsection{Interaction between a shock wave and a water column}

With this test case we want to show that the ADER-DG schemes with \emph{\aposteriori} subcell
Finite Volume limiter are capable of computing solutions of nonconservative hyperbolic systems 
not only in smooth regions, but can also robustly deal with shock waves while preserving the 
sharp profile that characterises these flow features. Specifically, we want to reproduce the 
results of the experiment of Igra and Takayama \cite{igra2002}, as it was done in \cite{Schmidmayer2017}. 

\begin{figure}[!b]
    \includegraphics[draft=false, scale=\figurescalefactor]{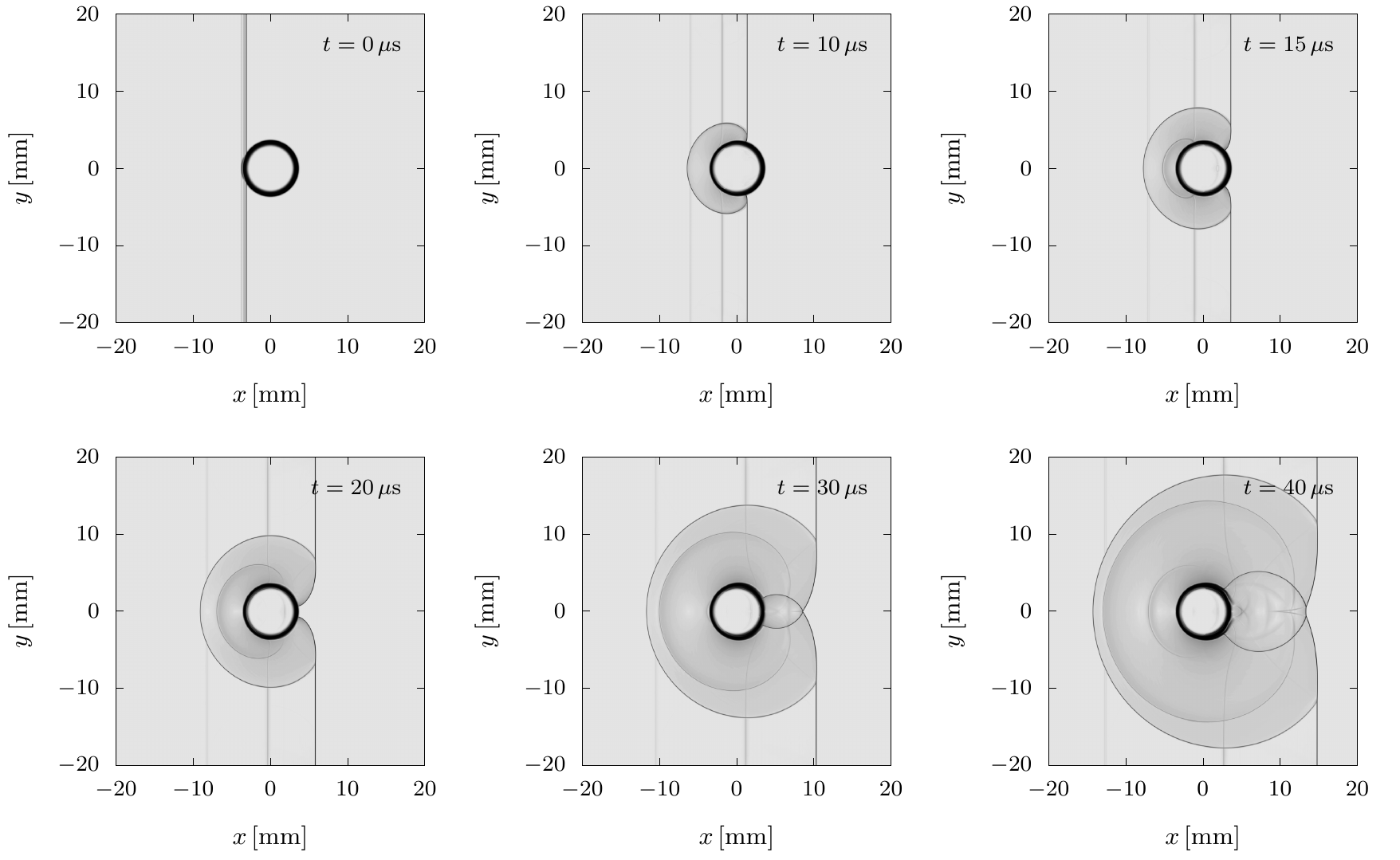}
    \caption{Numerical schlieren images of the early stages of the shock--water column interaction
    problem computed with a $\pnpm{3}{3}$ \aderdg scheme and a $\up{TVD}$ subcell limiter on a mesh of
    spacing $\Delta x = \Delta y = 0.0625\,\up{mm}$.}
    \label{fig:shock-droplet-early-schlieren.pdf}
\end{figure}

\begin{figure}[!b]
    \includegraphics[draft=false, scale=\figurescalefactor]{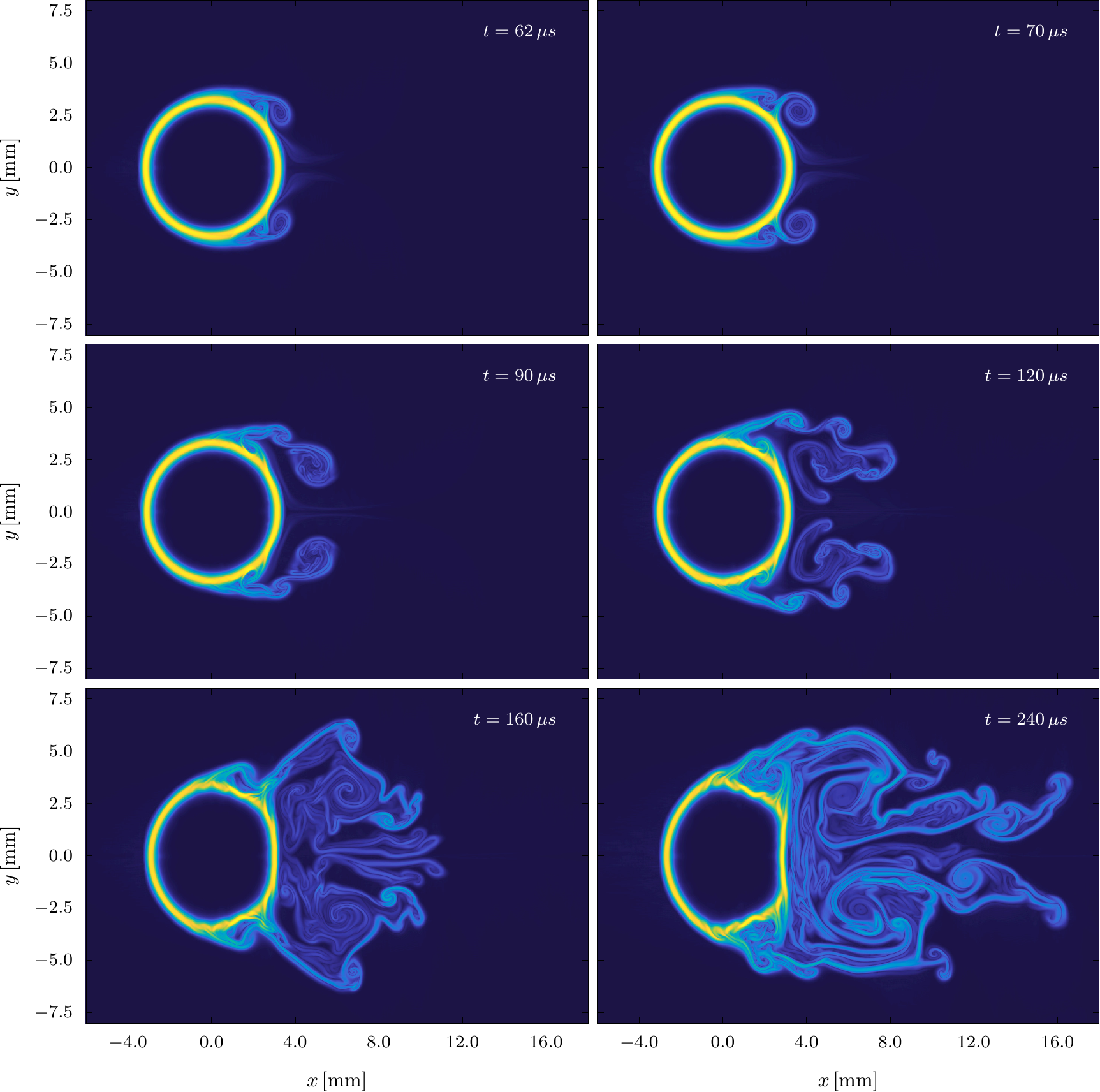}
    \caption{Visualisation of the interface transport by means of the filled contour plot of
    ${\sigma\,\norm{\vec{b}}}^{1/4}$ for the shock--water column interaction problem. The solution has
    been obtained with a $\pnpm{3}{3}$ \aderdg scheme and a $\up{TVD}$ subcell limiter on a uniform mesh
    of spacing $\Delta x = \Delta y = 0.0625\,\up{mm}$.}
    \label{fig:hst-shock-droplet-1002}
\end{figure}

The simulation setup is as follows: a cylindrical water column of radius $R = 3.2\,\up{mm}$ is
initialised at the origin of the computational domain following the exact solution given in Section
\ref{sec:exactsol}. The interface thickness parameter is $k_\epsilon = 1/20$ and the surface tension
coefficient is that of water, i.e. $\sigma = 0.072\,\up{N\,m^{-1}}$. Outside the water column, the 
pressure
is set to $p = p_\up{atm} = 10^5\,\up{Pa}$ and the liquid volume fraction is $\alpha_1 = 10^{-5}$,
while inside the droplet we have $\alpha_1 = 1 - 10^{-5}$. The density for water and air
are taken as $\rho_1 = 998.2\,\up{kg\,m^{-3}}$ and $\rho_2 = 1.18\,\up{kg\,m^{-3}}$, respectively.
The parameters for the equation of state are the usual ideal gas parameters for air $\Pi_2 =
0\,\up{Pa}$ and $\gamma_2 = 1.4$, while for the water we wanted to reproduce the correct speed of
sound in the pure liquid phase so we set $\Pi_1 = 4.7\!\times\!10^{8}\,\up{Pa}$ and $\gamma_1 =
4.7$. Since a perfectly pure phase is never present in our test, the speed of sound in water is
slightly smaller than the correct one, but still significantly larger than the speed of sound in air,
where the correct speed of sound is very well reproduced for $\alpha_1 = 10^{-5}$. A shock moving at
speed $S_\up{s} = 1.3\,a_\up{air}$ in the $x$ direction ($a_\up{air}$ being the speed of sound in
air, that is, we have a Mach 1.3 shock), with post-shock state computed following the jump relations
found in \cite{saurel2007}, is localised, at time $t = -1.0\,\mu \up{s}$, when the simulation is
started, at a distance $\delta = S_\up{s}\,\up{m\,s^{-1}}\cdot\,1.0\,\mu \up{s}$ (rounded to the
nearest element edge) from the \emph{nominal} edge of the droplet (see
Figure~\ref{fig:shock-droplet-early-schlieren.pdf} for a snapshot at time $t = 0\,\mu \up{s}$). A
result of this sharp initialization of the shock profile can be clearly observed in the numerical
schlieren images of Figure~\ref{fig:shock-droplet-early-schlieren.pdf}, in that two acoustic waves
due to the startup error can be seen travelling upstream in the post-shock region.


In order to produce the results discussed in this section, we ran, for convenience, two distinct 
simulations with different domain sizes, but with the same mesh spacing. 
One simulation deals with the early phase of the simulation, that is the impact between the 
shock and the water
column and the computational domain is the square $[-20\,\up{mm},\
20\,\up{mm}]\times[-20\,\up{mm},\ 20\,\up{mm}]$, while for the simulation on longer 
timescales we adopt a
rectangular domain $[-10\,\up{mm},\ 30\,\up{mm}]\times[-10\,\up{mm},\ 10\,\up{mm}]$. For
the solution we employ a fourth order ADER-DG $\pnpm{3}{3}$ scheme with primitive variable predictor
step, supplemented with a robust second order TVD limiter with primitive variable
reconstruction. The element size is the same for both simulations, since we use a grid of $640  
\times 640$ cells in the former case, and of $320 \times 640$ cells in the latter. The numerical fluxes
are computed with the HLL approximate Riemann solver.

In order to visualise the flow field, we plot the commonly used numerical schlieren pictures for the early 
stages of the simulation to highlight the shockwaves and aid comparison with the literature 
\cite{Schmidmayer2017, igra2002}, while, for the later stages of the simulation, we employ the key
variable of the model, that would be the interface field $\vec{b}$, to construct images that are very
rich in detail and show quite effectively the complex turbulent structures which develop in this
test problem, in a manner that is reminiscent of numerical schlieren pictures, since these are also 
nonlinearly scaled plots of the magnitude of a gradient. 

In Figure~\ref{fig:shock-droplet-early-schlieren.pdf}. one can see the first time instants of the
numerical experiment: discontinuities are very sharp, travel with the correct speed and in general
show good agreement with both the experimental data of \cite{igra2002} and the simulations of
\cite{Schmidmayer2017}. It is then notable that at time $t = 0\,\mu \up{s}$ some interaction can be
observed between the shock and the smoothing region of the water column, which extends symmetrically
towards the centre of the water column and towards the environment past the \emph{nominal} edge.

In
Figures~\ref{fig:hst-shock-droplet-1002}, at time $t = 62\,\up{ms}$ we can see the first vortical
structures developing around the water column and, at time $t = 120\,\mu \up{s}$,
Kelvin--Helmoltz-type \cite{kelvin1871, helmoltz1868} instabilities are clearly distinguishable,
while at time $t = 240\,\mu \up{s}$ one can also observe the presence of Richtmyer--Meshkov-type
instabilities \cite{richtmyer1961, meshkov1969}.

With this visualization method, the process of formation of filaments in the edge of the water
column, which then are drawn into the vortical flows in the wake of the obstacle, is quite apparent.

\subsection{Droplet transport in two and three space dimensions}

In this section, we conduct a systematic study of the stability and accuracy of the two new 
strongly hyperbolic systems of governing equations that have been proposed in this paper, which 
are both different from the original weakly hyperbolic model introduced in \sch. First, we set up a  
two-dimensional droplet in 
equilibrium, as prescribed by the exact solution given in Section~\ref{sec:exactsol}, in a uniform
velocity field with periodic boundary conditions, and track the time evolution of the
domain-averaged curl constraint violations. The problem is analogous as the one used for the
convergence study and is chosen because an exact solution for the problem is available, which allows
to assess the correctness of the results unequivocally. Differently from what has been done in the
convergence study, the sinusoidal density field given in Eq.~\eqref{eq:convergencesetup-rho}, is
replaced with two constant density values with ratio $\rho_1/\rho_2 = 1000$. In two space
dimensions, the same test is repeated for the original weakly hyperbolic model of \sch, for the new
hyperbolic formulation using the Godunov--Powell-type nonconservative products (denoted by GPNCP in the plots), 
and for another three 
runs with the new augmented hyperbolic GLM curl-cleaning system, with increasing values of the cleaning
speed $c_h$, namely choosing $c_h \in \{10,\ 20,\ 40\}\,\up{m\,s^{-1}}$. For each one of these five
choices, we let the computations run up to a final time $t_\up{end} = 20.0\,\up{ms}$, which
corresponds to 20 full advection cycles, first on a coarse mesh of $16^2$ cells, and then on a finer
grid counting $32^2$ elements, with the ADER-DG \pnpm{5}{5} scheme with ADER-WENO \pnpm{0}{2}
\aposteriori subcell limiter. The purpose of these runs is to verify
how the different formulations react to mesh refinement and how they compare for a given resolution.

Then we carry out another set of five runs, studying the advection of a three-dimensional droplet
with the ADER-DG \pnpm{3}{3} scheme with ADER-WENO \pnpm{0}{2} \aposteriori subcell limiter, 
on a coarse mesh of $16^3$ elements, to extend the previous
two-dimensional results to the full three-dimensional case.

\begin{figure}[!b]
    \includegraphics[draft=false, scale=\figurescalefactor]{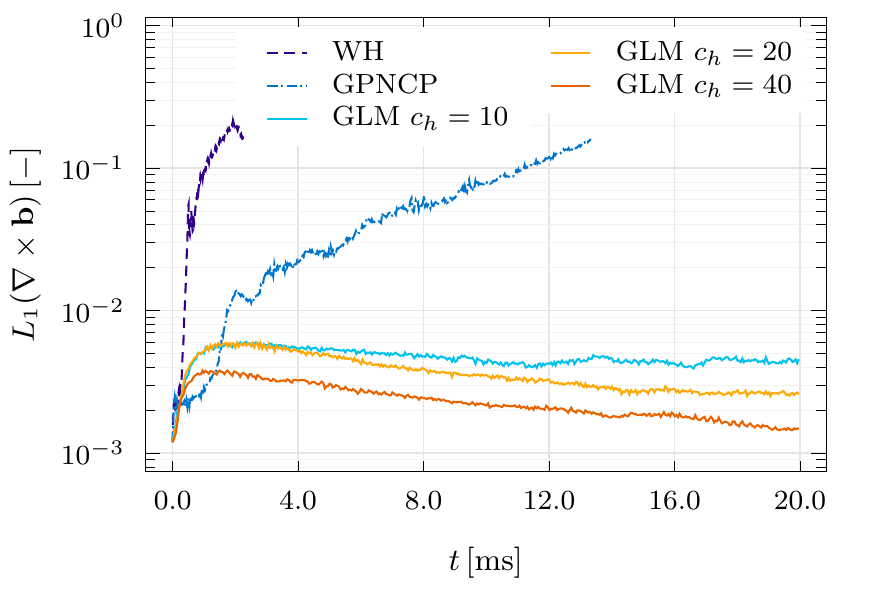}%
    \includegraphics[draft=false, scale=\figurescalefactor]{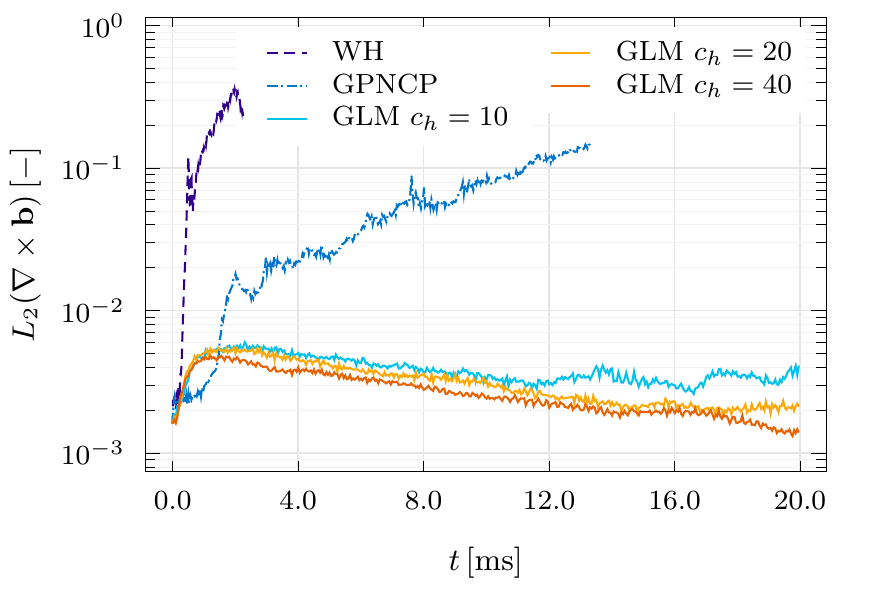}\\[4mm]
    \includegraphics[draft=false, scale=\figurescalefactor]{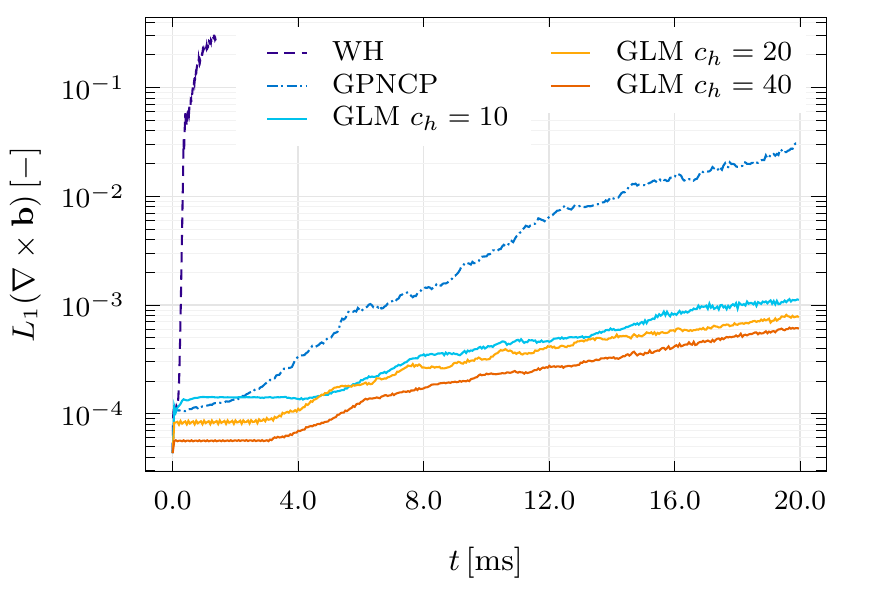}%
    \includegraphics[draft=false, scale=\figurescalefactor]{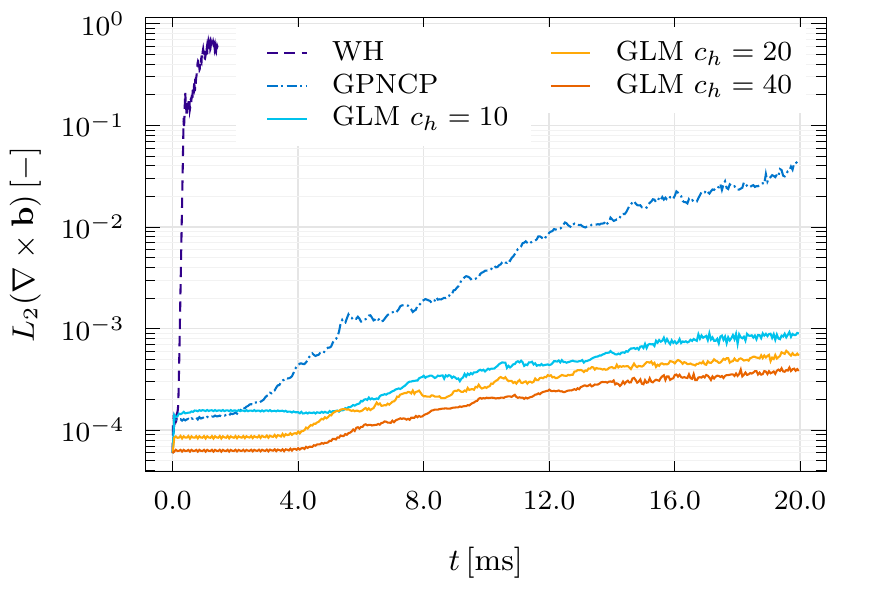}%
    \caption{Time evolution of the $L_1$ and $L_2$ norms of the curl constraint violations for the 
    two-dimensional droplet advection problem. In the top row, the results from a $\pnpm{5}{5}$ \aderdg
    scheme with \aderweno $\pnpm{0}{2}$ subcell limiter on a uniform coarse grid composed of $16^2$
    elements; in the bottom row, the results from the same method, but on a finer mesh composed of
    $32^2$ cells. In both cases, it is verified that curl errors are significantly reduced with the new GLM
    curl cleaning, with respect to those generated with the nonconservative Godunov--Powell-type formulation
    of the equations (GPNCP). In the GLM formulation, the constraint violations decrease with increasing 
    cleaning speed $c_h$, as expected. 
    Furthermore, on the coarser grid, the computation with the Godunov--Powell formulation fails after
    about 13 advection cycles ($13\,\up{ms}$). In no case stable results can be obtained from the
    original weakly hyperbolic model.}
    \label{fig:curl-errors-circle}
\end{figure}

\begin{figure}[!b]
    \includegraphics[draft=false, scale=\figurescalefactor]{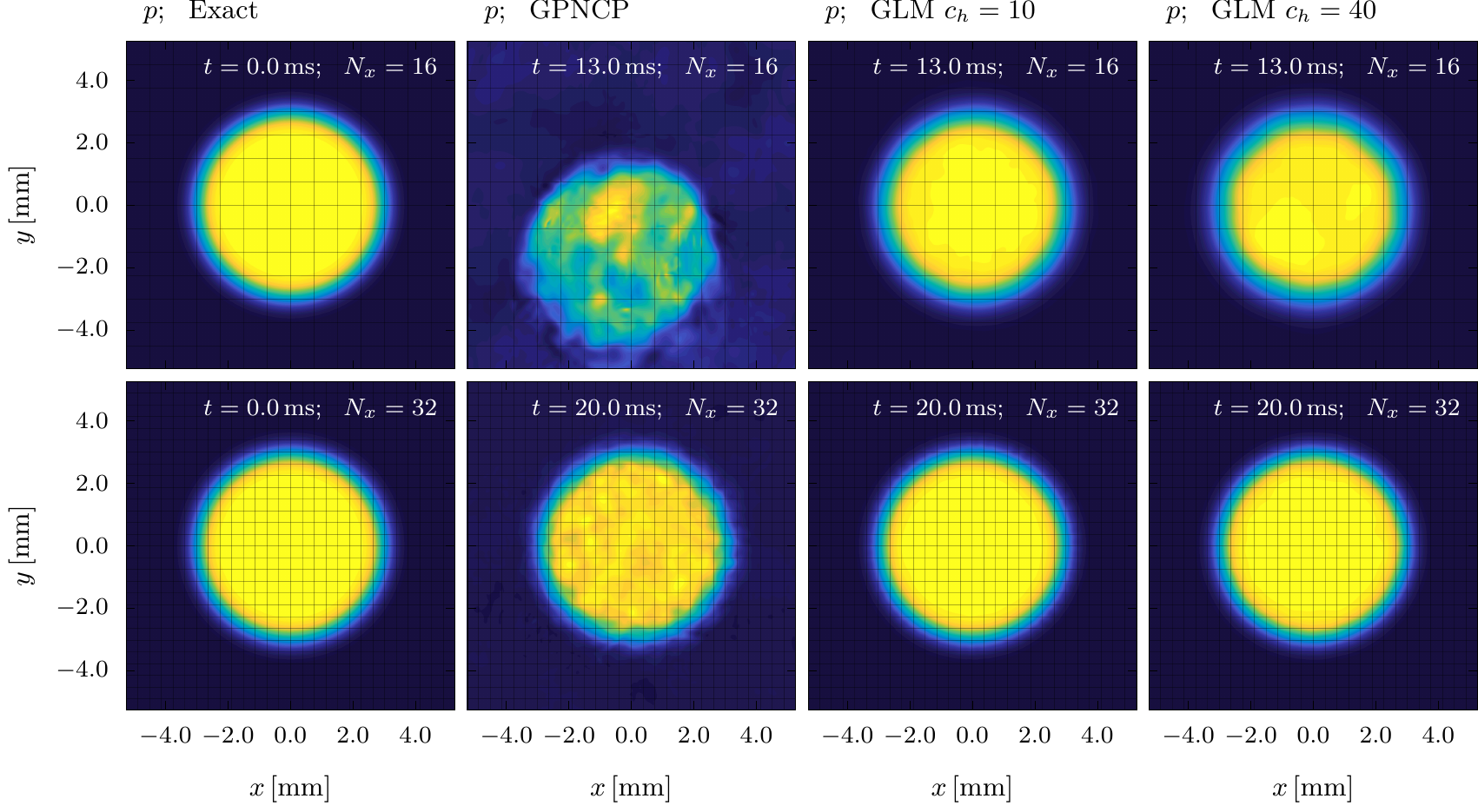}
    \caption{Numerical results for the two-dimensional droplet advection test problem. In the first 
    row,  we compare the results of the nonconservative Godunov--Powell model and
    of the augmented GLM curl cleaning system with two different values of cleaning speed $c_h$ 
    with the exact solution. The
    snapshots of the pressure field are taken at time $t=13.0\,\up{ms}$, which corresponds to thirteen
    advection cycles. The scheme used is ADER-DG \pnpm{5}{5} with ADER-WENO \pnpm{0}{2} subcell limiter
    and the mesh is composed of $16^2$ square control volumes. In the second row, the same comparison is
    carried out again, but on a finer mesh of $32^2$ elements at time $t = 20.0\,\up{ms}$, or after 20
    advection cycles. The results from the nonconservative model show a significant deviation
    from the exact solution of the problem, while the GLM curl cleaning approach yields very stable
    and accurate results: on the coarser mesh, some numerical diffusion is visible by comparing with the
    exact solution, but on the finer mesh numerical diffusion can be considered negligible.}
    \label{fig:twodimensional-advection-12}
\end{figure}

\begin{figure}[!b]
    \includegraphics[draft=false, scale=\figurescalefactor]{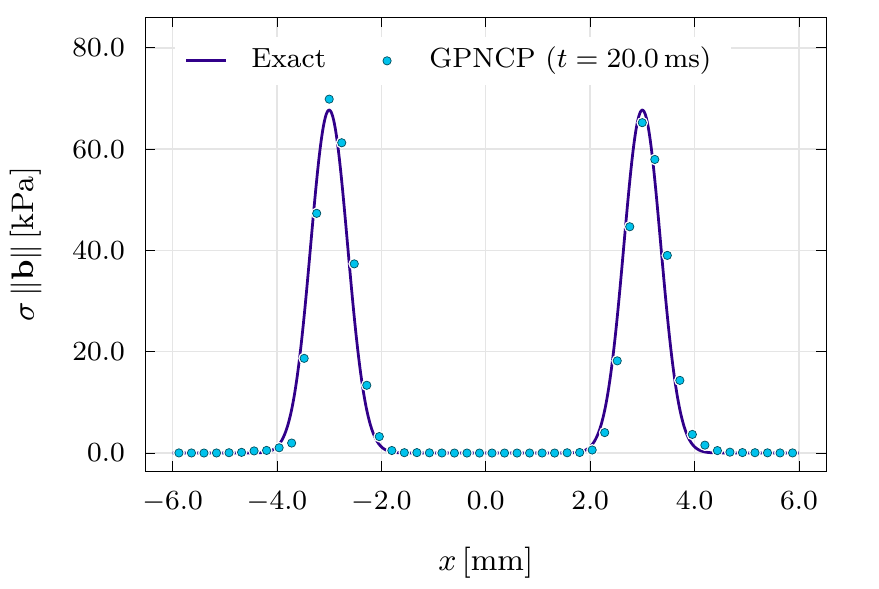}%
    \includegraphics[draft=false, scale=\figurescalefactor]{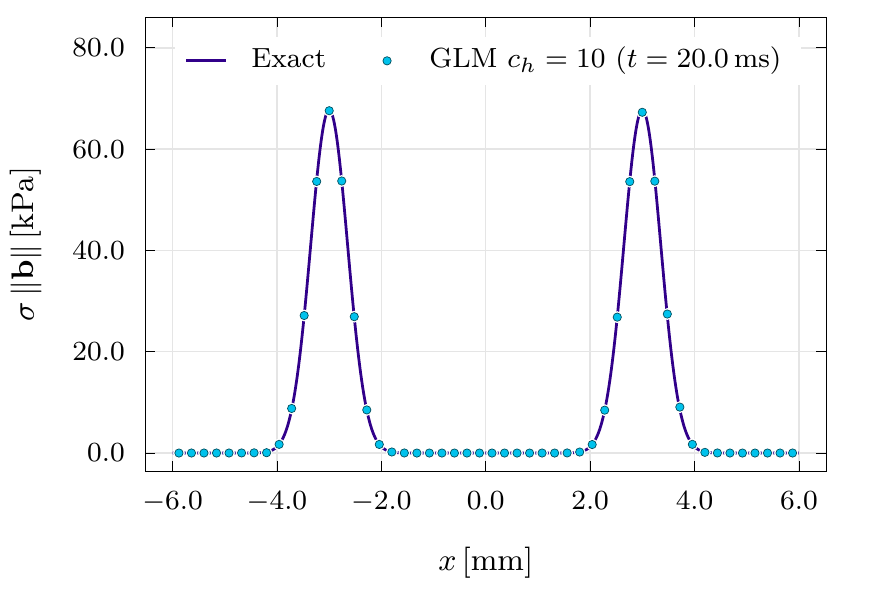}%
    \caption{One-dimensional cuts (50 uniform samples along the $x$ axis) of the interface
    energy $\sigma\,\norm{\vec{b}}$ for the two-dimensional droplet advection test problem. The 
    cuts are taken at time $t = 20.0\,\up{ms}$ (after 20 advection cycles)
    The scheme used is ADER-DG \pnpm{5}{5} with ADER-WENO \pnpm{0}{2} subcell limiter
    and the mesh is composed of $32^2$ square control volumes.}
    \label{fig:twodimensional-advection-3}
\end{figure}

\begin{figure}[!b]
    \includegraphics[draft=false, scale=\figurescalefactor]{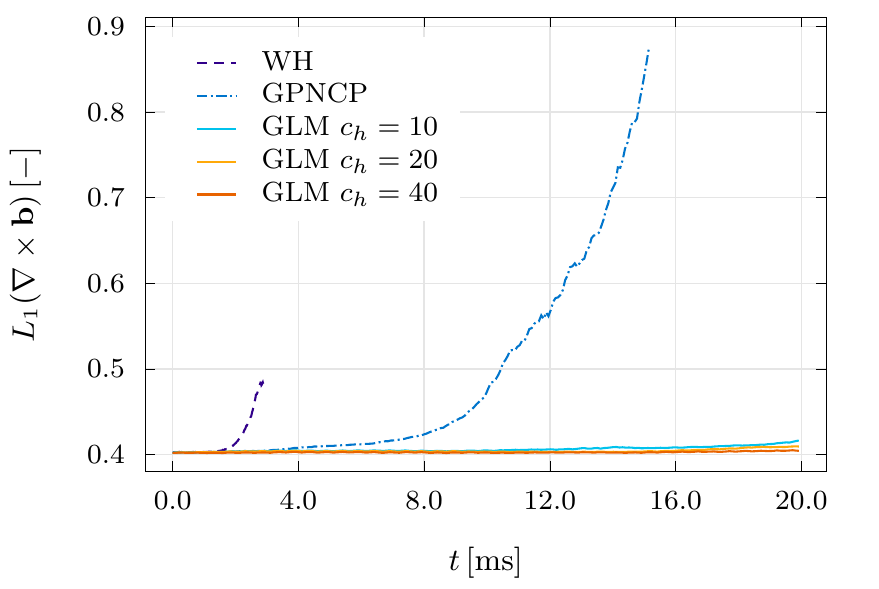}%
    \includegraphics[draft=false, scale=\figurescalefactor]{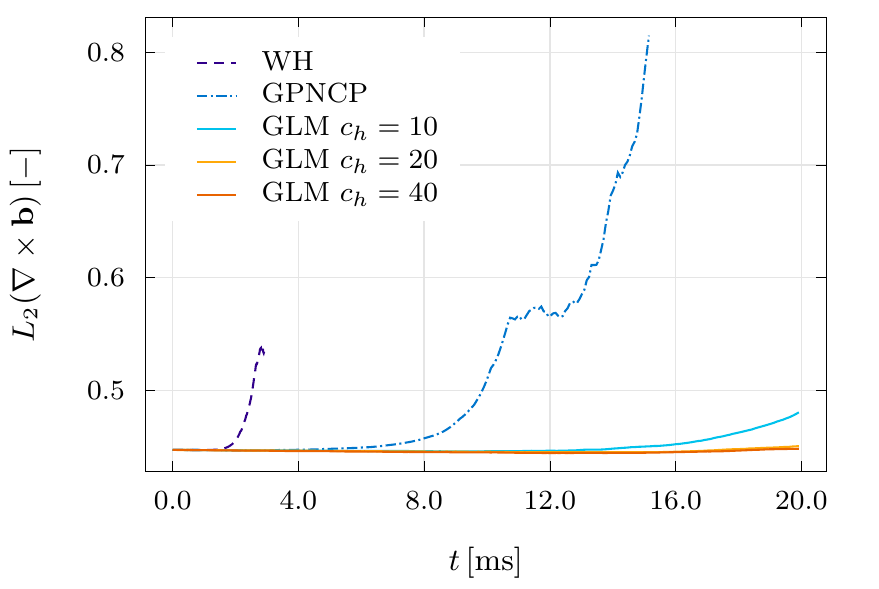}
    \caption{Time evolution of the $L_1$ and $L_2$ norms of the curl constraint violation for the
    three-dimensional droplet advection problem. %
    The results obtained for the two-dimensional
    experiments are confirmed. It is quite apparent that, employing the GLM curl cleaning, the
    constraint violation grows much slower than with the Godunov--Powell-type formulation of
    the system (GPNCP). In particular, in this latter case the computation fails after about 15 
    advection cycles ($15.0\,\up{ms}$), while the augmented GLM curl cleaning system shows much 
    better stability properties. 
    The simulations have been carried out using a $\pnpm{3}{3}$ \aderdg scheme with \aderweno $\pnpm{0}{2}$
    subcell limiter on a uniform grid composed of $16^3$ elements.}
    \label{fig:curl-errors-sphere}
\end{figure}

\begin{figure}[!b]
    \includegraphics[draft=false, 
    scale=\figurescalefactorb]{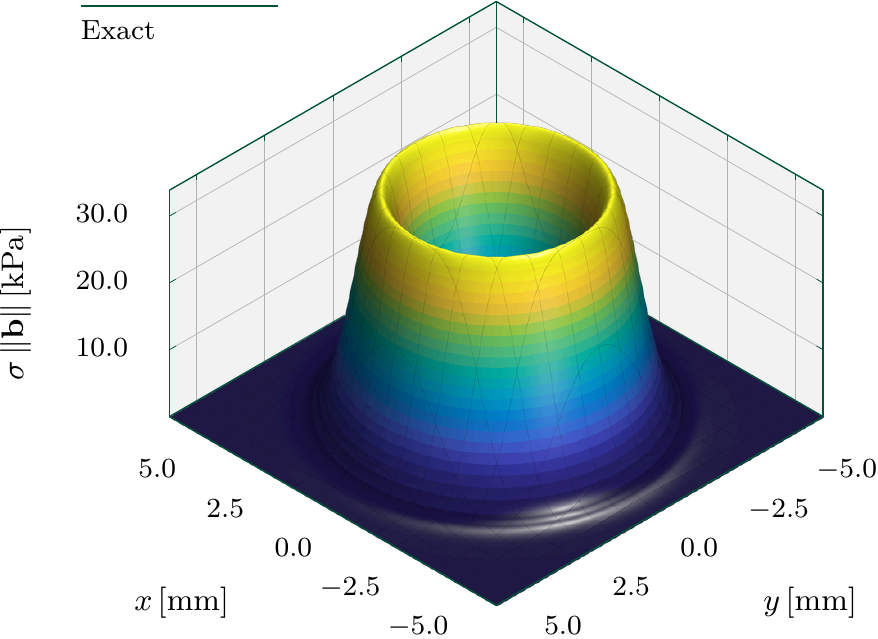}%
    \includegraphics[draft=false, 
    scale=\figurescalefactorb]{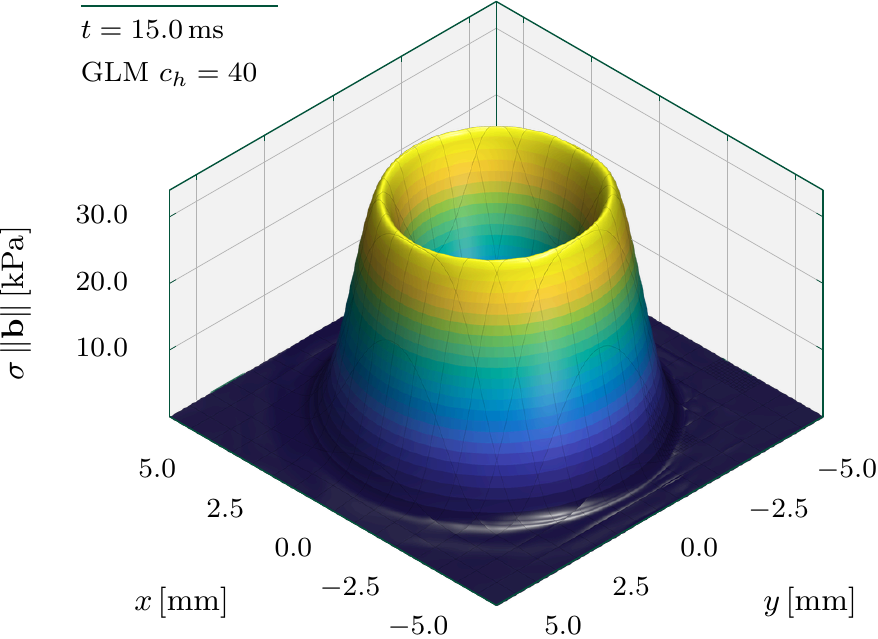}\\[8mm]
    \includegraphics[draft=false, 
    scale=\figurescalefactorb]{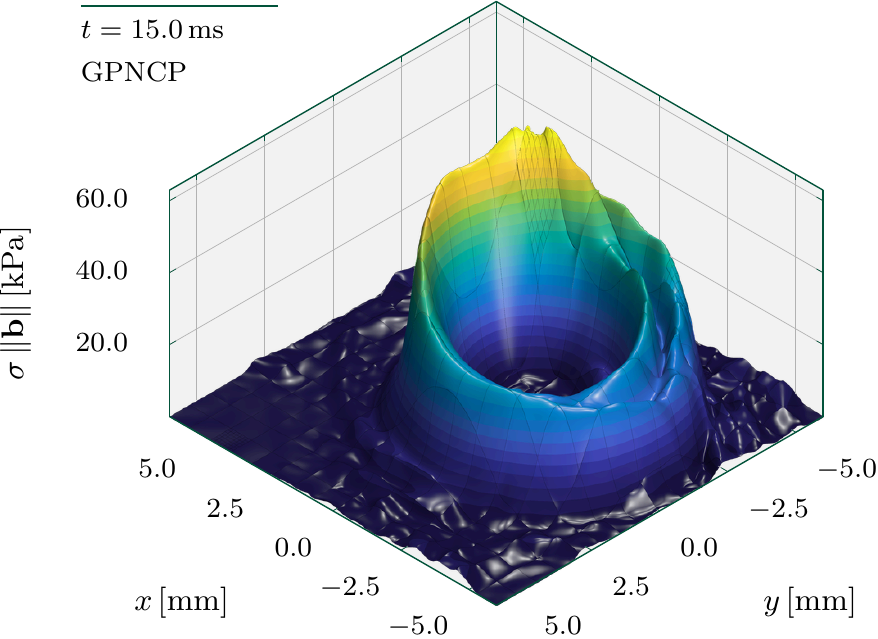}%
    \includegraphics[draft=false, 
    scale=\figurescalefactorb]{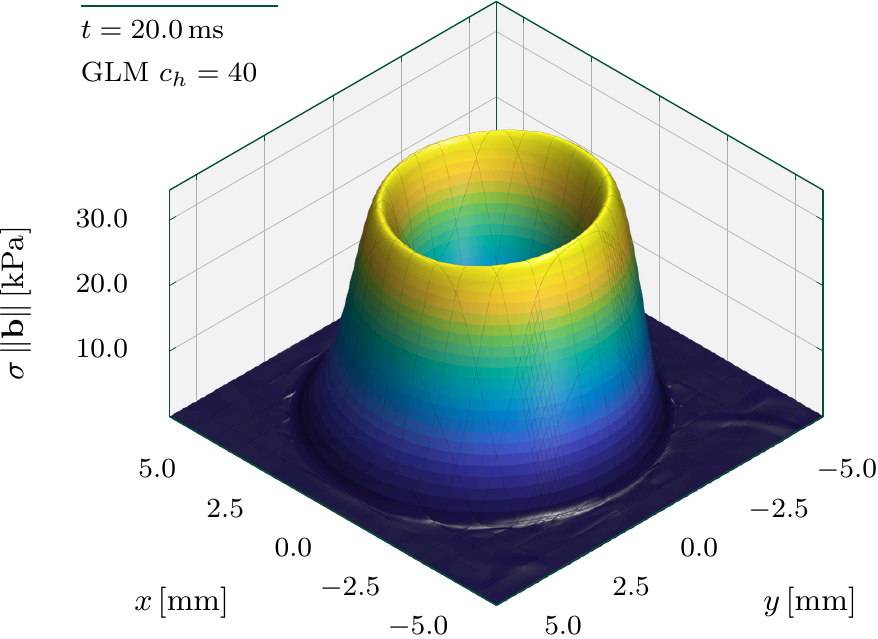}%
    \caption{Two-dimensional slices, at $z=0$, of the solution for the surface energy
    $\sigma\,\norm{\vec{b}}$, for the three-dimensional droplet advection test. The top left panel shows
    the exact solution of the problem. The top right plot is the numerical solution at time $t =
    15.0\,\up{ms}$, that is, after fifteen advection cycles, obtained using the GLM curl cleaning
    formulation of the model with cleaning speed $c_h = 40$, showing good agreement with the exact
    solution. The bottom left plot represents the solution at $t = 15.0\,\up{ms}$ obtained with the
    Godunov--Powell nonconservative formulation of the model (GPNCP); strong artefacts are visible. The
    bottom right plot shows the GLM curl cleaning solution after five additional advection cycles, with
    comparatively minor deformation of the interface. The results are obtained with a fourth order
    ADER-DG \pnpm{3}{3} and ADER-WENO \pnpm{0}{2} subcell limiter, on a very coarse mesh of $16^3$
    elements.}
    \label{fig:threedimensional-slices-2d}
\end{figure}

\begin{figure}[!b]
    \includegraphics[draft=false, scale=\figurescalefactor]{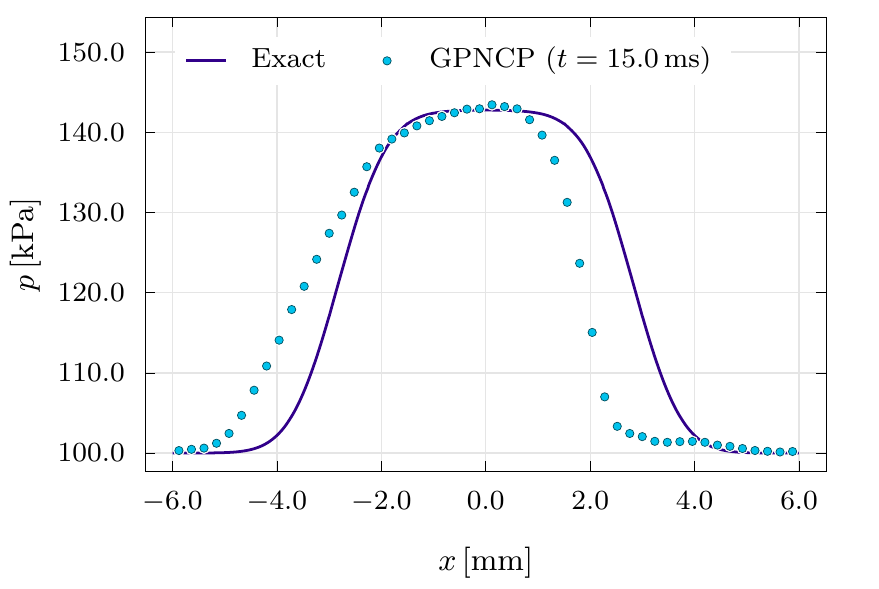}%
    \includegraphics[draft=false, scale=\figurescalefactor]{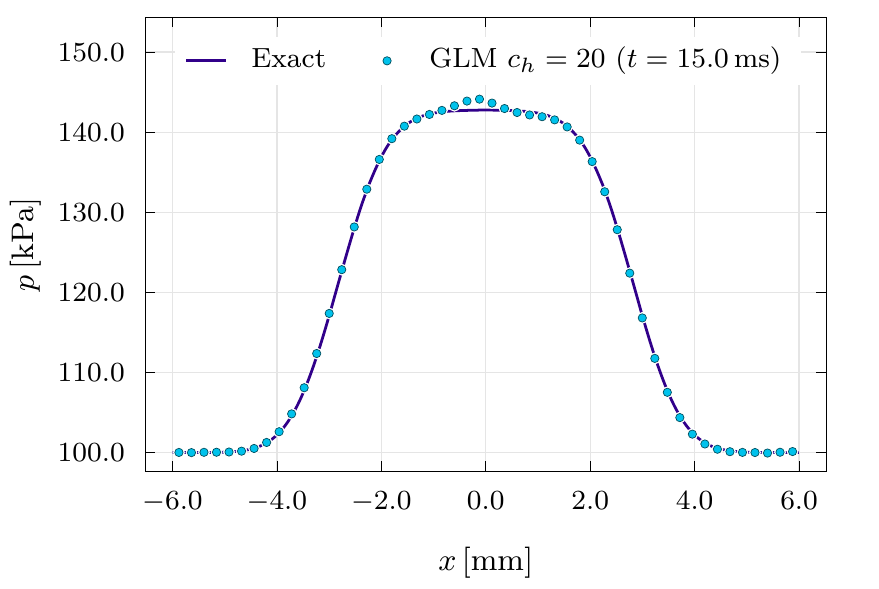}%
    \caption{One-dimensional cuts (50 uniform samples along the $x$ axis), at $z=y=0$, of the solution
    for the pressure field $p$, for the three-dimensional droplet advection problem. On the left, the
    solution obtained with the Godunov--Powell nonconservative formulation of the model (GPNCP). On the
    right, the solution from the GLM curl cleaning formulation, with cleaning speed $c_h = 20$. One can
    note that, at time $t = 15.0\,\up{ms}$, after fifteen full advection cycles, the nonconservative
    formulation significantly deviates from the exact solution derived in
    Section~\ref{sec:exactsol}, while the error is very well contained by the GLM curl cleaning. The
    results are obtained with a fourth order ADER-DG \pnpm{3}{3} and ADER-WENO \pnpm{0}{2} subcell
    limiter, on a very coarse mesh of $16^3$ elements.}
    \label{fig:threedimensional-slices-1d}
\end{figure}

The droplet has radius $R = 3\,\up{mm}$ and is centred at the origin of a square domain $\Omega =
[-6\,\up{mm},\ 6\,\up{mm}]\times[-6\,\up{mm},\ 6\,\up{mm}]$, the liquid and gas density are
respectively set to $\rho_1 = 1000\,\up{kg\,m^{-1}}$ and $\rho_2 = 1\,\up{kg\,m^{-1}}$ throughout
the domain. The volume fraction follows Eq.~\eqref{eq:convergencesetup-alpha}, with
{$\alpha_\up{min} = 0.01$} and {$\alpha_\up{max} = 0.99$}, and the interface field is given by
\eqref{eq:exactb}, with the dimensionless interface thickness parameter being $k_\epsilon = 1/6$
for the two-dimensional tests and $k_\epsilon = 1/3$ for the three-dimensional problem, additionally
setting $b_3 = 0$ for the two-dimensional runs. The pressure is initialised following the exact
solution \eqref{eq:compactpressureprofile}, with atmospheric pressure $p_\up{atm} = 100\,\up{kPa}$,
and the uniform velocity field components are $u_0 = 12\,\up{m\,s^{-1}}$, $v_0 =
12\,\up{m\,s^{-1}}$, and $w_0 = 12\,\up{m\,s^{-1}}$ in three space dimensions or $w_0 =
0\,\up{m\,s^{-1}}$ in two dimensions. The parameters for the equation of state are $\Pi_1 =
1\,\up{MPa}$, $\Pi_2 = 0$, $\gamma_1 = 4$, $\gamma_2 = 1.4$, and the surface tension coefficient is
set to $\sigma = 60\,\up{N\,m^{-1}}$.

The results are depicted in Figures~\ref{fig:curl-errors-circle}, 
\ref{fig:twodimensional-advection-12}, and \ref{fig:twodimensional-advection-3}.
In Figure~\ref{fig:curl-errors-circle}, we plot the time evolution of the normalised $L_1$ and 
$L_2$ norms of the 
curl constraint violations, defined as
\begin{equation}
\label{eq:errornorms}
    \left. L_1(\nabla\times\vec{b})\right. = \dfrac{\left.\displaystyle\int_\Omega 
        \sqrt{\,\transpose{(\nabla\times\vec{b})}\,\nabla\times\vec{b}\,}\right.}
        {\left.\displaystyle\int_\Omega \sqrt{\,\trace{\left[\transpose{(\nabla\vec{b})}\,
        \nabla\vec{b}\right]\,}}\right|_{t=0}},\qquad
    \left. L_2(\nabla\times\vec{b})\right. = \dfrac{\left.\displaystyle\int_\Omega 
        {\abs{\,\transpose{(\nabla\times\vec{b})}\,\nabla\times\vec{b}\,}}\right.}
        {\left.\displaystyle\int_\Omega {\abs{\,\trace{\left[\transpose{(\nabla\vec{b})}\,
        \nabla\vec{b}\right]\,}}\,}\right|_{t=0}}.
\end{equation}
We observe that in all cases the same trend is apparent: the curl error given by the weakly
hyperbolic model quickly grows until the computation terminates with unphysical values at 
rather early times, while the new strongly hyperbolic 
variants of the governing equations are stable, at least with increasing mesh refinement. It seems that not 
much can be done to improve the stability of the weakly hyperbolic model, which in the run with 
finer mesh blows up even earlier than with the coarse grid, most likely due to the smaller numerical 
dissipation of the scheme. In the long term, it is always true that the curl errors are lower with GLM curl 
cleaning 
than they are with the nonconservative Godunov--Powell-type model. One can also see that the higher 
the cleaning speed $c_h$ is, the smaller the constraint violations are. 
Moreover, on the fine mesh, the nonconservative Godunov--Powell system, while still generating much
larger errors than the augmented GLM model (clearly visible also in the pressure field shown in 
Figure~\ref{fig:twodimensional-advection-12}), could be solved for the full 20 advection cycles, as
opposed to only 13 on the coarse mesh.  

Concerning the effects of numerical dissipation, we can see that the curl errors for the GLM curl
cleaning simulations on the coarse grid decrease in time with the aid of numerical diffusion, which
reduces the overall steepness of the interface field. This effect can be easily quantified by
inspecting Figure~\ref{fig:twodimensional-advection-12} where it is apparent that with the coarse
mesh the pressure field after thirteen advection cycles is more diffused than in the initial
condition, while this effect is minimised by mesh refinement, as one can clearly see in
Figure~\ref{fig:twodimensional-advection-3}, where the profile of the interface field on the GLM
simulations is still in perfect agreement with the exact solution after 20 full advection
cycles. Regarding this simulation with the finer grid, the curl error timeseries no longer shows the
effects of numerical dissipation and in the first stages of the computation (up to about three to
four advection cycles) one can see that the curl errors are maintained at a very precise constant
value, suggesting that a sort of balance is established between the sources of the curl errors in 
the numerical scheme and their transport via the Maxwell-type curl cleaning waves of the augmented GLM system. 
Also, one can note that, for the run with cleaning speed $c_h =
40\,\up{m\,s^{-1}}$, in this early phase, the curl error is kept very close to its non-zero initial
value, which is given by the necessity of projecting the pressure profile on the
piecewise-polynomial \disgal data representation, even if evaluated at machine precision from an
exact formula.

In the three-dimensional tests, the effects of numerical diffusion are not seen because the
interface profile was chosen to be smoother than the one used for the two-dimensional simulations 
from the beginning. 
Otherwise, the same observations given for the two-dimensional case are valid, namely one can
construct a hierarchy of the simulations based on the entity of the curl-constraint violations, that
sees the weakly hyperbolic model break down very early, the Godunov--Powell-type symmetrisable model
being more stable, but more sensitive in the long term than the GLM cleaning simulations, which in
turn have lower errors for higher cleaning speeds. The timeseries of the constraint violations are
plotted in Figure~\ref{fig:curl-errors-sphere}, where the error is kept essentially equal to the
initial value with the GLM curl cleaning, while it grows rather quickly for the Godunov--Powell
formulation, for which the computation stops after completing 15 advection cycles.

In Figure~\ref{fig:threedimensional-slices-2d}, we show a set of two-dimensional slices of the
solution for the interface energy and we observe that, as for the analogous two-dimensional test,
the hyperbolic Godunov--Powell model shows severe degradation of the interface field after fifteen
advection cycles and the droplet is even shifted out of centre, as was the two-dimensional droplet in
the second panel of Figure~\ref{fig:twodimensional-advection-12}. At the same time instant, the GLM
curl cleaning formulation seems to adequately match the exact solution, despite using a rather
coarse mesh, and shows no spurious shift of the centre of mass of the droplet, as seen also in the
one-dimensional cuts of Figure~\ref{fig:threedimensional-slices-1d}.

\subsection{Oscillation of an elliptical water column} We continue our systematic comparison of the
different formulations of the hyperbolic surface tension model under investigation with a test
involving the oscillation of an elliptical water column, which, due to the elongated initial shape,
is not in mechanical equilibrium and tends to deform towards restoring a circular shape. The
phenomenon is of periodic nature since when the droplet has indeed reached a circular shape, it also
stores an amount of kinetic energy such that it starts to elongate again perpendicularly with
respect to the previous major axis, up to a maximum deformation, then deforming back to a circular
shape and finally to the initial configuration.

\subsubsection{Problem setup} 
\begin{figure}[!b]
    \includegraphics[draft=false, 
    scale=\figurescalefactorb]{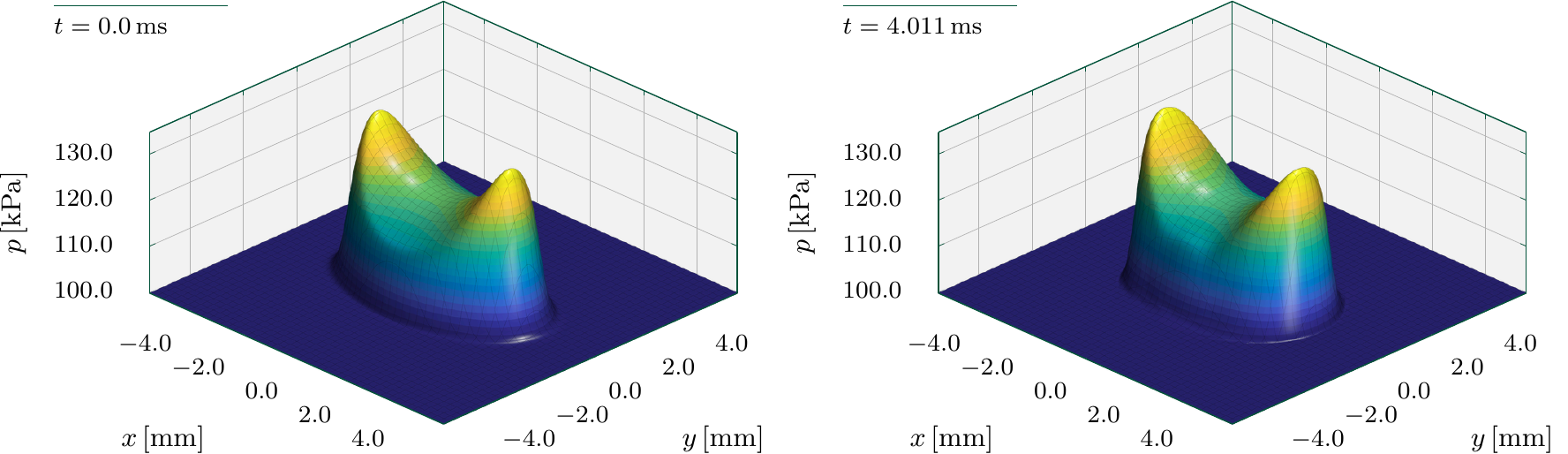}\\[8mm]
    \includegraphics[draft=false, scale=\figurescalefactorb]{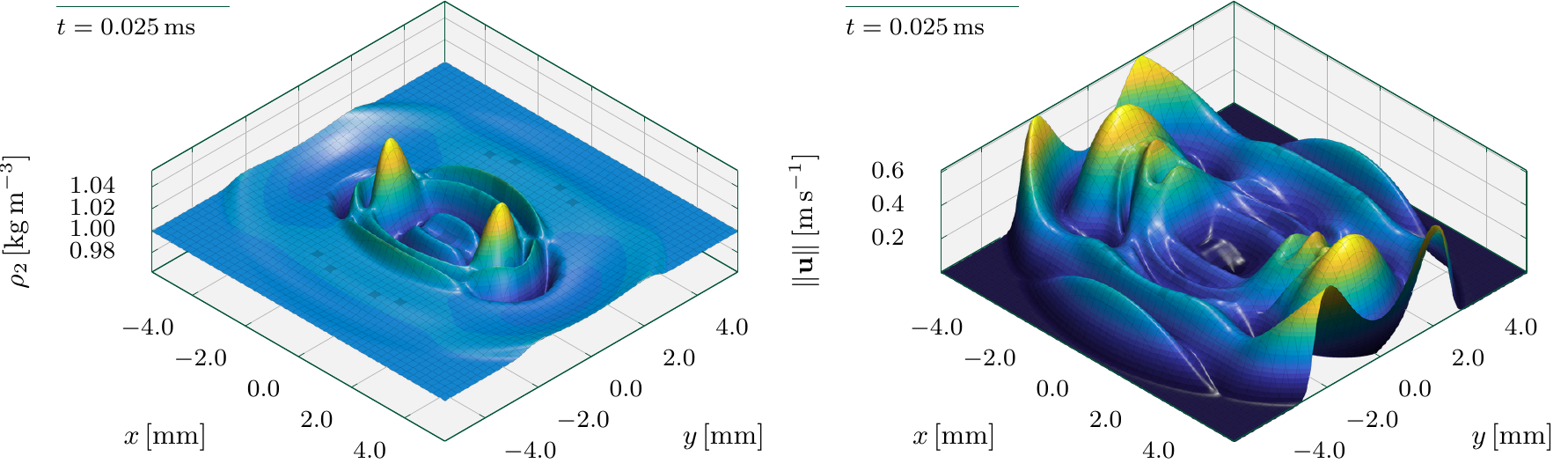}
    \caption{Initial stages of the elliptical droplet oscillation problem, In the top row, 
    oscillation periods ($t = 4.011\,\up{ms}$). In the bottom row, density of the gas phase (on the
    left) and flow speed distribution (on the right) at time $t=0.025\,\up{ms}$. The Godunov--Powell
    nonconservative formulation of the model was solved with a $\pnpm{7}{7}$ \aderdg scheme and
    \aderweno $\pnpm{0}{2}$ subcell limiter on a uniform grid of $50^2$ elements.}
    \label{fig:droplet_oscillation_presentation}
\end{figure}

For the description of the geometry of a smoothed elliptical
water column having a nominal interface defined by the parametric equation {${\vec{r}_\up{b} =
(x_\up{b},\ y_\up{b}) = (R_x\,\cos \psi,\ R_y\,\sin \psi)}$},

we introduce the following coordinates: for each point $(x,\ y)$ in the Cartesian plane the local
eccentric anomaly $\psi$ is defined implicitly by the formulas
\begin{equation}
        \cos^2\psi(x,\ y) = \frac{R_y^2\,x^2}{R_y^2\,x^2 + R_x^2\,y^2}, \qquad
        \sin^2\psi(x,\ y) = \frac{R_x^2\,y^2}{R_y^2\,x^2 + R_x^2\,y^2};
\end{equation}
for each point $(x,\ y)$, we can then define the nominal radius of the ellipse in the direction of
the local eccentric anomaly
\begin{equation}
    R_\psi(\psi(x,\ y)) = \sqrt{R_x^2\,\cos^2\psi + R_y^2\,\sin^2\psi}, 
\end{equation}
which would be the length of the segment running from the centre of the water column (located at the
origin of the reference system) to the intersection between the ellipse boundary and the line
connecting said generic point with the origin. Then we denote as usual with $r(x,\ y) = \sqrt{{x}^2
+ {y}^2}$ the distance of a generic point from the centre of the water column. 
Then, given its dimensionless form $\rs$ with respect to $R_\psi$, the colour function $c$, its
gradient $\vec{b} = [b_1,\ b_2,\ 0]$ and the liquid volume fraction $\alpha_1$ are given as
\begin{align}
    & c(x,\ y) = \frac{1}{2}\,\erfc{\left(\frac{r - R_\psi}{\epsilon}\right)},\\
    & \alpha_1(x,\ y) = \alpha_\up{min} + (\alpha_\up{max} - \alpha_\up{min})\,c(x,\ y),\\
    & b_1(x,\ y) = -\frac{x}{\sqrt{\pi}\,\epsilon\,r}\,%
        \exp{\left[-{\left(\frac{r - R_\psi}{\epsilon}\right)}^2\right]}\,%
        \left[1 - \left(1 - \frac{R_y^2}{R_x^2}\right)\frac{R_\psi}{r}\,\sin^2\psi\right],\\
    & b_2(x,\ y) = -\frac{y}{\sqrt{\pi}\,\epsilon\,r}\,%
        \exp{\left[-{\left(\frac{r - R_\psi}{\epsilon}\right)}^2\right]}\,%
        \left[1 - \left(1 - \frac{R_x^2}{R_y^2}\right)\frac{R_\psi}{r}\,\cos^2\psi\right],
\end{align}
while the pressure field is initialised as a local application of the solution for a cylindrical 
water column in the form
\begin{equation}
\begin{aligned}
    p\left(\rs\right) = p_\up{atm} + (d - 1)\,\frac{\sigma}{R}\,\int_{\rs}^{\infty}
    \frac{1}{\sqrt\pi\,\keps\,\rs^\prime}\,\exp\left[-{\left(\frac{\rs^\prime -
    1}{\keps}\right)}^2\right]\,\de{\rs^\prime},&\\
    \quad \text{with}\quad  r_\ast(x,\ y) = \frac{r(x,\ y)}{R_\psi(\psi(x,\ y))}, 
    \quad k_\epsilon = \dfrac{\epsilon}{R_\kappa},& 
\end{aligned}
\end{equation}
where by local we mean that an average curvature radius $R_\kappa$, defined at each point $(x,\ y)$
inside, on the nominal boundary, or outside of the droplet, is computed by averaging the curvature
along the nominal boundary of the ellipse, with a weight function inversely proportional to the
square of the distance from the interface ${\left[R_\psi(\psi) - r(x,\ y)\right]}^2$, so that we
have
\begin{equation}
    R_\kappa(x,\ y) = 
    {\left\{
    \dfrac{
        \displaystyle\int_{0}^{2\,\pi} \dfrac{
            1%
        }{
            {\left[R_\psi(\psi) - r(x,\ y)\right]}^2
        }\,\dfrac{\rx\,\ry}{\rx^2\,\sinq + \ry^2\,\cosq}\,\de{\psi}%
    }{
        \displaystyle\int_{0}^{2\,\pi} \dfrac{
            1%
        }{
            {\left[R_\psi(\psi) - r(x,\ y)\right]}^2
        }\,\sqrt{\rx^2\,\sinq + \ry^2\,\cosq}\,\de{\psi}%
    }
    \right\}}^{-1}.
\end{equation}
This averaging procedure yields a local curvature radius such that the initial pressure
configuration is similar to the one occurring at oscillation extrema, \textit{i.e.} at the end of
every half-period, when the kinetic energy of the droplet is zero, as it is set initially. Even if
based only on geometrical considerations, this initial condition is sufficient for individuating
very clearly only the main oscillation mode of the droplet, allowing to obtain a clean estimate of
the oscillation period. A comparison between the geometrically-derived initial pressure field and
the configuration after three oscillation periods is shown in
Figure~\ref{fig:droplet_oscillation_presentation}, together with the complex flow features that are
generated in the earliest instants of the simulation.

The density fields are set to the uniform values $\rho_1^0$ and $\rho_2^0$ throughout the
computational domain, as is the velocity field for which we set $\vec{u} = \transpose{[0,\ 0,\ 0]}$. The
numerical values employed for this test problem are: 
{$\rho_1^0 = 1000\,\up{kg\,m^{-3}}$}, {$\rho_2^0 = 1\,\up{kg\,m^{-3}}$}, 
{$p_\up{atm} = 100\,\up{kPa}$}, {$R_x = 3\,\up{mm}$}, {$R_y = 2\,\up{mm}$}, {$\alpha_\up{min} = 0.01$},
{$\alpha_\up{max} = 0.99$}, {$\sigma = 60\,\up{N\,m^{-1}}$}. The parameters for the
stiffened gas equation of state are: {$\Pi_1 = 1\,\up{MPa}$}, {$\Pi_2 = 0$}, {$\gamma_1 = 4$}, {$\gamma_2 =
1.4$}. The domain is the square $[-6\,\up{mm},\ 6\,\up{mm}]\times[-6\,\up{mm},\ 6\,\up{mm}]$ and additionally, the
initial condition is rotated counter-clockwise by 30 degrees, in order to avoid mesh alignment. 
In a first batch of tests, we set
{$\epsilon = 0.5\,\up{mm}$} and discretise the computational domain with $64^2$ square cells,
then solving with an ADER-DG $\pnpm{5}{5}$, ADER-WENO $\pnpm{0}{2}$ \finvol limiter and HLL flux. 
These simulations are intended to test the capability of the proposed
models in a dynamical setting where the interface deforms significantly under the effect of strong
surface tension, and verify that the GLM curl cleaning approach can deal with the violations of
involution constraints that such deformations generate.

In a second run, in order to study the sensitivity of results and in particular of the oscillation
period to the thickness of the diffuse interface region, we set {$\epsilon = 0.25\,\up{mm}$}, impose
no initial rotation of the droplet, thus aligning the two axis of the ellipse with the reference
frame, and compute the solution of the problem on a uniform grid of $50^2$ elements, with an eight
order ADER-DG $\pnpm{7}{7}$ scheme, ADER-WENO $\pnpm{0}{2}$ \finvol limiter, and Rusanov flux.

\begin{figure}[!b]
    \includegraphics[draft=false, scale=\figurescalefactor]{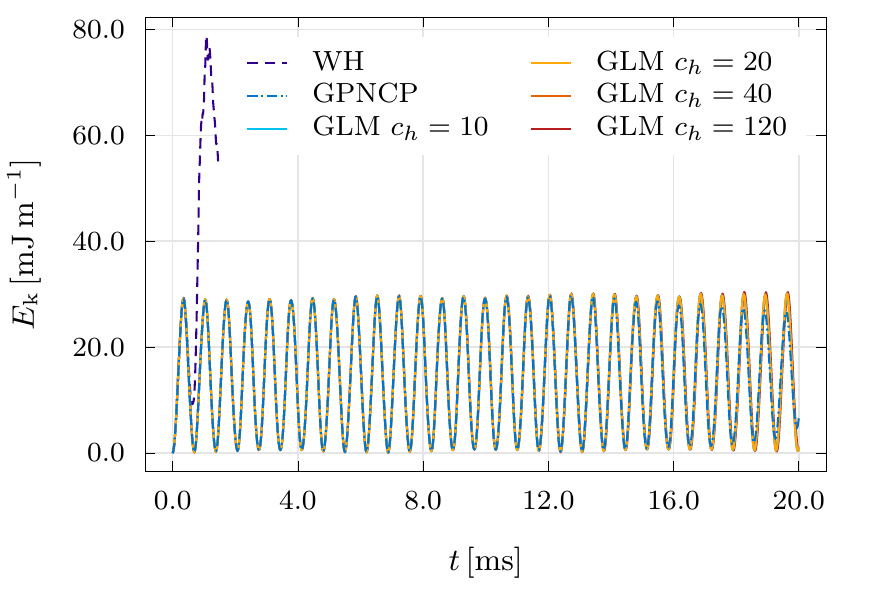}%
    \includegraphics[draft=false, scale=\figurescalefactor]{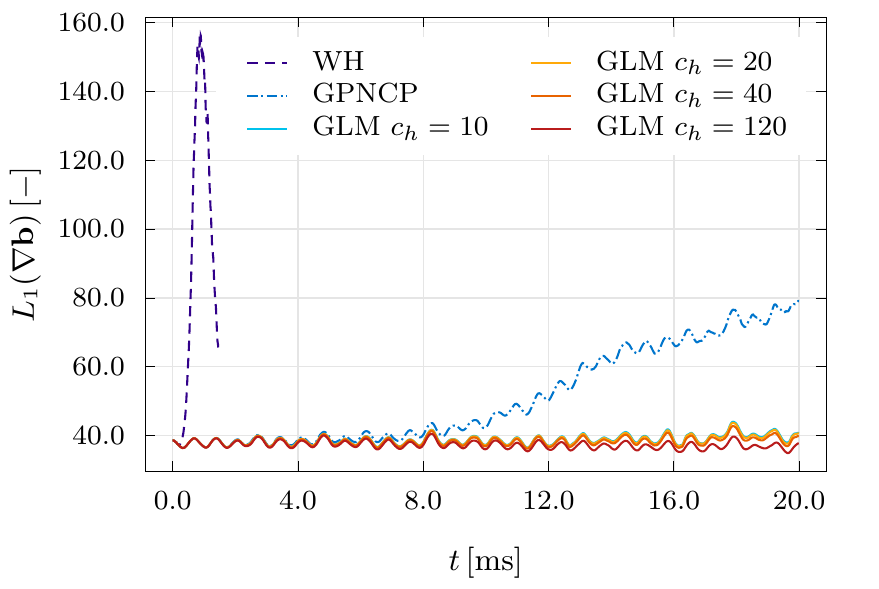}\\[4mm]
    \includegraphics[draft=false, scale=\figurescalefactor]{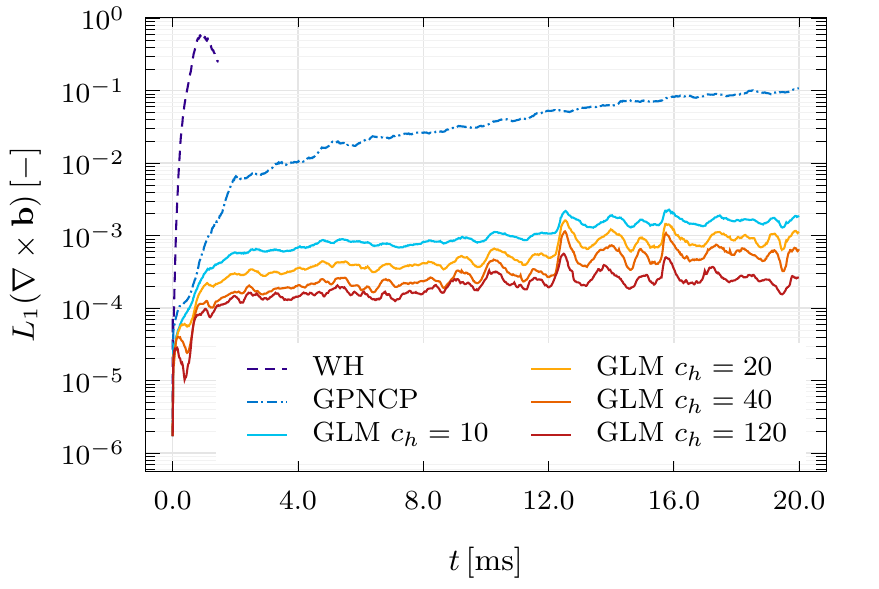}%
    \includegraphics[draft=false, scale=\figurescalefactor]{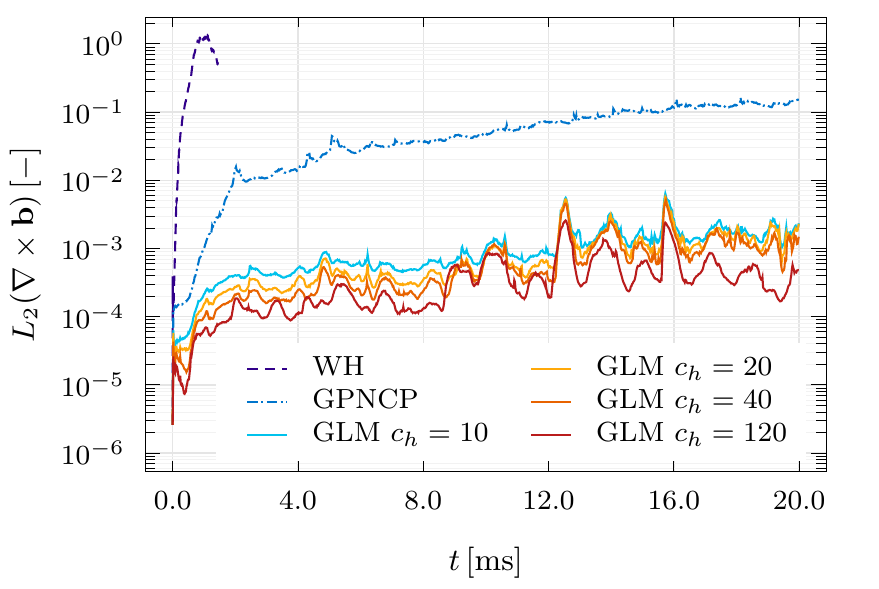}
    \caption{Time evolution of global dynamics and of curl errors for an oscillating elliptical water
    column. In the top row, left to right: the total kinetic energy $E_\up{k}$ and the $L_1$ norm of
    $\nabla \vec{b}$ over time. In the bottom row, the timeseries for the $L_1$ and $L_2$ norms of the
    curl constraint violation error. One can see that all the GLM curl cleaning simulation yield a
    constant oscillation period and kinetic energy is preserved quite well over more than
    $6\!\times\!10^5$ timestep, while the Godunov-Powell nonconservative formulation of the governing
    equations (GPNCP) shows signs of deterioration of the kinetic energy oscillation envelope. Moreover,
    curl errors decrease for increasing cleaning speed $c_h$, and in general GLM curl cleaning is
    effective in containing them, compared to the Godunov-Powell nonconservative formulation, and
    especially compared to the original weakly hyperbolic (WH) system, for which the computation fails
    rather early. The simulations employed a $\pnpm{5}{5}$ \aderdg scheme with \aderweno $\pnpm{0}{2}$
    subcell limiter on a uniform grid composed of $64^2$ elements.}
    \label{fig:curl-errors-ellipse}
\end{figure}

\subsubsection{Discussion of the results}
In this test, the \disgal scheme is supplemented by a third order \aderweno $\pnpm{0}{2}$
\finvol subcell limiter, which is an important ingredient for obtaining accurate results in this
test and for preserving the very complex smooth structures that arise in the flow (see
Figure~\ref{fig:droplet_oscillation_presentation}); in fact, even with the the extremely large value
we adopted for $\sigma$, the timescales associated with the theoretical oscillation period
$T_\up{p}^\up{a}$, given for small amplitude oscillations by the formula \cite{strutt1879, fyfe1988}
\begin{equation} \label{eq:periodsmalloscillation}
    T_\up{p}^\up{a} = 2\,\pi\,{\left[\sqrt{\frac{6\,\sigma}{(\rho_1 + \rho_2)\,R^3}}\right]}^{-1}, 
        \quad \text{with} \quad R = \frac{R_x + R_y}{2},
\end{equation}
are much larger than the timestep restriction for the numerical method in use and thus the task can
be regarded as a long-time integration problem. Specifically, to test the robustness of the
different formulations of the governing equations we evolved an oscillating elliptical droplet up to
a final time $t_\up{end} = 20.0\,\up{ms}$, which correspond to $605914$ timesteps for the GLM curl
cleaning simulation with $c_h = 10$, while the simulation with $c_h = 120$ required $649578$
timesteps, and the Godunov--Powell simulation reaches the final time in $637368$ timesteps. During
the run, the total kinetic energy of the droplet $E_\up{k}$ is tracked and subsequently employed to
measure the oscillation period of the droplet. Together with the kinetic energy, also the $L_1$ and
$L_2$ norms of the curl errors and of the gradient of the interface field $\nabla\vec{b}$ are
computed and stored. The norms of $\nabla\vec{b}$ can be taken as an indicator of the roughness of
the solution, which, with reference to Figure~\ref{fig:curl-errors-ellipse}, can in turn
indicate that the solution is developing spurious artefacts instead of maintaining its interface
field smooth; alternatively, since the interface field is sensitive to numerical diffusion, seeing
that $\nabla\vec{b}$ does not quickly decay, indicates that the scheme is not introducing excessive
artificial dissipation into the system. The timeseries for these integral quantities are shown in
Figure~\ref{fig:curl-errors-ellipse}, and again it appears that the GLM curl cleaning approach yields 
the best results: the Godunov--Powell simulation show some signs of deterioration of the kinetic
energy oscillation envelope and an increase in the average magnitude of the gradients of the
interface field. On the other side, the timeseries of kinetic energy obtained with GLM curl cleaning
does not show any signs of decay in the solution quality, and this is reflected in the fact that
violations of curl involutions are significantly lower with respect to the Godunov--Powell run.
Moreover, after fifteen oscillation periods, we cannot observe any decay in kinetic energy due to
numerical diffusion, which one would expect from lower order explicit methods for compressible
flow. Finally, it can be confirmed that the weakly hyperbolic formulation of the equations
is not well suited for solving time-dependent problems with high order Godunov-type schemes, as our
computations, employing a sixth order ADER-DG method on $64^2$ cells, blew up before a single full
oscillation period could be simulated.

\begin{figure}[!bp]
    \includegraphics[draft=false, scale=\figurescalefactorc]{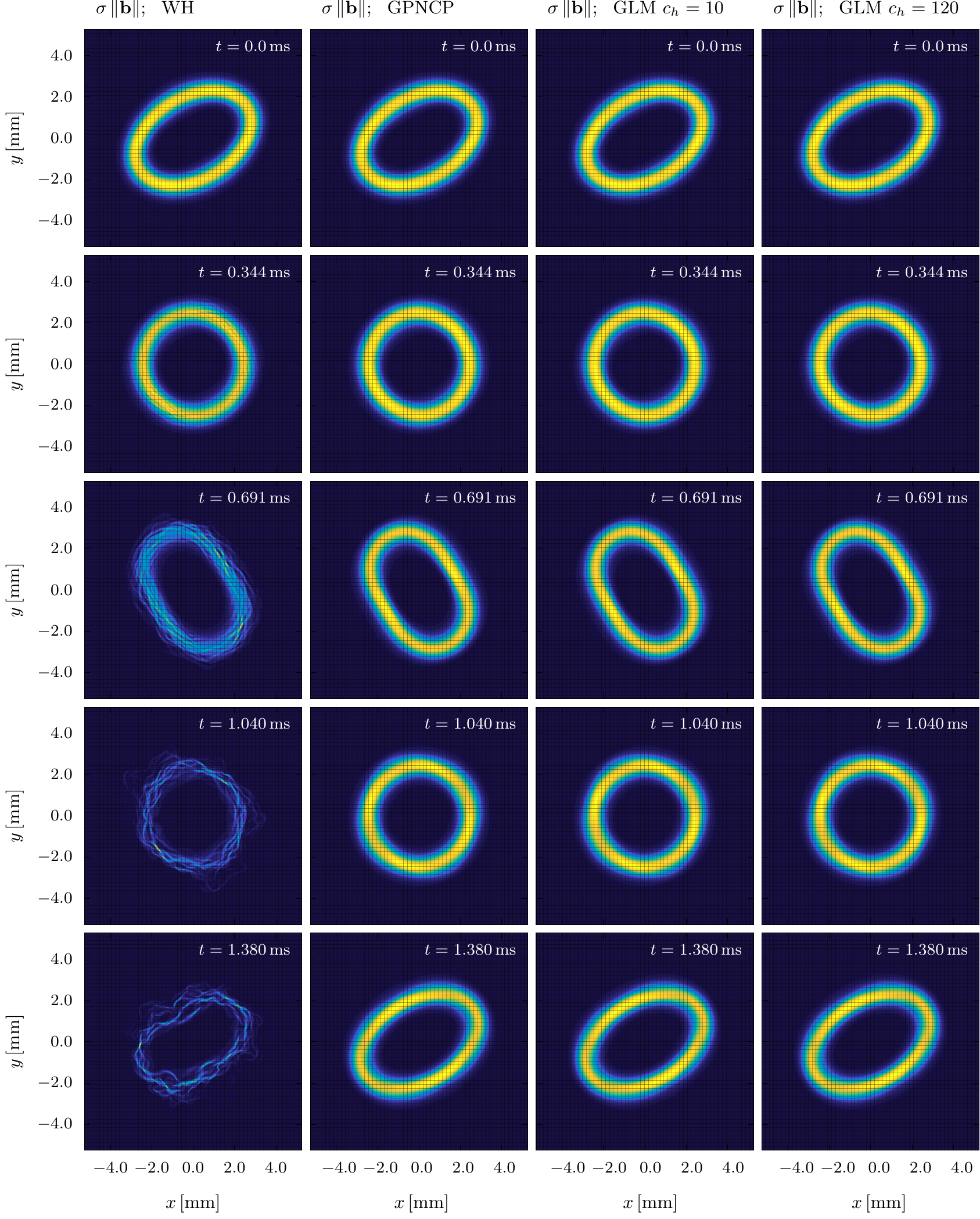}%
    \caption{Early stages of the oscillation of a two-dimensional elliptical droplet. The first column
    shows the quick deterioration of the interface field that is observed when solving the weakly
    hyperbolic formulation of the model (WH), the second shows that restoring hyperbolicity with the
    Godunov-Powell--type nonconservative products (GPNCP) prevents such ill behaviour, and the same is true
    for the GLM curl cleaning approach with different cleaning speeds $c_h$.}
    \label{fig:early-droplet}
\end{figure}

\begin{figure}[!bp]
    \includegraphics[draft=false, scale=\figurescalefactorc]{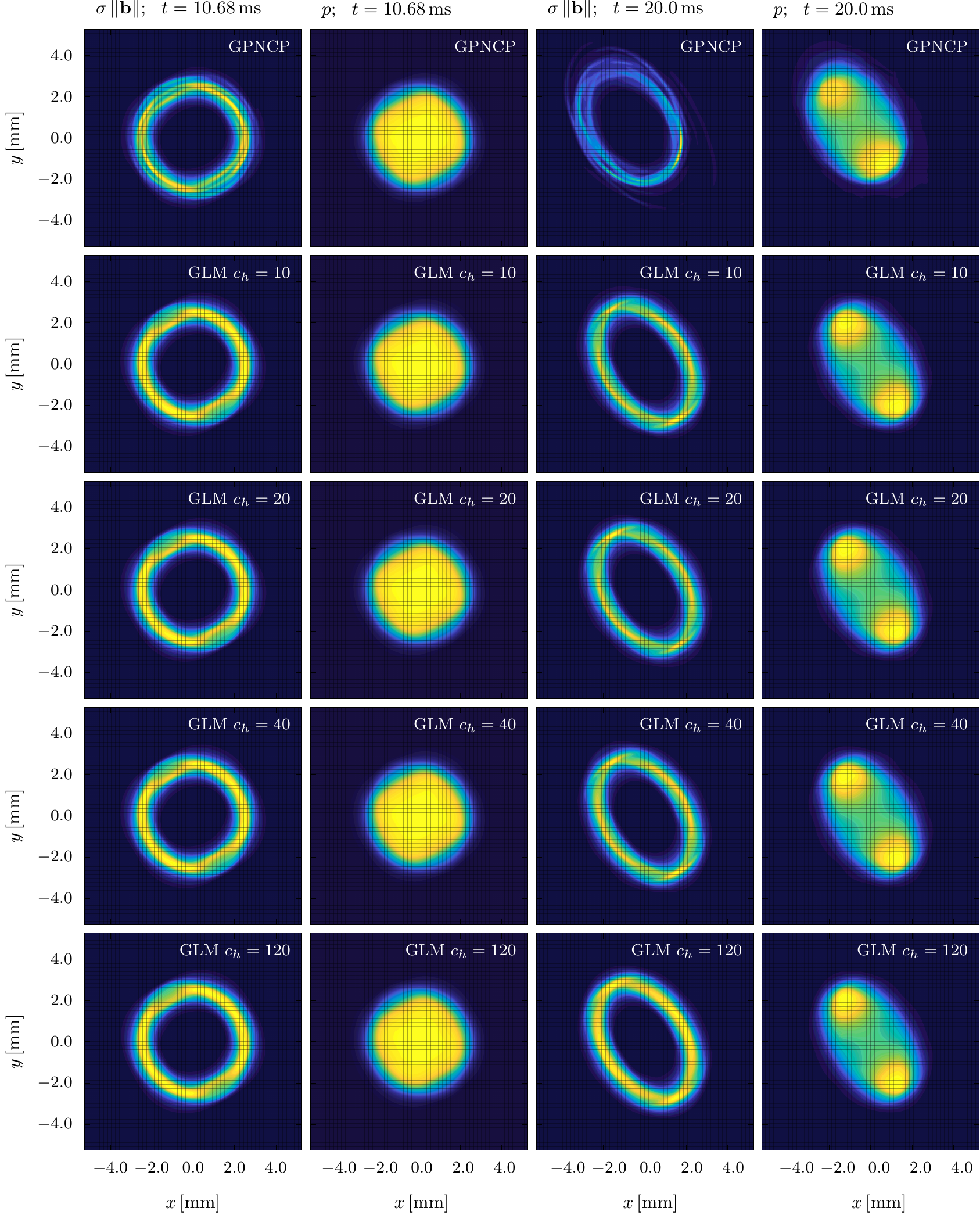}%
    \caption{Late stages of the oscillation of a two-dimensional elliptical droplet. In the first two
    columns we plot the surface energy $\sigma\,\norm{\vec{b}}$ and the pressure $p$ three quarters
    through the eighth oscillation cycle, and in the last two at about half of the fifteenth oscillation.
    The first row shows the results of the nonconservative Godunov--Powell-type model, with clear signs of
    deterioration of the solution, while in the other rows one can see the consistent results of GLM
    curl cleaning with several cleaning speeds $c_h$. }
    \label{fig:midlate-droplet}
\end{figure}

In Figure~\ref{fig:early-droplet}, we see that the simulation reproduces the expected global dynamic
behaviour in that the droplet can be observed achieving a circular shape at a quarter and at three
quarters of the oscillation period, when the maxima of kinetic energy are reached, while the minima
of kinetic energy, defining the half-period and end-period times, correspond to the maximum
elongation of the droplet in orthogonal directions. In these early stages of the simulation, we can
observe very little difference between the results from the GLM curl cleaning simulations with low
cleaning speed ($c_h = 10\,\up{m\,s^{-1}}$) or high cleaning speed ($c_h = 120\,\up{m\,s^{-1}}$), as
well as with respect to the solutions of the Godunov--Powell nonconservative variant of the model.
On the contrary, it is clear in the first column of Figure~\ref{fig:early-droplet}, that the
original weakly hyperbolic model does instead develop spurious filaments in the interface field
starting from the first oscillation period, leading to a very early breakdown of the computation.

In Figure~\ref{fig:midlate-droplet}, we compare the solutions, at two different time instants,
obtained from the Godunov--Powell-type formulation of the model \eqref{eq:model} and from the augmented
GLM curl cleaning system \eqref{eq:gavrilyuk.glm}, with several different values of cleaning speed
$c_h$. At both time instants, we can see only small differences among the simulations using GLM
curl cleaning, while it is clear in the first column of Figure~\ref{fig:early-droplet} that with the
nonconservative Godunov--Powell-type model some secondary subdivisions of the interface field are
starting to develop at time $t=10.68\,\up{ms}$, three quarters through the eighth oscillation cycle.
The effects on the pressure field are not yet visible at this stage, but one can see in the last two
columns of Figure~\ref{fig:midlate-droplet} that at the end time $t = 20\,\up{ms}$, about half of
the fifteenth oscillation cycle, these artefacts have further developed and have caused not only a
visible distortion of the pressure field, but also determined a bulk shift of the full droplet,
which is no longer centred at the origin of the computational domain, as already seen in the
three-dimensional droplet in Figures~\ref{fig:threedimensional-slices-2d} and
Figure~\ref{fig:threedimensional-slices-1d}.

Our numerical estimate of the oscillation period is obtained by solving the nonconservative 
Godunov--Powell-type 
formulation of the model with an ADER-DG $\pnpm{7}{7}$ scheme and \aderweno
$\pnpm{0}{2}$ subcell limiter on a mesh of $50^2$ cells; the interface thickness is set to
$\epsilon=2.5\times10^{-4}$. The deviation of the numerical estimate $T_\up{p} = 1.337\,\up{ms}$
with respect to the analytical prediction of $T_\up{p}^\up{a} = 1.3097\,\up{ms}$ is of $2.1\%$.
While not exact, the result can be considered satisfactory, since the difference can be attributed in
part to the fact that the reference formula \eqref{eq:periodsmalloscillation} was obtained under the
assumption of small amplitude oscillations, and is thus valid only for almost circular droplets. 
Another potential source of deviation from the predictions of linear theory is the diffuse interface
representation of the droplet, which distributes mass differently with respect to the ideal sharp
interface jump. In this regard, we computed the oscillation period also from another set of simulations,
namely those comparing the long-term behaviour of the Godunov--Powell and of the GLM curl cleaning
formulations of the model, employing an ADER-DG $\pnpm{5}{5}$ scheme and \aderweno $\pnpm{0}{2}$
subcell limiter on a mesh of $64^2$ elements, while the interface thickness is doubled with respect
to the previous run from $\epsilon=2.5\times10^{-4}$ to $\epsilon=5.0\times10^{-4}$. The estimated
period for this more diffuse droplet is
{$T_\up{p}^\textsc{gpncp} = 1.375\,\up{ms}$}
with the 
Godunov--Powell-type variant of the model,
while we computed {$T_\up{p}^\textsc{glm} = 1.377\,\up{ms}$} with the GLM curl cleaning formulation.
These estimates for the droplet with a thicker interface
correspond to a difference of
$2.8\%$ to $3.0\%$ with respect to the previous estimate $T_\up{p} = 1.337\,\up{ms}$, and deviate by
$5.0\%$ to $5.1\%$ from the small oscillations theory, despite the interface thickness being twice as large 
as the one used in the previous run.

\section{Summary and conclusions}
\label{sec:conclusions}
In this paper, we have presented two new \emph{strongly hyperbolic} reformulations of the originally 
\emph{weakly hyperbolic} two-phase flow model with surface tension \cite{Berry2008a, Schmidmayer2017}. 
Both reformulations heavily rely on the curl-free constraint \eqref{curl.b} of the interface vector 
field $\vec{b}$. The first of the two hyperbolic extensions was based on the theory of Symmetric
Hyperbolic and Thermodynamically Compatible (SHTC) systems \cite{God1961,Godunov1996,SHTC-GENERIC-CMAT},
and consists in modifying the momentum and energy equations by adding some \emph{symmetrising}
nonconservative terms, which are multiples of the curl involution constraint \eqref{curl.b} and 
thus, are formally zeros. 
The second reformulation of model \cite{Schmidmayer2017} was based on ideas developed in the 
context of numerical
Magnetohydrodynamics (MHD) and is more specifically based on the hyperbolic 
Generalized Lagrangian Multiplier (GLM) divergence-cleaning approach of Munz 
\etal\  \cite{MunzCleaning, Dedneretal}.  

We have then carefully compared the stability of the two reformulations, their accuracy, and consistency 
properties as well as showing the behaviour of the original weakly hyperbolic model.
The comparison seems to suggest that the weakly hyperbolic formulation is not suitable 
for a direct discretisation with general purpose explicit high order schemes, such as the 
ADER \disgal and ADER \finvol methods of this work, whereas the novel strongly hyperbolic reformulations 
can be shown to produce correct results, specifically for simulating the dynamics of oscillating droplets 
and they in general yield numerical results that are in agreement with analytical predictions. 

Of the two proposed models, the GLM approach uses hyperbolic constraint cleaning in order 
to minimise the violations of the curl involution constraint, obtaining a strongly hyperbolic and conservative model, while 
the first differs substantially from the weakly hyperbolic system \cite{Berry2008a, Schmidmayer2017}
in that it has a full set of eigenvectors, but no explicit enforcement of curl involution constraints
has been introduced. 
The results obtained with the GLM curl cleaning formulation are measurably better
than those of the first hyperbolic model here proposed, and this suggests that it is important to 
enforce such constraints in order to achieve reliable results in long time integration problems. 
In any case, it can be noted that for short times, or by increasing the resolution of the spatial 
discretisation, the nonconservative Godunov--Powell model can yield solutions that are comparable 
to those obtained using GLM curl cleaning (see for example Figures~\ref{fig:early-droplet} and 
\ref{fig:twodimensional-advection-3}). This is not the case for the weakly hyperbolic model, 
which does not seem to improve its behaviour significantly with mesh refinement, as apparent  
in Figure~\ref{fig:curl-errors-circle}

Another finding, that is completely independent of the model variant or of the numerical scheme in use, 
being based on an exact solution common to all three variants in consideration, is that the
pressure jump between the environment and region inside of a droplet is not particularly sensitive 
to the thickness of the smoothed interface, 
and in particular the error vanishes quadratically with respect to the width of the smoothing region.
This is also reflected in the bulk dynamical behaviour of oscillating droplets, in the sense that 
the oscillation period can be predicted with sufficient accuracy even with smoothed interfaces. 

Future developments will concern the extension of the mathematical model of this paper to the 
conservative multi-phase model of \cite{RomenskiTwoPhase2007,RomenskiTwoPhase2010} with 
pressure and velocity relaxation, as well as the development of exactly curl preserving schemes
on suitably staggered meshes, along the lines of exactly divergence-preserving schemes for the
Maxwell and MHD equations \cite{Yee66,DeVore,BalsaraSpicer1999,ADERdivB,BalsaraCED}, which need to be 
properly extended to curl-type involution constraints. We will also consider applications of the 
hyperbolic GLM curl cleaning approach to the first order hyperbolic reformulation 
of the nonlinear Schr\"odinger equation recently proposed in \cite{Dhaouadi2018}. 

\section*{Acknowledgments} 

The research presented in this paper has received financial support from the European Union's 
Horizon 2020 Research and Innovation Programme under the project \textit{ExaHyPE}, grant 
agreement number no. 671698 (call FET-HPC-1-2014).  

S.C. acknowledges the financial support received by the Deutsche Forschungsgemeinschaft (DFG) under 
the project \textit{Droplet Interaction Technologies (DROPIT)}, grant no. GRK 2160/1. 

I.P. acknowledges partial support by Agence Nationale de la Recherche
(FR) (grant ANR-11-LABX-0040-CIMI) under program ANR-11-IDEX-0002-02. 

M.D. also acknowledges the financial support received from the Italian Ministry of Education, University 
and Research (MIUR) in the frame of the Departments of Excellence Initiative 2018--2022 attributed to 
DICAM of the University of Trento (grant L. 232/2016) and in the frame of the PRIN 2017 project.
M.D. has also received funding from the University of Trento via the 
\textit{Strategic Initiative Modeling and Simulation}. 

\bibliographystyle{plain}
\bibliography{hst-bib}%
\end{document}